\documentclass[12pt,leqno,oneside]{amsart}
\usepackage{mathrsfs,dsfont}
\usepackage{amsmath,amstext,amsthm,amssymb,bbm}
\usepackage{charter}
\usepackage{typearea}
\usepackage{pdfsync}
\usepackage[backref=page]{hyperref}

\usepackage{listings}
\usepackage{enumitem}
\usepackage{comment}

\usepackage{color}

\usepackage{xfrac}

\usepackage[a4paper,width=16.2cm,top=2.5cm,bottom=2cm,
footskip=1cm]{geometry}

\newtheorem{Theorem}{Theorem}[section]
\newtheorem{Proposition}[Theorem]{Proposition}
\newtheorem{Lemma}[Theorem]{Lemma}
\newtheorem{Corollary}[Theorem]{Corollary}
\theoremstyle{definition}
\newtheorem{Definition}[Theorem]{Definition}
\newtheorem{Example}[Theorem]{Example}
\newtheorem{Remark}[Theorem]{Remark}

\numberwithin{equation}{section}

\newcommand{\dbar}{{\overline\partial}}

\DeclareMathOperator{\supp}{supp}

\DeclareMathOperator{\End}{End}
\DeclareMathOperator{\Vol}{Vol}
\DeclareMathOperator{\Div}{Div}
\DeclareMathOperator{\Bl}{Bl}

\newcommand{\field}[1]{\mathbb{#1}}

\newcommand{\R}{\field{R}}
\newcommand{\C}{\field{C}}
\newcommand{\N}{\field{N}}

\newcommand{\D}{\field{D}}

\newcommand{\cLL}{\mathcal{L}^{2}}
\newcommand{\E}{\mathbb{E}}

\newcommand{\V}{\mathrm{Var}}
\renewcommand{\P}{\mathbb{P}}

\begin{document}
\title{Gaussian 
holomorphic sections on noncompact complex manifolds}

\author{Alexander Drewitz, Bingxiao Liu and George Marinescu} 
\address{Universit{\"a}t zu K{\"o}ln,  Department Mathematik/Informatik,
    Weyertal 86-90,   50931 K{\"o}ln, Germany}
    \email{adrewitz@uni-koeln.de}
\email{bingxiao.liu@uni-koeln.de}
\email{gmarines@math.uni-koeln.de}
\thanks{The authors are partially supported by the 
DFG Priority Program 2265 `Random Geometric Systems' (Project-ID 422743078).}
\thanks{G.\ M.\ is partially supported by DFG funded projects SFB/TRR 191
(Project-ID 281071066-TRR 191), 
and the ANR-DFG project QuaSiDy (Project-ID 490843120).
}

\date{\today}
\begin{abstract}
We give two constructions of Gaussian-like random holomorphic 
sections of a Hermitian holomorphic line bundle $(L,h_{L})$ on a Hermitian 
complex manifold $(X,\Theta)$. In particular, we are  interested in 
the case where the space of $\cLL$-holomorphic sections 
$H^{0}_{(2)}(X,L)$ is infinite dimensional. 
We first provide a general construction
of Gaussian random holomorphic sections of $L$, which, if $\dim 
H^{0}_{(2)}(X,L)=\infty$, are almost never $\cLL$-integrable on $X$. 
The second construction combines the abstract Wiener space theory
with the Berezin-Toeplitz quantization and yields a random 
$\cLL$-holomorphic section. 
Furthermore, we study their random zeros in the context of semiclassical 
limits, including their equidistribution, large deviation estimates and hole 
probabilities.
\end{abstract}

\maketitle

\tableofcontents

\section{Introduction}\label{introduction}

Let $(X,J,\Theta)$ be a connected complex 
$n$-dimensional manifold without boundary, where 
$J$ denotes the complex structure and $\Theta$ is a Hermitian form. 
To $\Theta$ one can associate a $J$-invariant Riemannian metric 
$g^{TX}(\cdot,\cdot)=\Theta(\cdot,J\cdot)$. Let $L$ be a holomorphic line bundle over $X$, and let $h_{L}$ 
be a smooth Hermitian metric on $L$. We denote 
the corresponding Chern curvature form of $L$ by $R^L$.

In this paper we aim to study the zeros of certain
random holomorphic sections in $H^0(X,L)$, especially in the case of 
noncompact $X$. When $X$ is noncompact, we are mainly concerned with the 
subspace of $H^0(X,L)$ consisting of $\cLL$-integrable holomorphic 
sections, denoted by $H^0_{(2)}(X,L)$, which is a separable Hilbert 
space equipped with 
the $\cLL$-metric. We set $d:=\dim_{\C} H^0_{(2)}(X,L)$. Note that 
without any further assumptions on $X$ or $L$, the dimension $d$ can be $0$, a 
positive integer, or $\infty$. Our main interest is in 
the case of $d=\infty$, where some natural constructions of random 
sections with $d<\infty$ -- such as the Gaussian probability 
measure on $H^0_{(2)}(X,L)$ given by the $\cLL$-metric -- will fail in 
this case. To tackle this difficulty, we will provide two different approaches of 
constructing a random holomorphic section from the infinite 
dimensional $H^0_{(2)}(X,L)$ which both are natural as extensions of the finite 
dimensional case.

The first approach is a direct generalization of the study of random 
holomorphic functions on $\C^{n}$ to the context of complex geometry. The random holomorphic functions 
given by power series 
on $\C$ as well as the distribution of their  zeros (or other values) have been studied by Littlewood-Offord 
\cite{LOff45,LOff48}, Offord \cite{Off65,Off67,Off95}, and by 
Edelman-Kostlan \cite{EK95, EK96}, etc. Then for 
Gaussian random holomorphic functions, the results have further been extended by 
Sodin \cite{Sodin00},
Sodin-Tsirelson \cite{STr, STr05, STr06}, and then, on $\C^{n}$, by 
Zrebiec \cite{MR2369936}. In particular, the general Gaussian random 
holomorphic functions on the domains in $\C$ (also known under the 
name Gaussian analytic functions, GAFs) have been investigated vastly (cf.\ 
\cite{HKPV2009}) from 
probabilistic perspectives, serving as examples of the point processes on 
$\C$.

In fact, one can trace back to the work of Paley and Zygmund 
\cite{PZ1930} and Paley-Wiener-Zygmund \cite{PWZ33} for the construction of 
general random functions as well as their properties, including the study on the 
Rademacher series, random Fourier series, random Taylor series, etc. 
A general framework would be to construct random variables in a Banach 
or Hilbert space of functions, and we refer to the book of Kahane 
\cite{Kahane85} as well as the references therein for this purpose.

Now we explain our first approach in which we construct a Gaussian random section in 
terms of an 
orthonormal basis of $H^{0}_{(2)}(X,L)$, but its distribution, as 
a holomorphic section, is independent of the choice of such basis (by Proposition 
\ref{prop:uniqueness}). More concretely, if $\{S_{j}\}^{d}_{j=1}$ is an orthonormal basis of 
$H^{0}_{(2)}(X,L)$ with respect to the $\cLL$-metric, and if
$\{\eta_{j}\}_{j=1}^{d}$ denotes a sequence of independent and 
identically distributed (i.i.d.) standard
complex Gaussian variables, then we can define a random holomorphic 
section of $L$ via
\begin{equation}
	\psi^{S}_{\eta}:=\sum_{j=1}^{d}\eta_{j}S_{j},
	\label{eq:intro1.1}
\end{equation}
by using elementary properties of of the Bergman kernel 
associated with $H^{0}_{(2)}(X,L)$ (cf.\ Proposition \ref{prop:convergence}).
We will call $\psi^{S}_{\eta}$ a standard Gaussian random 
holomorphic section of $L$.

The above facts are nontrivial when $d=\infty$. In 
particular, in this case, it turns out that $\psi^{S}_{\eta}$ as constructed in 
\eqref{eq:intro1.1} is almost surely non-$\cLL$-integrable over $X$ 
(cf.\ Lemma \ref{lem:L2}). 
Such observation exhibits the abundance of holomorphic sections 
of certain holomorphic line bundles on a noncompact complex manifold. 
In the case of the Bargmann-Fock space on $\C^{n}$ (cf.\ Example 
\ref{ex:2.9}), $\psi^{S}_{\eta}$ is just a
Gaussian holomorphic function on $\C^{n}$ as mentioned before. If 
$d<\infty$ the above construction is equivalent to endowing 
$H^{0}_{(2)}(X,L)$ with the standard Gaussian probability measure associated 
to the $\cLL$ inner product.

From the above observations, a naturally ensuing and interesting 
question is how to randomize $\cLL$-holomorphic sections 
in a natural way, or equivalently, how to construct 
Gaussian probability measures on $H^{0}_{(2)}(X,L)$ \textit{in 
a geometric way}. Our second 
approach provides an answer to this question by 
combining the abstract Wiener space 
approach from probability theory with the Toeplitz operator machinery from 
geometric quantization.

One simple way to 
understand this approach would be as follows: fix a nonzero element 
$\mathbf{a}=(a_{j})_{j=1}^{d}\in\ell^{2}(\C);$ instead of 
\eqref{eq:intro1.1}, we 
define
\begin{equation}
	\psi^{S}_{\mathbf{a},\eta}:= \sum_{j=1}^{d}\eta_{j}a_{j}S_{j}.
	\label{eq:intro1.2}
\end{equation}
Then $\psi^{S}_{\mathbf{a},\eta}$ is almost surely 
$\cLL$-integrable holomorphic section of $L$ on $X$, which follows 
from the claim $\P(\sum_{j}|a_{j}|^{2}|\eta_{j}|^{2}<\infty)=1$. In our method, 
the (point) spectra of certain Toeplitz operators will play the role 
of the sequence
$\mathbf{a}=(a_{j})_{j=1}^{d}$, which have  significant applications in 
the theory of geometric quantization.

The approach to construct random sections 
\eqref{eq:intro1.2} is by considering an injective Hilbert-Schmidt Toeplitz operator 
$T_{f}$ on $H^{0}_{(2)}(X,L)$ associated with certain positive function 
$f$ on $X$ (for instance, a smooth positive function with compact 
support). This then defines a measurable norm
$\|T_{f}\cdot\|$ on 
$H^{0}_{(2)}(X,L)$ (cf.\ Definition \ref{def:4.1a}). 
As a consequence of the theory of abstract Wiener space by Gross 
\cite{Gross67}, we can construct in a unique way a Gaussian-like 
probability measure $\P_{f}$ on 
$H^{0}_{(2)}(X,L)$ associated with $T_{f}$. This way, the random $\cLL$-holomorphic section 
following the probability law $\P_{f}$ is exactly given as in 
\eqref{eq:intro1.2}, where each $a_{j}>0$ is an eigenvalue of $T_{f}$ and 
the orthonormal basis $\{S_{j}\}_{j=1}^{d}$ is such that
\begin{equation}
	T_{f}S_{j}=a_{j}S_{j}.
\end{equation}
For a brief introduction to Gross' abstract Wiener spaces 
we refer to \cite[Example 1.25]{Janson1997}.

On top of the constructions of random holomorphic sections outlined above, we 
aim to study the distributions of their zeros as 
$(1,1)$-currents on $X$ in the framework of semiclassical limits, i.e., 
considering the random holomorphic sections of the 
sequence of high tensor powers $(L^{p},h^{p}):=(L^{\otimes 
p},h_{L}^{\otimes p})$, $p\in\N$, of a given positive Hermitian line 
bundle $(L,h_{L})$. As $p\rightarrow \infty$, the number $h:=1/p$, 
playing the role of the Planck 
constant, tends to $0$.

For this purpose, we need to make further assumptions on $(X,J,\Theta)$ and $(L,h_{L})$, which will be made precise
later on. Then we consider the sequence of random sections $\psi^{S_{p}}_{\eta}$ 
constructed as in \eqref{eq:intro1.1} from the Hilbert spaces 
$H^{0}_{(2)}(X,L^{p})$, $p\in\N$. Set $d_{p}:=\dim_{\C} 
H^{0}_{(2)}(X,L^{p})\in\N\cup\{\infty\}$. As $p\rightarrow\infty$, 
the equidistribution of the normailzed zeros of $\psi^{S_{p}}_{\eta}$ is 
expected on $X$, where the limit is given by the first Chern form 
$c_{1}(L,h_{L})$. From this scope, we 
will 
extend in this paper the classical results on the random zeros for compact 
K\"{a}hler manifolds to the general noncompact setting.

The equidistribution of zeros of the $\mathrm{SU}(2)$-polynomials as their
degree tends to infinity was obtained by Bogomolny, Bohigas and 
Leboeuf in their paper \cite{BBL96}, where the inverse of the degree 
plays the role of Planck constant $h$. A relevant problem on the $2$-torus in the context of 
quantum chaotic dynamics was also studied in the paper of Nonnenmacher-Voros 
\cite{NoVo:98}. Then Shiffman and Zelditch 
\cite{ShZ99} extended further these results to the case of compact K\"{a}hler manifold 
equipped with a prequantum line bundle by considering the high tensor 
powers explained as above. They also discussed the 
equidistribution of the zeros of quantum ergodic eigensections. One 
key ingredient in their approach is the asymptotic expansion of the 
associated 
Bergman kernel (cf.\ \cite{Ti90}, \cite{Zel1998}, \cite{MM07} and the references 
therein). Dinh 
and Sibony \cite{DS06} introduced a different approach using \
ideas from the complex 
dynamics, which also gives an estimate for the speed of convergence of 
the distributions of random zeros. Subsequently, Dinh, Marinescu and Schmidt 
\cite{DMS}
extended such results to the noncompact setting, where they needed to assume 
$d_{p}=\mathcal{O}(p^{n})$ for $p\gg 0$. 
Along these lines, there are also plenty of 
generalizations to different geometric or 
probabilistic settings, cf.\ \cite{BCM, BL13, CM11,CMM,CMN15,CMN17,DMM}. We refer to 
the survey papers \cite{Zel2001} and \cite{BCHM} for more details and 
references on this topic.

In particular, for the compact K\"{a}hler 
manifold, Shiffman, Zelditch and Zrebiec \cite{SZZ} established the large 
deviation estimates for the random zeros of Gaussian holomorphic 
sections as the tensor power $p$ grows to infinity, and as a consequence, they 
obtained the expected exponential decay of the hole 
probabilities, which are the probabilities of that the Gaussian random holomorphic sections do not vanish 
on a given domain in $X$. In our previous paper \cite{DLM:21}, we 
generalized their results to the noncompact setting, especially the 
case of Riemannian surfaces with cusps, under the assumption $d_{p}=\mathcal{O}(p^{n})$ for $p\gg 0$. Then in this paper, we can 
finally complete the last piece of the puzzle such that the results also 
extend to the Gaussian random
holomorphic sections $\psi^{S_{p}}_{\eta}$ without assuming 
$d_{p}<\infty$.

Under this semiclassical 
setting, the use of the Toeplitz 
operators in our construction of the random $\cLL$-holomorphic 
sections becomes a natural way. The family of 
Toeplitz operators $T_{f,p}\in\mathrm{End}(H^{0}_{(2)}(X,L^{p}))$, 
$p\in\N$, are called Berezin-Toeplitz quantization of a given 
real smooth function $f:X\rightarrow 
\R$ which is also known as a classical observable in classical 
mechanics with phase space $(X,\Theta)$ (cf.\ \cite{BMS94}). Such operators are central 
object in the study of geometric quantization on K\"{a}hler or, in 
general, symplectic manifolds. For more details, we refer to the papers of Ma and 
Marinescu \cite{MM08,MM11,MM12} and their book \cite[Chapter 7]{MM07}.

Here we introduce a class of 
functions $f$ on $X$ such that $T_{f,p}$ is Hilbert-Schmidt for all $p\gg 
0$. Associated to a positive smooth function $f$ in this class, we construct canonically a sequence of probability spaces 
$(H^{0}_{(2)}(X,L^{p}),\P_{f,p})$, $p\gg 0$. Then we are concerned with the 
asymptotic behaviors of the zeros of random 
$\cLL$-holomorphic sections as $p\rightarrow\infty$.
Their limit as $(1,1)$-currents will be given by $c_{1}(L,h_{L})$ 
but only on the support of $f$ (with vanishing points of order at most 2). 
When we consider the random zeros inside a small ball of the 
Planck scale (i.e., radius $\sim \frac{1}{\sqrt{p}}$), the precise 
values of $f$ can make 
differences on the fluctuations of random zeros. A further 
interesting question would be to describe the asymptotic distribution of random 
zeros outside the support of $f$.

Our approach to the above results 
relies on the asymptotic expansion of the on-diagonal Schwartz kernel 
of the operator
$T^{2}_{f,p}=T_{f,p}\circ T_{f,p}$, 
as $p\rightarrow \infty$, whose 
first several terms are computed explicitly in \cite[Chapter 
7]{MM07} and in \cite{MM12} (for the K\"{a}hler case with a prequantum line 
bundle). Note that in \cite{MM12}, the formulae are stated for a compact 
K\"{a}hler manifold, but their computations are actually local, so 
that the formulae extend to certain cases of noncompact manifolds. In 
particular, we can apply them to the case considered in \cite[Section 7.5]{MM07} and the 
case of bounded geometry discussed in \cite{MM15} and \cite{Finski22a}. 

In the next four sections we provide the setting and formulate our main results.

\subsection{Zeros of Gaussian random holomorphic sections}\label{s1.1a} 
Let us start with a Hermitian holomorphic line bundle $(L,h_{L})$ on a 
(paracompact) complex 
manifold $(X,J,\Theta)$ with arbitrary 
$d=\dim_{\C} H^0_{(2)}(X,L)\geq1.$

For $s\in H^{0}(X,L)\setminus\{0\}$, let $Z(s)$ denote the set of zeros of $s$, which is 
a purely $1$-codimensional analytic subset of $X$. The divisor 
$\Div(s)$ of $s$ is then defined as the formal sum
\begin{equation}
	\Div(s)=\sum_{V \subset Z(s)} \mathrm{ord}_{V}(s) V,
	\label{eq:divisor}
\end{equation}
where $V$ runs over all the irreducible analytic hypersurfaces 
contained in $Z(s)$, and $\mathrm{ord}_{V}(s)\in \N^{+}$ denotes the 
vanishing order of $s$ along $V$. For any analytic hypersurface $V\subset X$,
we deonte by $[V]$ the current of integration on $V$, defined 
by $\varphi\mapsto\int_V\varphi$, where $\varphi$ runs
in space $\Omega^{(n-1,n-1)}_{0}(X)$ 
of $(n-1,n-1)$-form with compact support in $X$. The current
of integration (with multiplicities) on the divisor $\Div(s)$
is defined by 
\begin{equation}
[\Div(s)]=\sum_{V \subset Z(s)} \mathrm{ord}_{V}(s) [V],
	\label{eq:divisorint}
\end{equation}

Our first result concerns the expectation of the currents of integration
on the zero-divisors of the Gaussian random holomorphic section $\psi^{S}_{\eta}$ defined 
in \eqref{eq:intro1.1}, as a current on $X$, i.e.\ of the 
random $(1,1)$-current $[\Div(\psi^{S}_{\eta})]$. For any 
test form $\varphi\in\Omega^{(n-1,n-1)}_{0}(X)$, the random variable 
$\langle [\Div(\psi^{S}_{\eta})],\varphi\rangle$ is measurable 
(cf.\ \cite[proof of 
Proposition 4.2]{CM11}). If the random 
variable $\langle 
[\Div(\psi^{S}_{\eta})],\varphi\rangle$ is integrable 
for any test form $\varphi$, then the linear map 
$$\varphi\mapsto \E\left[\langle 
[\Div(\psi^{S}_{\eta})],\varphi\rangle\right]\,,\quad
\varphi\in\Omega^{(n-1,n-1)}_{0}(X),$$ defines a 
$(1,1)$-current on $X$, which is called the expectation of 
$[\Div(\psi^{S}_{\eta})]$ and denoted by 
$\E[[\Div(\psi^{S}_{\eta})]]$.

Next we define the \emph{Fubini-Study current}
$\gamma(L,h_L)$ on $X$.  Let 
\begin{equation}\label{e:BP}
P:\mathcal{L}^2(X,L)\to H_{(2)}^0(X,L)
\end{equation}
be the $\cLL$-orthogonal projection, called the \emph{Bergman projection}.
It has a smooth Schwartz kernel $P(x,y)$, called the \emph{Bergman 
kernel}, cf. Subsection \ref{S:HLB-BK}. 
The Bergman kernel function $X\ni x\mapsto P(x,x)$ is a 
non-negative smooth function on $X$, and 
the function $\log P(x,x)$ is locally integrable on $X$.
We set 
\begin{equation}
\gamma(L,h_{L})=c_{1}(L,h_{L})+\frac{\sqrt{-1}}{2\pi}
\partial\bar{\partial}\log P(x,x),
	\label{eq:intro2.23}
\end{equation}
where $c_{1}(L,h_{L})$ is the Chern form of $(L,h_{L})$.

\begin{Theorem}\label{thm:expectation}
	Assume that $d\geq 1$. Then the expectation of the random variable
	$[\Div(\psi^{S}_{\eta})]$ exists as a $(1,1)$-current on $X$. Moreover, we have 
	\begin{equation}
		\E[[\Div(\psi^{S}_{\eta})]]=\gamma(L,h_{L})
		\label{eq:2.26}
	\end{equation}
	as an
	identity of $(1,1)$-currents.
\end{Theorem}

In the case $d<\infty$, \eqref{eq:2.26} was already 
known for line bundles with empty base locus 
(cf.\ \cite[Lemma 3.1]{ShZ99}) and in several situations
when the metric $h_L$ or the base $X$ are singular 
(see e.g.\ \cite[Proposition 4.2]{CM11}, \cite[Theorem 1.4]{CM13}). 
When $d=\infty$ analogues of this result 
are known in the context of random holomorphic 
functions on $\C^{m};$ for instance, Edelman and Kostlan \cite[Sections 7 \& 
8]{EK95} studied the expectations 
of complex zeros of random power series 
(in their paper, they mainly aimed to study the distribution of real 
zeros). Other interesting examples from complex geometry, where our 
Theorem \ref{thm:expectation} applies, are given in Subsection 
\ref{ss:GeoExa}.
\subsection{High tensor powers of $L$: equidistribution and large 
deviations}\label{s1.2a}
We are interested in the semiclassical limit of the zeros of 
the Gaussian holomorphic sections when we replace $L$ by its high 
tensor powers. For this purpose,  we need to make further assumptions on 
$(X,J,\Theta)$ and on $(L,h_{L})$ as follows. 
We assume that Riemannian metric $g^{TX}$ is complete and
there exist $C,C_{0},\varepsilon>0$ such that on $X$,
\begin{equation}
\sqrt{-1}R^{L}\geq \varepsilon \Theta,\;\;
\sqrt{-1}R^{\det}\geq -C_{0} \Theta,\;\;|\partial \Theta|_{g^{TX}}\leq C\,,
	\label{eq:intro3.0.1}
\end{equation}
where $R^{\det}$ be the curvature of the holomorphic connection 
$\nabla^{\det}$ on $K_X^*=\det(T^{(1,0)}(X))$.
 
In this case, by \cite[Chapter 
6]{MM07}, the on-diagonal 
Bergman kernels 
$P_{p}(x,x)$ have an asymptotic expansion in the tensor power $p$, which 
is uniform on any given compact subset of $X$. As a consequence, we 
have the convergence of currents 
\begin{equation}
	\frac{1}{p}\gamma(L^{p},h^{p})\rightarrow c_{1}(L,h_{L}) \quad \text{ as }p\rightarrow \infty.
	\label{eq:intro1.3}
\end{equation}

In the following, we denote by $\psi^{S_{p}}_{\eta}$ the Gaussian random 
holomorphic section (as in \eqref{eq:intro1.1}) constructed from an orthonormal basis 
$S_{p}=\{S^{p}_{j}\}_{j=1}^{d_{p}}$ of $H^{0}_{(2)}(X,L^{p})$.
As is natural, before formulating our concentration estimates, we begin with stating findings for the limit of the expectations 
$\E[[\Div(\psi^{S_{p}}_{\eta})]]$. While the results are novel in our specific setting and formulated precisely in Theorems \ref{thm:asym} 
and \ref{thm:2.15} below, we roughly speaking prove the following: 
\begin{itemize}
	\item as $p\rightarrow \infty$, 
	$\frac{1}{p}\E[[\Div(\psi^{S_{p}}_{\eta})]]\rightarrow 
	c_{1}(L,h_{L})$; 
	\item for each $\varphi\in 
	\Omega^{(n-1,n-1)}_{0}(X)$, we have that
	\begin{equation}
		\P\Big(\lim_{p\rightarrow \infty} \frac{1}{p}\langle
		[\Div(\psi^{S_{p}}_{\eta})],\varphi\rangle= \langle 
		c_{1}(L,h_{L}),\varphi\rangle\Big)=1.
		\label{eq:intro2.43}
	\end{equation}
\end{itemize}
It is clear that the first point is a consequence of Theorem 
\ref{thm:expectation} in combination with \eqref{eq:intro1.3}. The almost sure 
convergence in the second point is deduced by means of the Bergman kernel. 

With these equidistribution results on the random zeros at 
our disposal, a natural next step is to investigate the speed 
of convergence in terms of large deviation estimates as in \cite{SZZ} and 
\cite{DLM:21}, but with the possibility $d_{p}=\infty$.

\begin{Theorem}\label{thm:6.2}
We assume that Riemannian metric $g^{TX}$ is complete and 
\eqref{eq:intro3.0.1} holds.
	If $U$ is a relatively compact open subset of $X$, then for any 
	$\delta >0$ and $\varphi\in \Omega^{(n-1,n-1)}_{0}(U)$, there 
	exists a constant $c=c(U,\delta,\varphi)>0$ such that for $p\in\N$,
	we have
	\begin{equation}
		\P\Big( \ 
		\Big|\Big\langle\frac{1}{p}[\Div(\psi^{S_{p}}_{\eta})] 
		-  c_{1}(L,h),\varphi\Big\rangle\Big|>\delta\  \Big)\leq 
		e^{-c\,p^{n+1}}.
		\label{eq:6.0.7}
	\end{equation}
\end{Theorem}
Another 
natural question is then the validity Central Limit Theorem for the 
distribution of zeros of $\psi^{S_{p}}_{\eta}$ as $p\rightarrow\infty$, which will be 
touched upon in Remark \ref{rm:asympnormality}. 

Since $c_{1}(L,h_{L})$ is positive, $\frac{1}{n!}c_{1}(L,h_{L})^{n}$ also 
defines a positive volume element on $X$. If $U\subset X$ is open, set
\begin{equation}
\Vol^{L}_{2n}(U)=\int_{U} \frac{1}{n!}{c_{1}(L,h_{L})^{n}}.
	\label{eq:1.6.8DLM}
\end{equation}

For $s_{p}\in H^{0}(X,L^{p})\setminus\{0\}$
we define the $(2n-2)$-dimensional volume with 
respect to $c_{1}(L,h_{L})$ of the divisor $\Div(s_{p})$ 
(cf.\ \eqref{eq:divisor})
in an open subset $U\subset X$ as follows:
\begin{equation}
	\Vol^{L}_{2n-2}\big(\Div(s_{p})\cap U\big)=\sum_{V\subset 
	Z(s_{p})}\mathrm{ord}_{V}(s_{p})\int_{V\cap U} 
	\frac{c_{1}(L,h_{L})^{n-1}}{(n-1)!}\,\cdot
	\label{eq:6.1.9}
\end{equation}
If we use this volume to measure the size of the zeros of $s_{p}$ in $U$, 
then Theorem \ref{thm:6.2} leads to the following result.
\begin{Theorem}\label{thm:6.3}
We assume that Riemannian metric $g^{TX}$ is complete and 
\eqref{eq:intro3.0.1} holds.
	If $U$ is a nonempty relatively compact open subset of $X$ such that 
	$\partial U$ has zero measure in $X$, 
	then for any $\delta >0$, there 
	exists a constant $c_{U,\delta}>0$ such that for $p$ 
	large enough,
	we have
	\begin{equation}
		\P\Big( \ 
		\Big|\frac{1}{p}\Vol^{L}_{2n-2}(\Div(\psi^{S_{p}}_{\eta})\cap U) 
		- n\Vol^{L}_{2n}(U) \Big|>\delta \ \Big)\leq 
		e^{-c_{U,\delta} p^{n+1}}.
		\label{eq:6.1.10}
	\end{equation}
	In addition, there exists a 
	constant $C_{U}>0$ such that for $p> 0$,
	\begin{equation}
		\P\big(\Div(\psi^{S_{p}}_{\eta})\cap U=\varnothing\big)
		\leq e^{-C_{U}p^{n+1}}\,.
		\label{eq:1.6.14paris}
	\end{equation}
\end{Theorem}

The proofs of the above two theorems will be provided in Subsection 
\ref{subsection3.2}. One essential ingredient for these proofs is 
Proposition \ref{prop:1.3.3}, for which we need a more refined 
investigation of the local sup-norms of holomorphic sections on $X$ 
(cf.\ Subsection \ref{subsection3.3}).

The probability in \eqref{eq:1.6.14paris} is referred to as \textit{hole 
probability} of the random section $\psi^{S_{p}}_{\eta}$ on the subset $U$. This 
estimate then provides us with an upper bound for the hole 
probabilities for $p> 0$. 
In \cite[Theorem 1.4]{SZZ} and \cite[Proposition 1.7]{DLM:21}, under additional assumptions on $U$, a 
lower bound of the form $e^{-C'_{U}p^{n+1}}$ for the hole probabilities was 
proved. In general though, such a lower bound remains unclear in the case 
$d_{p}=\infty$.

In the case of the Bargmann-Fock space, for the standard Gaussian 
random holomorphic function on $\C^{n}$ (cf.\ \eqref{eq:2.36}),  
the two-sided bound on the hole probabilities when 
$U=\mathbb{B}(0,r)$ as $r\rightarrow \infty$ was proved by 
Sodin-Tsirelson (for $\C$, \cite[Theorem 1]{STr05}) and by Zrebiec (for 
$\C^{n}$, \cite[Theorem 1.2]{MR2369936}). In Subsection 
\ref{ss3.4ss}, we will explain how to recover their results from our general
results being specialized to the scaled Bargmann-Fock spaces.

\subsection{Random $\mathcal{L}^{2}$-holomorphic sections and 
Toeplitz operators}\label{s1.3a}

In the setting of Section \ref{s1.1a} we introduce 
for a bounded function $f$ on $X$ the associated Toeplitz operator 
$T_{f}$ defined by 
$T_{f}: H^{0}_{(2)}(X,L)\ni S\mapsto P(fS)\in H^{0}_{(2)}(X,L)$, where $P$ is the Bergman projection
\eqref{e:BP} (see Definition \ref{def:4.4} for further details).

If 
$f$ is smooth and also satisfies
	\begin{equation}
		\int_{X}|f(x)|P(x,x)\mathrm{dV}(x)<\infty,
		\label{eq:intro5.8paris}
	\end{equation}
	then the operator $T_{f}$ is Hilbert-Schmidt (cf.\ Proposition 
	\ref{prop:4.7a}). If in addition $f$ is a real nonnegative function 
	(which shall not be identical zero), then $T_{f}$ is injective. 
	
	For such  nonnegative smooth function $f$, we get a Hilbert 
	metric 
	$\langle T_{f}\,\cdot,T_{f}\,\cdot\rangle_{\cLL(X,L)}$ on 
	$H^{0}_{(2)}(X,L)$, which is a measurable norm in the sense of 
	Gross (cf.\ \cite{Gross67}). Let $\mathcal{B}_{f}(X,L)$ be the Hilbert 
	space given as the completion of $H^{0}_{(2)}(X,L)$ under this 
	measurable norm. The theory of abstract Wiener spaces implies that for $f$ as above given, there exists a unique Gaussian probability measure 
$\mathcal{P}_{f}$ on $\mathcal{B}_{f}(X,L)$ such that it extends the 
Gaussian probability measure on any finite dimensional subspace of 
$\mathrm{Im}(T_{f})$ associated with the standard $\cLL$-metric. 

The injective linear operator $T_{f}$ extends to an isometry of 
Hilbert spaces
\begin{equation}
	\widehat{T}_{f}:\big(\mathcal{B}_{f}(X,L),\|T_{f}\cdot\|\big)\rightarrow \big(H^{0}_{(2)}(X,L),\|\cdot\|_{\cLL(X,L)}\big).
	\label{eq:intro1.20a}
\end{equation}
After taking the 
pushforward of $\mathcal{P}_{f}$ by 
$\widehat{T}_{f}$, we obtain a Gaussian probability measure $\P_{f}$ on 
$H^{0}_{(2)}(X,L)$.

The on-diagonal restriction $T^{2}_{f}(x,x)$ of the Schwartz 
kernel of $T^{2}_{f}=T_{f}\circ T_{f}$ is locally integrable 
on $X$ (cf.\ Lemma \ref{lm:4.13a}). As an analog of \eqref{eq:intro2.23}, we define a closed 
positive $(1,1)$-current on $X$ as 
\begin{equation}
	\gamma_{f}(L,h_{L})=c_{1}(L,h^{L})+\frac{\sqrt{-1}}{2\pi}\partial\bar{\partial}\log T^{2}_{f}(x,x).
	\label{eq:intro4.37ss}
\end{equation}
In Subsection \ref{ss4.4a}, we prove the following result for the expectation of the random zeros 
of $\cLL$-holomorphic section.

\begin{Theorem}\label{thm:4.15a}
	Denote by $s$ the identity on $(H^{0}_{(2)}(X,L),\P_{f})$ and consider
	the random variable $[\Div(s)]$ taking values in the space of  $(1,1)$-currents on $X$.
 Then
	\begin{equation}
		\mathbb{E}^{\P_{f}}\left[[\Div(s)]\right]=\gamma_{f}(L,h_{L}).
		\label{eq:3.1.1}
	\end{equation}
\end{Theorem}

\begin{Remark}\label{rk:4.12b}
	During our writing of this paper, we became aware of the work of 
	Ancona and Le Floch \cite{AF:2022} on random sections under the Toeplitz 
	operator $T_{f}$ for the case of compact K\"{a}hler $X$. For compact $X$ 
one has $d<\infty$ (we assume that $d>0$), and in this case the random 
	section $s$ in 
	$H^{0}(X,L)$ with the probability measure $\P_{f}$ defined above 
	has the same distribution as the random section $T_{f}s'$ 
	considered by Ancona and Le Floch, where $s'$ is the 
	random section in $H^{0}(X,L)$ with the standard Gaussian 
	probability measure given by the $\cLL$-metric. 
\end{Remark}
\subsection{High tensor powers of $L$: equidistribution on the 
support of $f$}\label{s1.4a}
To consider the semiclassical limit in the noncompact setting, we 
need to make the same assumptions as in Subsection \ref{s1.2a}. For 
simplicity, in this subsection we only consider a nontrivial nonnegative smooth function $f$ 
on $X$ with compact support. Note that our results hold for a general 
class of nonnegative smooth functions $f$ that are not required to 
have compact support (cf.\ Subsections \ref{ss5.1a} and \ref{ss:5.2b}).

Since $f$ has compact support, condition
\eqref{eq:intro5.8paris} is satisfied for the line bundle $L^{p}$
for each $p$. This way, we 
can construct a sequence of probability spaces 
$(H^{0}_{(2)}(X,L^{p}),\P_{f,p})$ using the corresponding Toeplitz 
operator $T_{f,p}$. We denote by $\mathbf{S}_{f,p}$ be the identity map on the 
canonical probability space $(H^{0}_{(2)}(X,L^{p}),\P_{f,p})$.

In Theorems \ref{thm:5.4a} and \ref{thm:5.5b}, we prove the general 
version of the following 
results.
\begin{Theorem}\label{thm:1.6a}
(1) Let $U$ be an open subset of $X$ such that 
$f>0$ on $U.$ Then, as 
$p\rightarrow\infty$, we have the weak convergence of $(1,1)$-currents
on $U$,
\begin{equation}
\frac{1}{p}\E^{\P_{f,p}}[[\Div(\mathbf{S}_{f,p})]|_{U}]
\rightarrow c_{1}(L,h_{L})|_{U}\,.
\label{eq:intro5.14a}
\end{equation}
	 

\noindent
(2) Moreover, for any $\varphi\in 
\Omega^{(n-1,n-1)}_{0}(\overline{U})$, we have
\begin{equation}
\P\left(\lim_{p\rightarrow 
\infty}\frac{1}{p}\langle[\Div(\mathbf{S}_{f,p})],
\varphi\rangle=\langle c_{1}(L,h_{L}),\varphi\rangle\right)=1.
	\label{eq:intro5.20a}
\end{equation}
\end{Theorem}

In fact, $f$ might vanish on some points in $\supp f$. Since $f\geq 
0$, then the smallest vanishing order of $f$ at a vanishing point is 
$2$. If we assume further that $\Theta$ is K\"{a}hler and that 
$(L,h_{L})$ is a prequantum line bundle (i.e., 
$c_{1}(L,h_{L})=\Theta$), then we can allow such kind of vanishing 
points in the set $U$ in Theorem \ref{thm:1.6a}.

\begin{Theorem}\label{thm:intro5.8a}
	For $f$ as above, assume furthermore the prequantum line bundle condition for $(L, h_{L})$. Let $U$ be an open subset of 
	$\supp f$ be such that $f$ only vanishes up to order $2$ 
in $U$ with nonzero $\Delta f$ at the vanishing points. Then as $p\rightarrow \infty$,
\begin{itemize}
	\item we have the weak convergence of $(1,1)$-currents on $U$
	\begin{equation}
		\frac{1}{p}\E^{\P_{f,p}}[[\Div(\mathbf{S}_{f,p})]|_{U}]\rightarrow c_{1}(L,h_{L})|_{U}.
	\end{equation}
	\item for any $\varphi\in 
\Omega^{(n-1,n-1)}_{0}(\overline{U})$, we have
\begin{equation}
	\P\left(\lim_{p\rightarrow 
	\infty}\frac{1}{p}\langle[\Div(\mathbf{S}_{f,p})],\varphi\rangle=\langle c_{1}(L,h_{L}),\varphi\rangle\right)=1.
\end{equation}
\end{itemize}
\end{Theorem}
A general version of the above theorem is provided in Theorem 
\ref{thm:5.8a}.

One important ingredient in the proofs of the above results is the 
following identity from Theorem \ref{thm:4.15a},
\begin{equation}
	\E^{\P_{f,p}}[[\Div(\mathbf{S}_{p,f})]]-pc_{1}(L,h_{L})=\frac{\sqrt{-1}}{2\pi}\partial\overline{\partial}\log(T^{2}_{f,p}(x,x)).
\end{equation}
Then considering the zeros in the small geodesic ball 
$B(x,R/\sqrt{p})$ centered at $x$ via pairing with a test form $\varphi\in 
	\Omega^{(n-1,n-1)}_{0}(X)$, in Subsection \ref{ss5.3ab}, our 
	computations (especially by Theorem \ref{thm:5.14a}) show that
\begin{equation}
	\begin{split}
		\Big\langle\E^{\P_{f,p}}[[\Div(\mathbf{S}_{f,p})]]-pc_{1}(L,h_{L}),\chi_{B\big(x,\frac{R}{\sqrt{p}}\big)}\varphi\Big\rangle
		=\begin{cases}
\mathcal{O}(p^{-n}), \!\!\!\!&\text{if $f(x)>0$;}\\
\mathcal{O}(p^{-n+1}),\!\!\!\!&\text{if $f(x)=0$, $\Delta f(x)<0$,}
\end{cases}
	\end{split}
	\label{eq:5.64a}
\end{equation}
where $\chi_{B(x,R/\sqrt{p})}$ is the indicator function of the set 
$B(x,R/\sqrt{p})$, and we have the explicit formulae for the 
coefficients of $p^{-n}$ and of $p^{-n+1}$ in the above estimates. The different powers in \eqref{eq:5.64a} show 
that, in the Planck scale, our random zeros can have higher fluctuations near a 
vanishing point of $f$ of order $2$ than near the nonvanishing points.

At last, in Subsection \ref{ss5.5}, we consider a not necessarily 
nonnegative real smooth function $f$ satisfying \eqref{eq:intro5.8paris} for $L^{p}$, $p\gg 
0$. In this case, $T_{f,p}$ 
might not be injective, and with suitable conditions on the vanishing 
points of $f$, we can still extend Theorem \ref{thm:intro5.8a} to this 
case.

The next four sections of this paper correspond exactly to the above four 
subsections describing the main results: the first two sections deal with
Gaussian random holomorphic sections, and the last two sections deal with
random $\cLL$-holomorphic sections using the Toeplitz operators.
\bigskip

{\bf Acknowledgment:} We gratefully acknowledge support of 
DFG Priority Program 2265 \lq Random Geometric Systems\rq. 
The authors thank Prof. Xiaonan Ma for useful discussions.

\section{Gaussian random holomorphic sections}\label{section2}
In this section, we define the Gaussian random holomorphic section of $L$ and 
study its zeros as a $(1,1)$-current on $X$.

While some results proved in this section are not new in the 
special case of random functions or power series, to the best of our 
knowledge, we were not able to locate these results for holomorphic 
sections with $d=\infty$ in the litterature. 

\subsection{Holomorphic line bundles and Bergman kernels}\label{S:HLB-BK}
Let $(X,J,\Theta)$ be a connected  
$n$-dimensional manifold (without boundary) where 
$J$ denotes the complex structure and $\Theta$ is a Hermitian form. 
To $\Theta$ one can associate a $J$-invariant Riemannian metric 
$g^{TX}(\cdot,\cdot)=\Theta(\cdot,J\cdot)$. Let $L$ be a holomorphic line bundle over $X$, and let $h_{L}$ 
be a smooth Hermitian metric on $L$. We denote 
the corresponding Chern curvature form of $L$ by $R^L$, and the first 
Chern form of $(L, h_{L})$ is denoted by
\begin{equation}
	c_{1}(L, h_{L})=\frac{\sqrt{-1}}{2\pi}R^{L}.
\end{equation}

Let $\mathscr{C}_0^{\infty}(X,L)$ denote the space of 
compactly supported smooth sections of $L$ on $X$. Associated with the metrics 
$g^{TX}$ and $h_{L}$, we define the 
$\mathcal{L}^{2}$-inner product as follows, for $s_{1}, s_{2}\in 
\mathscr{C}_0^{\infty}(X,L)$,
\begin{equation}\label{norm}
\langle s_1,s_2\rangle_{\mathcal{L}^{2}(X,L)}:=\int_X\langle 
s_1(x),s_2(x)\rangle_{h_L} \mathrm{dV}(x),
\end{equation}
where $\mathrm{dV}=\frac{1}{n!}\Theta^n$ is 
the volume form induced by $\Theta.$ We also let $\mathcal{L}^2(X,L)$ 
be the separable Hilbert space obtained by completing  $\mathscr{C}_0^{\infty}(X,L)$
with respect to the norm $\|\cdot\|_{\mathcal{L}^{2}(X,L)}$ induced by 
\eqref{norm}. 
Let $H^{0}(X,L)$ denote the vector space of holomorphic sections of 
$L$ over $X$. Set
\begin{equation}
	H_{(2)}^0(X,L):= \mathcal{L}^2(X,L) \cap H^0(X,L).
\end{equation}
It follows from the Cauchy estimates for holomorphic functions 
that for every compact set $K\subset X$ there exists $C_K>0$ such that
\begin{equation}
	\sup_{x\in K}|s(x)|\leq C_K\|s\|_{\mathcal{L}^2(X,L)}\ \ \ 
	\text{for}\ s\in H_{(2)}^0(X,L),
	\label{eq:1.3}
\end{equation}  
which in turn implies that $H_{(2)}^0(X,L)$ is a closed subspace of 
$\mathcal{L}^2(X,L)$. Moreover,  
$H_{(2)}^0(X,L)$ is a separable Hilbert space with induced 
$\mathcal{L}^{2}$-metric (cf.\ \cite[p.\ 60]{Weil:58}). 

The evaluation functional $H_{(2)}^0(X,L)\ni S\mapsto S(x)$ is continuous by
\eqref{eq:1.3}, so by Riesz representation theorem for each $x\in X$ there exists 
$P(x,\cdot)\in \mathcal{L}^2(X,L_x\otimes L^*)$ 
such that
$$s(x)=\int_X P(x,y)s(y)\, \mathrm{dV}(y)\,,\quad
\text{for all $s\in H_{(2)}^0(X,L).$}$$ 
Set
\begin{equation}
	d=\dim H_{(2)}^0(X,L)\in \N\cup\{\infty\}.
\end{equation}
If $X$ is compact, then $d<\infty$. If $d\geq 1$, consider an orthonormal basis 
$\{S_j\}_{j=1}^{d}$ of $H_{(2)}^0(X,L)$. Then the series
$\sum_{j=1}^{d}S_j(x)\otimes (S_j(y))^*$ converges 
uniformly on every compact together with all its derivatives
(cf.\ \cite[Proposition 2.4]{At:76}, \cite[Remark 1.4.3]{MM07}, 
\cite[p.\ 63]{Weil:58}). In particular, $P(x,y)$ is smooth on $X\times X$. 
It follows that 
\begin{equation}
	P(x,y)=\sum_{j=1}^{d}S_j(x)\otimes (S_j(y))^*\,.
	\label{eq:kernel}
\end{equation}
We obtain thus for the Bergman projection \eqref{e:BP},
$$(P s)(x)=\int_X P(x,y)s(y)\,\mathrm{dV}(y),$$ 
i.e., $P(x,y)$ is the integral kernel 
of the Bergman projection.
Recall that the line bundle $L\boxtimes L^*$ on 
$X\times X$ has fibres $(L\boxtimes L^*)_{(x,y)}:=L_{x}\otimes L^{\ast}_{y}$
for $(x,y)\in X\times X$.
The section $P(\cdot,\cdot)$ of $L\boxtimes L^*\to X\times X$ 
is called \textit{Bergman kernel}. 

The canonical identification $L_x\otimes L_x^*=\End(L_x)=\C$,
$s\otimes s^*\mapsto s^*(s)=|s|^2_{h_L}$ allows to identify
$P(x,x)$ to the smooth function
\begin{equation}
	P(x,x)=\sum_{j=1}^{d}|S_j(x)|^2_{h_L},
	\label{eq:kernelfct}
\end{equation}
called the \emph{Bergman kernel function}.
We deduce that
$d=\int_X P(x,x)\, \mathrm{dV}(x)\in \N\cup\{\infty\}.$
Hence,  the Bergman kernel function is the 
dimensional density of $H_{(2)}^0(X,L)$. If $d=0$, then the above 
considerations are trivially true. 

\subsection{Gaussian random holomorphic sections}\label{section2.1}

The results proved in this subsection are extensions of the well-known results for random power series or random analytic 
functions on $\C^{n}$ (cf.\ \cite{Kahane85} or \cite[Section 3]{EK95}) to the complex geometric setting. We include details of the proofs for the sake 
of completeness.

Let $\eta=\{\eta_{j}\}_{j\in\N}$ be a sequence of i.i.d.\ centered real or complex 
Gaussian random variables and denote by $\P$ and $\E$ the underlying probability measure and its expectation.

For $d\geq 1$, let $S=\{S_{j}\}_{j=1}^{d}$ be an orthonormal basis of 
$H^{0}_{(2)}(X,L)$. Define
\begin{equation}
	\psi^{S}_{\eta}(x)=\sum_{j=1}^{d}\eta_{j}S_{j}(x).
	\label{eq:2.2}
\end{equation}
If $d=0$, we simply set $\psi^{S}_{\eta}\equiv 0$.
\begin{Proposition}\label{prop:convergence}
	The section $\psi^{S}_{\eta}$ is almost surely a holomorphic 
	section of $L$ on $X$.
\end{Proposition}
\begin{proof}
	If $d$ is finite, the claim is clearly true. Hence, it remains to prove it 
	for the case $d=\infty.$ In this case, $X$ is noncompact. Let $\{K_{i}\}_{i\in\N}$ be an increasing sequence of compact subsets of 
	$X$ such that $X=\cup_{i\in\N} K_{i}$. We can take each $K_{i}$ 
	to be the closure of a relatively compact open subset $U_{i}$ of $X$. 
	Then to prove this proposition, we only need to show that for 
	each $i$, $\psi^{S}_{\eta}$ is almost surely a holomorphic 
	section of $L$ on $U_{i}$.

	Let $K$ be a compact subset of $X$, and let $U$ be an open 
	relatively compact neighborhood of $K$. 
	Similarly to \eqref{eq:1.3}, there 
	exists a constant $C_{U}>0$ such that for $s\in 
	H^{0}_{(2)}(X,L)$, 
	\begin{equation}
		\sup_{x\in K} |s(x)|_{h_{L}}\leq 
		C_{U}\|s\|_{\cLL(\overline{U},L)}.
		\label{eq:2.2a}
	\end{equation}
By \eqref{eq:kernelfct}, we have 
	\begin{equation}
	\sum_{j=1}^{d}\|S_{j}\|^{2}_{\cLL(\overline{U},L)}=\int_{x\in\overline{U}}P(x,x)\mathrm{dV}(x)<+\infty.
		\label{eq:2.3}
	\end{equation}

	For $j\in\N_{>0}$, $x\in X$, set
	\begin{equation}
		X_{j}(x)=\eta_{j}S_{j}(x),
		\label{eq:2.5}
	\end{equation}
	it is an $L_{x}$-valued random variable. Since $\eta_{j}$ is 
	centered, we infer
	\begin{equation}
		\E[X_{j}(x)]=0\in L_{x}.
	\end{equation}
 	It is then consistent to define the variance as
	\begin{equation*}
		\V(X_{j}(x))=\E[|X_{j}(x)|_{h_{L}}^{2}],
	\end{equation*}
	and we can compute
	\begin{equation}
		\V(X_{j}(x))=\E[|X_{j}(x)|_{h_{L}}^{2}]=
		\V(\eta_{j})|S_{j}(x)|^{2}_{h_{L}}=
		\V(\eta_{1})|S_{j}(x)|^{2}_{h_{L}}.
	\end{equation}	
	We next prove that for any $k\in\N, N\in\N_{>0}$ and for $r>0$, we have
	\begin{equation}
		\P\Big(\sup_{\ell=1,\ldots, N}\sup_{x\in 
		K}\big|\sum_{j=1}^{\ell}X_{k+j}(x)\big|_{h_{L}}>r\Big)<\frac{C^{2}_{U}\V(\eta_{1})}{r^{2}}\sum_{j=1}^{N}\|S_{k+j}\|^{2}_{\cLL(\overline{U},L)}.
		\label{eq:2.8}
	\end{equation}
	For this purpose, define the stochastic process 
	\begin{equation}
	Y_{\ell}=\big\|\sum_{j=1}^{\ell}X_{k+j}\big\|^2_{\cLL(\overline{U},L)},\quad \ell =1, \ldots, N,
	\end{equation}
	 and observe that by virtue of \eqref{eq:2.2a}, we have 
	 \begin{equation} \label{eq:supL2Bd}
	 \sup_{x\in 
		K}\big|\sum_{j=1}^{\ell}X_{k+j}(x)\big|_{h_{L}}\leq 
		C_{U}Y_{\ell}^\frac12.
	\end{equation}
	As a consequence, we have 
	\begin{equation}
	\P\Big(\sup_{\ell=1,\ldots, N}\sup_{x\in 
		K}\big|\sum_{j=1}^{\ell}X_{k+j}(x)\big|_{h_{L}}>r\Big) \le 
				\P\Big(\sup_{\ell=1,\ldots, N} C_{U}^2 Y_{\ell} >r^2\Big).
	\end{equation}
	Now the process $(Y_{\ell})$, $\ell =1, \ldots, N,$ is a submartingale with respect to the filtration $(\mathcal F_\ell),$ where 
	\begin{equation}
	\mathcal F_\ell = \sigma \big(\langle X_{k+i},X_{k+j} \rangle_{\cLL(\overline{U},L)},\  i,j =1, \ldots, \ell \big).
	\end{equation}
Therefore, Doob's submartingale inequality 
(see e.g.\ \cite[Lemma 11.1]{Kl-14}) yields 
	 \begin{equation}
	 \P\Big(\sup_{\ell=1, \ldots, N} Y_{\ell} > \frac{r^2}{C_{U}^2}\Big) \le C_{U}^2 \frac{\E[Y_N]}{r^2},
	 \end{equation}
	 which immediately entails \eqref{eq:2.8}.
	Now, letting $N\rightarrow +\infty$ in \eqref{eq:2.8}, we get
	\begin{equation}
		\P\Big(\sup_{\ell\in\N_{>0}}\sup_{x\in 
		K}|\sum_{j=1}^{\ell}X_{k+j}(x)|_{h_{L}}>r \Big)\leq \frac{C^{2}_{U}\V(\eta_{1})}{r^{2}}\sum_{j=1}^{+\infty}\|S_{k+j}\|^{2}_{\cLL(\overline{U},L)}.
		\label{eq:2.17}
	\end{equation}
	Then taking the limit of \eqref{eq:2.17} as $k\rightarrow \infty$, and
	using \eqref{eq:2.3}, we infer
	\begin{equation}
		\P\Big(\limsup_{k\rightarrow+\infty}\sup_{\ell\in\N_{>0}}\sup_{x\in 
		K}|\sum_{j=1}^{\ell}X_{k+j}(x)|_{h_{L}}>r\Big)=0.
		\label{eq:2.18}
	\end{equation}
	Therefore, a union bound along the sequence of $r = \frac1n$ immediately supplies us with
		\begin{equation}
		\P\Big(\limsup_{k\rightarrow+\infty}\sup_{\ell\in\N_{>0}}\sup_{x\in 
		K}|\sum_{j=1}^{\ell}X_{k+j}(x)|_{h_{L}}>0\Big)=0.
		\label{eq:2.19}
	\end{equation}
	If we take $V$ to be a relatively compact open subset of $X$, and 
	take $K=\overline{V}$, then by \eqref{eq:2.19}, the sum 
	$\sum_{j=1}^{\infty} X_{j}$ is almost surely uniformly convergent 
	on $K$, so that it almost surely defines a holomorphic section on 
	$V$. This completes the proof of our proposition.
\end{proof}
For the purpose of the following definition, 
we note at this point that 
a standard complex Gaussian is a random variable having the distribution 
$\frac{1}{\sqrt{2}} (X + \sqrt{-1}Y)$, 
where $X$ and $Y$ are standard real Gaussian variables.
\begin{Definition}\label{def:flatgaussian}
The random section $\psi^{S}_{\eta}$ defined in \eqref{eq:2.2} 
is called a standard Gaussian random holomorphic section of $L$ over 
$X$ if $\eta=\{\eta_{j}\}_{j\in\N}$ is a sequence of 
i.i.d.\ standard complex Gaussian random variables.
\end{Definition}
Now we prove that the distribution of a standard Gaussian random holomorphic section
$\psi^{S}_{\eta}$ does not depend on the choice of the orthonormal basis.
\begin{Proposition}\label{prop:uniqueness}
Assume that $d\geq 1$, and assume that $\eta=\{\eta_{j}\}_{j=1}^{d}$ 
is a sequence of i.i.d.\ standard complex 
 Gaussian random variables.
	If $S'=\{S'_{j}\}_{j=1}^{d}$ is another choice of orthonormal basis 
	of $H^{0}_{(2)}(X,L)$, then $\psi^{S'}_{\eta}$ and 
	$\psi^{S}_{\eta}$ have the same distribution as random 
	holomorphic sections.
\end{Proposition}
\begin{proof}
	It is sufficient to a sequence 
	$\eta'=\{\eta'_{j}\}_{j=1}^{d}$ of i.i.d.\ standard complex Gaussian random 
	variables such that a.s. 
	$\psi^{S'}_{\eta}=\psi^{S}_{\eta'}$.
	
	Let $\ell^{2}(\C)$ denote the Hilbert space of 
	$\ell^{2}$-summable complex sequences. If 
	$u=(u_{j})_{j\in\N}\in\ell^{2}(\C)$, set 
	\begin{equation}
		(\eta,u)_{\ell^{2}}=\sum_{j\in\N} \eta_{j}\bar{u}_{j}.
		\label{eq:2.20}
	\end{equation}
	By Kolmogorov's Three-Series Theorem (cf.\ \cite{Va90}), the sum in \eqref{eq:2.20} 
	is almost surely convergent, so that $(\eta,u)_{\ell^{2}}$ is a 
	well-defined random variable. By the property of Gaussian random
	variable, we conclude that $(\eta,u)_{\ell^{2}}$ is a centered 
	complex Gaussian random variable with variance 
	$|u|^{2}_{\ell^{2}}$. In particular, if 
	$|u|_{\ell^{2}}=1$, then $(\eta,u)_{\ell^{2}}$ has the same 
	distribution as $\eta_1$. Moreover, if nonzero $u,v\in\ell^{2}$ is such that $(u,v)_{\ell^{2}}=0$, 
	then $(\eta,u)_{\ell^{2}}$ and $(\eta,v)_{\ell^{2}}$ are 
	independent.
	
	Take $(a_{ij}\in\C)_{i,j\in\N}$ such that for each $i$,
	\begin{equation}
		S'_{i}=\sum_{j\in\N}a_{ij}S_{j}.
		\label{eq:2.21}
	\end{equation}
	For $j\in\N$, set $b_{j}=(\bar{a}_{ij})_{i\in\N}$. Then 
	$b_{j}\in\ell^{2}(\C)$ is with norm $1$, moreover, if $j\neq j'$, 
	then $(b_{j},b_{j'})_{\ell^{2}}=0$. Now define
	\begin{equation}
		\eta'_{j}=(\eta,b_j)_{\ell^{2}}.
		\label{eq:2.22}
	\end{equation}
	Then $\eta'=(\eta'_{j})_{j\in\N}$ is a sequence of i.i.d.\ centered 
	Gaussian random variables with the same distribution as $\eta$. 
	By definition, we get that almost surely,
	\begin{equation}
		\psi^{S'}_{\eta}=\psi^{S}_{\eta'}.
		\label{eq:2.23aaa}
	\end{equation}
	Therefore, $\psi^{S'}_{\eta}$ and $\psi^{S}_{\eta}$ have the same 
	distribution.
\end{proof}

\begin{Remark}\label{rem:Frechet}
(a) When $d=\infty$, note that by taking a sequence of compact subset $\{K_{i}\}_{i\in 
	\N}$ as in the proof of Proposition \ref{prop:convergence}, we 
	can define a sequence of semi-norms for $H^{0}(X,L)$, hence a 
	Fr\'{e}chet distance, so that $H^{0}(X,L)$ is a Fr\'{e}chet 
	space. In Proposition \ref{prop:convergence}, we actually prove that $\psi^{S}_{\eta}$ is a random variable taking values in the 
	Fr\'{e}chet space
	$H^{0}(X,L)$.

\noindent
(b) In the proof of Proposition \ref{prop:convergence}, we do 
	not use the Gaussianity of the $\eta_{j}$ in an essential way. Hence,  we can 
	work with any sequence $\eta$ of pairwise uncorrelated centered random variables with 
	uniformly bounded variance. In that case, however, the distribution of the random section 
	$\psi^{S}_{\eta}$ might depend on the choice of the 
	basis $S$. Generally one needs suitable moment conditions on $\eta$ 
	to obtain the good behaviours such as the universality results of the zeros of $\psi^{S}_{\eta}$, we 
	refer to \cite{KZ14}, \cite{BCM}, \cite{DLM:21} for the related 
	details.
\end{Remark}

\begin{Lemma}\label{lem:L2}
	If $d=\infty$, then with probability one, $\psi^{S}_{\eta}$ is not $\cLL$-integrable on 
	$X$.
\end{Lemma}
\begin{proof}
	The event that $\psi^{S}_{\eta}$ is $\cLL$-integrable is equivalent
	to the event  $\{\sum_{j=1}^{\infty} 
	|\eta_{j}|^{2}< \infty\}$. But e.g.\ by the law of large numbers, we infer 
	\begin{equation}
		\P\Big(\sum_{j=1}^{\infty} |\eta_{j}|^{2}< \infty\Big)=0,
	\end{equation}
	and the statement of the lemma follows.
\end{proof}

\subsection{Expectation of random zeros: proof of {Theorem \ref{thm:expectation}}}\label{section2.2}
In the sequel we always assume $d=\dim H^{0}_{(2)}(X,L)\geq 1$. 
We start with some considerations about the Fubini-Study currents.
\begin{Lemma}\label{lem:localintegral}
Assume that $d\geq 1$.
Then the function $X\ni x\mapsto \log P(x,x) \in 
\{-\infty\}\cup \R$ is locally $\mathcal{L}^{1}$-integrable on 
	$X$ with respect to $\mathrm{dV}$. Thus  
	$\dfrac{i}{2\pi}\partial\bar{\partial}\log P(x,x)$ 
	defines a $(1,1)$-current on $X$.
\end{Lemma}
\begin{proof}
Let $e_{L}:U\to L$ be a local holomorphic frame of $L$.
Let $\varphi\in\mathscr{C}^\infty(U)$ be the local weight
of $h_L$ with respect to $e_L$, that is, $|e_L|^2_{h_L}=e^{-2\varphi}$
on $U$. We consider an orthonormal basis $\{S_j\}_{j=1}^d$ of
$H^0_{(2)}(X,L)$ and write $S_{j}(x)=f_{j}(x)e_{L}(x)$, 
$x\in U$, where $f_{j}$ a non-trivial holomorphic functions on $U$. 
Then $P(x,x)=\sum_{j=1}^d|S_j|^2_{h_L}=\sum_{j=1}^d|f_j|^2e^{-2\varphi}$
on $U$, hence
\begin{gather}
\log P(x,x)=\log\Big(\sum_{j=1}^d|f_j|^2\Big)-2\varphi\,.\label{e:logP}
\end{gather}
The series $\sum_{j=1}^d|f_j|^2$ converges locally uniformly on $U$,
thus $\log\!\big(\sum_{j=1}^d|f_j|^2\big)$ is a plurisubharmonic
function that is not identically $-\infty$, hence locally integrable.
\end{proof}
Lemma \ref{lem:localintegral} shows that the 
Fubini-Study currents \eqref{eq:intro2.23} are well defined.
Note that 
$c_1(L,h_L)|_U=\frac{\sqrt{-1}}{\pi}\partial\overline{\partial}\varphi$.
By applying $\partial\overline{\partial}$ on both sides of
\eqref{e:logP} and taking into account \eqref{eq:intro2.23} we see
that 
\begin{equation}\label{e:FS}
\gamma(L,h_{L})\big|_U=\frac{\sqrt{-1}}{2\pi}\partial\overline{\partial}
\log\big(\sum_{j=1}^d|f_{j}|^{2}\big)\,,
\end{equation}
thus $\gamma(L,h_{L})$ is a closed positive $(1,1)$-current.
The base locus of $H^{0}_{(2)}(X,L)$ is the proper analytic set
\begin{equation}
	\Bl(X,L):=\big\{x\in X\;|\; \text{$s(x)=0$ for all 
	$s\in H^{0}_{(2)}(X,L)$}\big\}.
	\label{eq:intro2.1}
\end{equation}
Thus $\{x\in X:P(x,x)=0\}=\Bl(X,L)$. Hence $\gamma(L,h_{L})$
is a smooth form if $\Bl(X,L)=\varnothing$. 
Note that if $X$ is compact and $\Bl(X,L)=\varnothing$, then
$\gamma(L,h_{L})$
is the pullback of the Fubini-Study form on the projective space
by the Kodaira map defined by $H^0(X,L)$. 
This justifies the name of Fubini-Study currents.

Now we are ready to prove Theorem \ref{thm:expectation}. 
Let $\eta=\{\eta_{j}\}_{j=1}^{d}$ is a sequence of 
i.i.d.\ standard complex Gaussian random variables.
Let $\psi^{S}_{\eta}$ be the random holomorphic section defined in 
\eqref{eq:2.2} and let $[\Div(\psi^{S}_{\eta})]$ denote the 
$(1,1)$-current given by its zeros (cf.\ \eqref{eq:divisor}).

\begin{proof}[Proof of {Theorem \ref{thm:expectation}}]
	By the Lelong-Poincar\'{e} formula (cf.\ \cite[Theorem 
	2.3.3]{MM07}), for $s\in H^{0}(X,L)$, we have
\begin{equation}
	[\Div(s)]=\frac{\sqrt{-1}}{2\pi}
	\partial\bar{\partial}\log|s|^{2}_{h_{L}}+ c_{1}(L,h_{L}).
	\label{eq:2.25}
\end{equation}

	Fix a test form $\varphi\in \Omega^{(n-1,n-1)}_{0}(X)$, and we 
	evaluate $\E[\langle [\Div(\psi^{S}_{\eta})], 
	\varphi\rangle]$. Apply \eqref{eq:2.25} to $\psi^{S}_{\eta}$, we 
	get 
	\begin{equation}
		\begin{split}
			\langle [\Div(\psi^{S}_{\eta})],\varphi\rangle &= 
			\int_{X}\big(\frac{\sqrt{-1}}{2\pi}\partial\bar{\partial}\log|\psi^{S}_{\eta}|^{2}_{h_{L}}+ c_{1}(L,h_{L})\big)\wedge \varphi\\
			&=\int_{X}c_{1}(L,h_{L})\wedge\varphi+\frac{\sqrt{-1}}{2\pi}\int_{X}\log|\psi^{S}_{\eta}|_{h_{L}}^{2} \partial\bar{\partial}\varphi\\
			&=\int_{X}c_{1}(L,h_{L})\wedge\varphi+\frac{\sqrt{-1}}{2\pi}\int_{X\backslash \text{Bl}(X,L)}\log|\psi^{S}_{\eta}|_{h_{L}}^{2} \partial\bar{\partial}\varphi.
		\end{split}
		\label{eq:2.27}
	\end{equation}	
For $x\in X\backslash \text{Bl}(X,L)$, we have $P(x,x)\neq 0$, 
and let $e_{L}(x)$ be a unit vector of $L$ at $x$, define
	\begin{equation}
		b(x)=\big(P(x,x)^{-1/2}S_{j}(x)/e_{L}(x)\big)_{j\in\N}\in 
		\ell^{2}(\C).
		\label{eq:2.28}
	\end{equation}
	We have $|b(x)|_{\ell^{2}}=1$.  Note that
	\begin{equation}
		P(x,x)^{-1/2}\psi^{S}_{\eta}=(\eta,\overline{b(x)})_{\ell^{2}}\,e_{L}(x).
		\label{eq:2.29}
	\end{equation}
	Then 
	\begin{equation}
		\E[\log|P(x,x)^{-1/2}\psi^{S}_{\eta}|_{h_{L}}^{2}]=\E[\log\big|(\eta,\overline{b(x)})_{\ell^{2}}\big|^{2}]=\E[\log|\eta_{1}|^{2}]
		\label{eq:2.30}
	\end{equation}
	Note that $\E[|\log|\eta_{1}|^{2}|]<\infty$. By 
	Lemma \ref{lem:localintegral}, $\log P(x,x)$ is locally 
	integrable on $X$, then we can apply the Fubini's theorem to the 
	following integrals so that
	\begin{equation}
		\begin{split}
			&\E\big[\int_{X\backslash \text{Bl}(X,L)}\log|P(x,x)^{-1/2}\psi^{S}_{\eta}|_{h_{L}}^{2} 
		\partial\bar{\partial} \varphi\big]\\
		&=\int_{X\backslash \text{Bl}(X,L)}\E\big[\log|P(x,x)^{-1/2}\psi^{S}_{\eta}|_{h_{L}}^{2}\big]
		\partial\bar{\partial} \varphi\\
		&=\E[\log|\eta_{1}|^{2}]\int_{X}\partial\bar{\partial} 
		\varphi=0.
		\end{split}
		\label{eq:2.31}
	\end{equation}
	Then by \eqref{eq:2.27}, we get
	\begin{equation}
		\begin{split}
			 \E\left[\langle[\Div(\psi^{S}_{\eta})],\varphi\rangle 
			 \right]
			&=\int_{X}c_{1}(L,h_{L})\wedge\varphi+\frac{\sqrt{-1}}{2\pi}\int_{X\backslash \text{Bl}(X,L)}\log{P(x,x)}\cdot \partial\bar{\partial}\varphi\\
			&=\big\langle 
			c_{1}(L,h_{L})+\frac{\sqrt{-1}}{2\pi}\partial\bar{\partial}\log(P(x,x)), \varphi\big\rangle.
		\end{split}
		\label{eq:2.40ss}
	\end{equation}
This completes the proof.
\end{proof}

\subsection{Geometric examples}\label{ss:GeoExa}
We present in this subsection some interesting examples of 
Bergman spaces and Fubini-Study currents where
our result apply.
We start with some simple observations.

(i) If $P(x,x)>0$ (equivalently, $x$ is not in the base locus
of $H^0_{(2)}(X,L)$), then the $(1,1)$-form
$\sqrt{-1}\partial\dbar\log P(x,x)$ is smooth in a neighborhood of $x$, 
and hence $\gamma(L,h_L)$, too.
In particular, if $\Bl(X,L)=\varnothing$,
then $\gamma(L,h_L)$ is smooth.

(ii) If $P(x,x)>0$ let $s_0\in H^0_{(2)}(X,L)$ with $s_0(x)\neq0$.
Assume that there exist $s_1,\ldots,s_n\in H^0_{(0)}(X,L)$
such that $d(s_1/s_0)(x),\ldots, d(s_n/s_0)(x)$ are linearly 
independent (that is, sections of $H^0_{(2)}(X,L)$ give local coordinates
at $x$). Then $\sqrt{-1}\partial\dbar\log P(x,x)$ is strictly positive
near $x$. 

(iii) Thus, if $\Bl(X,L)=\varnothing$ and sections
of $H^0_{(2)}(X,L)$ give local coordinates
at any point in $X$, then $\sqrt{-1}\partial\dbar\log P(x,x)$
defines a K\"ahler metric on $X$. 

\begin{Example}[Bergman metric]\label{Ex:BM}
We consider the case when $L$ is the canonical bundle
$K_X$ of $X$ (cf.\ \cite{Weil:58}). The space of holomorphic sections of $K_X$ is the 
space $H^{n,0}(X)$ of holomorphic
$(n,0)$-forms. Such a form can be written in local 
coordinates $(z_1,\ldots,z_n)$ as 
$f(z)dz_1\wedge\ldots\wedge dz_n$, with $f$ a
holomorphic function. 
We say that a measurable $(n,0)$-form $\beta$ is an $\cLL$ section
of $K_X$ if 
\begin{equation}
\label{e:L2form}
\|\beta\|^2:=2^{-n}(\sqrt{-1})^{n^2}
\int_X\beta\wedge\overline\beta<\infty.
\end{equation} 
We denote by $H^{n,0}_{(2)}(X)$ the space of $\cLL$ holomorphic
$(n,0)$-forms. We have $H^{n,0}_{(2)}(X)=H^{0}_{(2)}(X,K_X)$,
where the right-hand side is defined with respect to an arbitrary
metric $\Theta$ on $X$ and
the metric on $K_X$ is induced by $\Theta$.

We assume that $H^{n,0}_{(2)}(X)\neq\{0\}$ and let 
$\{\beta_j\}_{j=1}^d$ be an orthonormal basis of
$H^{n,0}_{(2)}(X)$. 
In local coordinates $(U;z_1,\ldots,z_n)$ 
write $\beta_j=f_j(z)dz_1\wedge\ldots\wedge dz_n$.
According to \eqref{e:FS} the Fubini-Study current
is given on $U$ by  
$\gamma(K_X,h_{K_X})|_U=\frac{\sqrt{-1}}{2\pi}\partial\dbar\log(\sum_j |f_j|^2)$.
If the the Fubini-Study current is actually a K\"ahler metric
on $X$, then it is called the \emph{Bergman metric} of $X$. 
We will denote it by $\omega_B$.
The metric $\omega_B$ is invariant by the group
of biholomorphic transformations of $X$.

If $X$ is an open set in $\C^n$, the canonical bundle is trivial, so we identify
the space $H^{n,0}_{(2)}(X)$ of $\cLL$-holomorphic 
$(n,0)$-forms to the space $H^{0}_{(2)}(X)$ holomorphic functions
which are $\cLL$ with respect to the Lebesgue measure. 
There is a vast literature on Bergman spaces and kernels on domains
in $\C^n$, see e.g. \cite{HKZ00,JP93} and the references therein.
\end{Example}

To give concrete examples let us recall the definition of Stein manifolds, which
are  interesting due to their rich function-theoretical structure
\cite{GR04}. 
For a complex manifold $X$, let 
$\mathcal{O}(X)$ denote the space of all holomorphic functions on $X$.
\begin{Definition}
A complex manifold $X$ is called Stein if the following two 
conditions are satisfied:
(1) $X$ is homomorphically convex, i.e., for every compact 
		subset $K\subset X$, its holomorphically convex hull 
		$\widehat{K}=\big\{z\in X: |f(z)|\leq \sup_{w\in K}|f(w)|, 
			\forall\, f\in\mathcal{O}(X)\big\}$
		 is compact.
(2) $X$ is holomorphically separable, i.e., if $x\neq y$ in 
		$X$, then there exists $f\in \mathcal{O}(X)$ such that 
		$f(x)\neq f(y)$.
\end{Definition}
Let $L\to X$ be a holomorphic line bundle.
The cohomology vanishing theorem for coherent 
analytic sheaves on Stein manifolds (Cartan's theorem B, cf.\ \cite{GR04})
yields the following: 

(i) The holomorphic
sections $H^0(X,L)$ give local coordinates at each point of $X$.

(ii) For any closed discrete set 
$A=\{p_k:k\in\N\}$ and any family $\{v_k\in L_{p_k}:k\in\N\}$ there exists
$s\in H^0(X,L)$ with $s(p_k)=v_k$ for all $k\in\N$.
In particular, for each $p\in X$ the evaluation map
$H^0(X,L)\to L_p$ is surjective and we have $\dim H^0(X,L)=\infty$. 

\begin{Example}
Let $X$ be a Stein manifold and $D\Subset X$ be a relatively
compact domain. We consider a Hermitian metric on $X$
whose associated $(1,1)$-form is denoted by $\Theta$.
Let $\mathrm{dV}_\Theta=\Theta^n/n!$ be the volume form induced by $\Theta$,
where $\dim X=n$. Let $(L,h_L)$ be a Hermitian holomorphic line bundle.
Consider the space $\cLL(D,L,h_L,\mathrm{dV}_\Theta)$ of measurable
sections $S$ of $L$ over $D$ satisfying $\int_D|S|^2_{h_L}\mathrm{dV}_\Theta<\infty$
and let $H^0_{(2)}(D,L,h_L,\mathrm{dV}_\Theta)=\cLL(D,L,h_L,\mathrm{dV}_\Theta)\cap H^0(X,L)$.
The restriction map $H^0(X,L)\to H^0_{(2)}(D,L,h_L,\mathrm{dV}_\Theta)$
is well-defined and injective. 
We deduce that the space $H^0_{(2)}(D,L,h_L,\mathrm{dV}_\Theta)$ is infinite
dimensional, has empty base locus and sections of this space
give local coordinates at any point of $D$. Therefore, $\gamma(L,h_L)$
is smooth on $X$ and if $(L,h_L)$ is semipositive (i.e.\ $c_1(L,h_L)$
is positive semidefinite), it is a K\"ahler form.

We deduce from Theorem \ref{thm:expectation} and 
the discussion from Example 
\ref{Ex:BM} the following.
\begin{Corollary}
For any relatively
compact domain $D\Subset X$ in a Stein manifold 
the expectation of the zero divisors of the standard Gaussian random 
holomorphic $(n,0)$-forms defined from the $\cLL$-holomorphic $(n,0)$-forms 
on $D$ is given by the Bergman metric
on $D$. If $D\Subset\C^n$ this is true for
standard Gaussian random 
holomorphic functions defined from the $\cLL$-holomorphic functions
on $D$.
\end{Corollary}
One of the simplest examples is the unit disc $\D\subset\C$
endowed with the Lebesgue measure.
Then $P(z,z)=\frac{1}{\pi(1-|z|^2)^2}$ and the Bergman metric
\begin{equation}
	\omega_B=\frac{\sqrt{-1}}{\pi}\frac{dz\wedge 
	d\bar{z}}{(1-|z|^{2})^{2}}
\end{equation}
is the hyperbolic metric (up to a constant factor) on the disc. 
We see on this example that the Bergman metric explodes
for $|z|\to1$, so the zeros accumulate
towards the boundary of $\D$. This is a more general phenomenon,
in the sense that the Bergman metric is complete if $D$ is a
domain of holomorphy in $\C^n$.
\end{Example}

\begin{Example}[Bargmann-Fock space: flat Gaussian holomorphic function]\label{ex:2.9}
	Let $L$ be the trivial line bundle on $\C^{n}$ but we equip it with 
	the Hermitian metric $h_{L}$ such that $|1|^{2}_{h_{L},z}=e^{-|z|^{2}}$, 
	$z\in\C^{n}$. In this case,
	\begin{equation}
		R^{L}=\sum_{j=1}^{n}dz_{j}\wedge d\bar{z}_{j}.
		\label{eq:2.32}
	\end{equation}
	We endow $\C^{n}$ with the flat metric 
	$\Theta=\frac{\sqrt{-1}}{2\pi}\sum_{j=1}^{n}dz_{j}\wedge d\bar{z}_{j}$, 
	then
	\begin{equation}
		\mathrm{dV}_{\Theta}=\frac{1}{\pi^{n}}\Pi_{j=1}^{n}dx_{j}\wedge dy_{j}.
		\label{eq:2.33}
	\end{equation}
	For a 
	multi-index
	$\alpha=(\alpha_{1},\ldots,\alpha_{n})\in \N^{n}$, we write
	\begin{equation}
		S_{\alpha}(z)=\frac{z_{1}^{\alpha_{1}}\ldots 
		z_{n}^{\alpha_{n}}}{\sqrt{\alpha_{1}!\ldots\alpha_{n}!}}.
		\label{eq:2.34}
	\end{equation}
A straightforward calculation then confirms that 
	$\{S_{\alpha}\}_{\alpha\in\N^{n}}$ forms an orthonormal basis of 
	$H^{0}_{(2)}(\C^{n},L)$.
In this case, we have
	\begin{equation}
		P(z,z)=\sum_{\alpha\in \N^{n}} \frac{|z_{1}|^{2\alpha_{1}}\ldots 
		|z_{n}|^{2\alpha_{n}}}{\alpha_{1}!\ldots\alpha_{n}!}e^{-|z|^{2}}=1.
		\label{eq:2.35}
	\end{equation}
Denoting by $\eta=(\eta_{\alpha})_{\alpha\in\N^{n}}$ a family of 
i.i.d.\ standard complex 
Gaussian random variables, we define the standard Gaussian random 
holomorphic function on $\C^{n}$ as
\begin{equation}
	\psi^{S}_{\eta}=\sum_{\alpha\in\N^{n}} \eta_{\alpha}S_{\alpha}.
	\label{eq:2.36}
\end{equation}
By Theorem \ref{thm:expectation}, we have
\begin{equation}
	\E[[\Div(\psi^{S}_{\eta})]]=\gamma(L, h_{L})=\frac{\sqrt{-1}}{2\pi}\sum_{j=1}^{n}dz_{j}\wedge d\bar{z}_{j}.
	\label{eq:2.37}
\end{equation}
\end{Example}

\section{Equidistribution and large deviation for high 
tensor powers of line bundles}\label{section3}
In 
the sequel, assume that $\eta=\{\eta_{j}\}_{j\in\N}$ is a sequence of 
i.i.d.\ standard complex Gaussian random variables, note that 
$\V(\eta_{1})=1$. 

In this section, we consider the setting of Subsection \ref{s1.2a}, 
in particular, we assume \eqref{eq:intro3.0.1}. Let $\dot{R}^{L}\in\mathrm{End}(T^{(1,0)}X)$ such that $x\in X$,
for $u,v\in T_{x}^{(1,0)}X$,
\begin{equation}
	R^{L}_{x}(u,v)=g^{TX}_{x}(\dot{R}^{L}u,v).
	\label{eq:3.3a}
\end{equation}
By \eqref{eq:intro3.0.1}, we have
$a_{0}(x)=\det \dot{R}^{L}_{x}\geq \varepsilon^{n}$.

\subsection{Equidistribution of zeros of Gaussian random holomorphic sections}\label{ss3.1a}
We consider the sequence of Hilbert spaces $H^{0}_{(2)}(X,L^{p})$,
$p\in\N$ large. Set 
\begin{equation}
	d_{p}=\dim H^{0}_{(2)}(X,L^{p})\in 
\N\cup\{\infty\}.
	\label{eq:2.39}
\end{equation}
We equip $L^{p}$ with the induced Hermitian metric $h_{p}:=h^{\otimes 
p}_{L}$.
Let $P_{p}$ denote the orthogonal projection from $\cLL(X,L^{p})$ 
onto $H^{0}_{(2)}(X,L^{p})$, and let $P_{p}$ 
denote the corresponding Bergman kernel on $X$ with 
respect to $\mathrm{dV}(x)=\frac{\Theta^{n}}{n!}$.


For $p\in\N_{>0}$, let $\psi^{S_{p}}_{\eta}$ be a standard Gaussian 
random holomorphic section constructed from
	$H^{0}_{(2)}(X,L^{p}),$ i.e., for $\{S^{p}_{j}\}_{j=1}^{d_{p}}$ 
	an orthonormal basis of $H^{0}_{(2)}(X,L^{p})$ with respect to 
	the $\cLL$-metric, and set
	\begin{equation}
		\psi^{S_{p}}_{\eta}=\sum_{j=1}^{d_{p}}\eta_{j} S^{p}_{j}.
	\end{equation}

\begin{Theorem}\label{thm:asym}
We assume that Riemannian metric $g^{TX}$ is complete and 
\eqref{eq:intro3.0.1} holds. 
Then as $p\rightarrow +\infty$, we have the weak convergence 
	\begin{equation}
		\frac{1}{p}\E[[\Div(\psi^{S_{p}}_{\eta})]]\rightarrow 
		c_{1}(L,h_{L})
		\label{eq:2.40}
	\end{equation}
	of 
	$(1,1)$-currents, i.e., for any $\varphi\in 
	\Omega^{(n-1,n-1)}_{0}(X)$, as $p\rightarrow +\infty$,
	\begin{equation}
		\Big \langle\frac{1}{p}\E[[\Div(\psi^{S_{p}}_{\eta})]],\varphi\Big \rangle\rightarrow 
		\langle c_{1}(L,h_{L}),\varphi\rangle.
		\label{eq:2.42}
	\end{equation}
	On any given compact subset $K\subset X$ and for sufficiently large 
	$p\in \N$, $\dfrac{1}{p}\E[[\Div(\psi^{S_{p}}_{\eta})]]$ 
	is a smooth $(1,1)$-form, and \eqref{eq:2.40} holds
	in the $\mathscr{C}^{\ell}(K)$-norm for any $\ell\in\N$.
\end{Theorem}
\begin{proof}
	By Theorem \ref{thm:expectation} for $L^{p}$, 
	 \begin{equation}
		 \E[[\Div(\psi^{S_{p}}_{\eta})]]=\gamma(L^{p},h_{p}),
		\label{eq:2.41}
	\end{equation}
	where $\gamma(L^{p}, h_{p})$ is the corresponding Fubini-Study 
	current defined via \eqref{eq:intro2.23}.
	
	By \cite[Theorem 
	 6.1.1]{MM07}, for a given 
	 compact subset $K$ of $X$ and for sufficiently large $p$, 
	 $P_{p}(x,x)>0$, so that
	$\gamma(L^{p}, h_{p})$ is a smooth 
	$(1,1)$-form on $K$. 
	Applying \cite[Corollary 6.1.2]{MM07} to $\gamma(L^{p},h_{p})$, 
	for any $l\in\N$, we 
	have the $\mathscr{C}^{l}(K)$-convergence as $p\rightarrow 
	+\infty$,
	\begin{equation}
		\frac{1}{p}\gamma(L^{p},h_{p})\rightarrow c_{1}(L,h_{L}).
	\end{equation}
	In particular, for any given test form $\varphi\in 
	\Omega^{(n-1,n-1)}_{0}(X)$, we get \eqref{eq:2.42}.
\end{proof}

The convergence in \eqref{eq:2.40} can be improved by imposing further
geometric assumptions, for instance the assumption of bounded geometry. 
We say that $(X,J,\Theta)$, $(L,h_{L})$ have bounded geometry if $J$, 
$g^{TX}$, $R^{L}$ and their derivatives of any order are uniformly 
bounded on $X$ in the norm induced by $g^{TX}$, and the injective 
radius of $(X,g^{TX})$ is strictly positive. One important example of complex manifolds of bounded geometry is the Galois 
coverings of a compact K\"{a}hler manifold $M$ by the Deck 
transformations, and taking the line bundle to be the pull-back a 
positive holomorphic line bundle on $M$.

We recall the following results 
proved in \cite[Theorem 3]{MM15}.

\begin{Theorem}[\cite{MM15}]
	Under the assumptions of bounded geometry and of \eqref{eq:intro3.0.1}, 
	we have the expansion 
	\begin{equation}
		P_{p}(x,x)=\frac{a_{0}(x)}{(2\pi)^{n}}p^{n}+\mathcal{O}(p^{n-1})
		\label{eq:3.10a}
	\end{equation}
	 in the $\mathscr{C}^\infty$-topology on $X$.
	 
	Moreover, there exists $p_{0}\in\N$ such that for all $p>p_{0}$, $X$ is 
	holomorphically convex with respect to the bundle $L^{p}$ and 
	$H^{0}_{(2)}(X,L^{p})$ separates points and gives local 
	coordinates on $X$. 
\end{Theorem}

As a consequence, we get the 
following results.
\begin{Proposition}\label{prop:boundedasym}
Assume bounded geometry as well as
\eqref{eq:intro3.0.1}. Writing $\psi^{S_{p}}_{\eta}$ 
for the Gaussian random section constructed from 
	$H^{0}_{(2)}(X,L^{p})$, then for sufficiently large $p$, 
	$\E[[\Div(\psi^{S_{p}}_{\eta})]]$ is a smooth 
	$(1,1)$-form on $X$. Then we have  
	\begin{equation}
\frac{1}{p}\,\E\left[[\Div(\psi^{S_{p}}_{\eta})]\right]\rightarrow 
c_{1}(L,h_{L}), \:\:\text{as $p\rightarrow +\infty$,
	in the  $\mathscr{C}^{\infty}$-topology on $X$.}
\label{eq:2.40bis}
\end{equation}
\end{Proposition}

\begin{Remark}
Note that under the assumption of bounded geometry and for $X$ 
noncompact, we have $d_{p}=\infty$, $p\gg 0$.
\end{Remark}

\begin{Example}[Scaled Bargmann-Fock spaces]\label{ex:3.6}
	We consider the line bundle $(L,h_{L})$ on $\C^{n}$ from
	Example \ref{ex:2.9}, which satisfies the above 
	assumptions. For $p\geq 1$, an orthonormal basis of 
	$H^{0}_{(2)}(\C^{n},L^{p})$ is given by the family
	\begin{equation}
		S^{p}_{\alpha}(z)=p^{\frac{n}{2}}S_{\alpha}(\sqrt{p}z), \quad
		\alpha\in\N^{n}.
		\label{eq:3.8ss}
	\end{equation}
	Then the Bergman kernel function is given
	\begin{equation}
		P_{p}(z,z)\equiv p^{n}.
		\label{eq:ppconstant}
	\end{equation}
	
	Recall the flat Gaussian random holomorphic function 
	$\psi^{S}_{\eta}$ on $\C$ is defined by \eqref{eq:2.26}. Then for $p\geq 
	1$, we have
	\begin{equation}
		\psi^{S_{p}}_{\eta}(z)=p^{n/2}\psi^{S}_{\eta}(\sqrt{p}z).
		\label{eq:flatguassian}
	\end{equation}
	A direct computation then shows that
	\begin{equation}
		\frac{1}{p}\E[[\Div(\psi^{S_{p}}_{\eta})]]=\E[[\Div(\psi^{S}_{\eta})]]=\frac{\sqrt{-1}}{2\pi}\sum_{j=1}^{n}dz_{j}\wedge d\bar{z}_{j}.
	\end{equation}
\end{Example}

\begin{Theorem}\label{thm:2.15}
	Let $(L,h_{L})$ and $(X,\Theta)$ be as in Theorem \ref{thm:asym}. 
	For any given test form $\varphi\in 
	\Omega^{(n-1,n-1)}_{0}(X)$, we have
	\begin{equation}
		\P\Big(\lim_{p\rightarrow+\infty} \frac{1}{p}\big\langle
		[\Div(\psi^{S_{p}}_{\eta})],\varphi\big\rangle= \langle 
		c_{1}(L,h_{L}),\varphi\rangle\Big)=1.
		\label{eq:2.43}
	\end{equation}
\end{Theorem}
\begin{proof}
To prove this theorem, we mainly follow the arguments from 
proof of \cite[Theorem 5.3.3]{MM07}, 
and the possibility of infinite dimension does not lead to complications in this setting. 
Fix a non-trivial test form $\varphi\in 
\Omega^{(n-1,n-1)}_{0}(X)$.
	Note that from the proof of Theorem \ref{thm:asym}, we have the convergence 
	\begin{equation}
		\lim_{p\rightarrow\infty}\Big\langle 
		\frac{1}{p}\gamma(L^{p},h_{p}),\varphi\Big\rangle=\langle 
		c_{1}(L,h_{L}),\varphi\rangle.
	\end{equation}
Defining the random variable
\begin{equation}
	Y_{p}=\frac{1}{p}\Big\langle 
	[\Div(\psi^{S_{p}}_{\eta})]-\gamma(L^{p},h_{p}),\varphi\Big\rangle,
\end{equation}
 statement \eqref{eq:2.43} is equivalent to proving the almost 
sure convergence
\begin{equation}
	Y_{p}\rightarrow 0.
\end{equation}
For any $x\in \supp \varphi$, let $e_{L}(x)$ denote a unit vector of 
$(L_{x},h_{L,x})$. Set
\begin{equation}
	b_{p}(x)=(P_{p}(x,x)^{-1/2}S^{p}_{j}(x)/e_{L}^{\otimes 
p}(x))\in\ell^{2}(\C).
\end{equation}
Then $\eta\cdot b_{p}(x)$ is a standard complex Gaussian 
variable. The covariance 
matrix of the Gaussian 
vector $(\eta\cdot b_{p}(x),\eta\cdot b_{p}(y))$ depends 
smoothly on $(x,y)\in \supp\varphi\times\supp\varphi$.

For $v=(v_{1},v_{2})\in\C^{2}$ with $\|v\|=1$, we consider the 
integral
\begin{equation}
	\rho(v):=\frac{1}{4\pi^{2}}\int_{\C^{2}}e^{-\frac{1}{2}(|z_{1}|^{2}+|z_{2}|^{2})}\big|\log|z_{1}|\cdot\log|v_{1}z_{1}+v_{2}z_{2}|\big|\mathrm{dV}(z).
\end{equation}
The computations in \cite[Eqs.\ (5.3.13) to (5.3.15)]{MM07} then show that
\begin{equation}
	C:=\sup_{v\in\C^{2},\|v\|=1 }\rho(v)<\infty,
\end{equation}
so for $x,y\in\supp\varphi$ we have
\begin{equation}
\E\Big[\Big|\log\big|P_{p}(x,x)^{-1/2}\sum_{j}\eta_{j}S^{p}_{j}(x)
\big|_{h_{p}}\log\big|P_{p}(y,y)^{-1/2}
\sum_{j}\eta_{j}S^{p}_{j}(y)\big|_{h_{p}}\Big|\Big]\leq C.
	\label{eq:3.26a}
\end{equation}
Note that
\begin{equation}
\E[|Y_{p}|^{2}]=\frac{1}{p^{2}}\E\big[\left|\langle 
[\Div(\psi^{S_{p}}_{\eta})],\varphi\rangle\right|^{2}\big]-\frac{1}{p^{2}}\big|\langle \gamma(L^{p},h_{p}),\varphi\rangle\big|^{2}.
\end{equation}
Then by \eqref{eq:2.25}, \eqref{eq:3.26a} and the 
Fubini-Tonelli theorem we infer that
\begin{equation}
\begin{split}
\E[|Y_{p}|^{2}]=\frac{1}{\pi^{2}p^{2}}\int_{X\times X} &
\big(\partial\overline{\partial}\varphi(x)\big)
\big(\overline{\partial\overline{\partial}\varphi(y)}\big)\\
&\E\Big[\log\big|P_{p}(x,x)^{-1/2}\sum_{j}\eta_{j}S^{p}_{j}(x)
\big|_{h_{p}}\log\big|P_{p}(y,y)^{-1/2}
\sum_{j}\eta_{j}S^{p}_{j}(y)\big|_{h_{p}}\Big],
\end{split}
\label{eq:3.22a}
\end{equation}

By \eqref{eq:3.26a} we conclude
\begin{equation}
	\E[|Y_{p}|^{2}]=\mathcal{O}\Big(\frac{1}{p^{2}}\Big).
\end{equation}
Hence $\E[\sum_{p\geq 1}|Y_{p}|^{2}]=\sum_{p\geq 1}\E[|Y_{p}|^{2}]<\infty$,
thus $Y_{p}\rightarrow0$ almost surely.
\end{proof}
\begin{Corollary}\label{cor:3.8a}
 If $\Theta$ is a K\"ahler form and $\int_{X} 
	c_{1}(L,h_{L})\wedge \Theta^{n-1}<\infty$, or if $\int_{X} 
	c_{1}(L,h_{L})^{n}<\infty$, then 
	\begin{equation}
		\P\Big(\lim_{p\rightarrow+\infty} \frac{1}{p}
		[\Div(\psi^{S_{p}}_{\eta})]= 
		c_{1}(L,h_{L})\Big)=1,
		\label{eq:2.47}
	\end{equation}
	where the limit is taken with respect to 
	the weak convergence of $(1,1)$-currents on $X$.
\end{Corollary}
\begin{proof}
	Due to the assumptions, there exists a constant 
	$C>0$ such that for all $\varphi\in 
	\Omega^{(n-1,n-1)}_{0}(X)$ and $s_{p}\in H^{0}(X,L^{p})$,
	\begin{equation}
		\frac{1}{p}|\langle[\Div(s_{p})],\varphi\rangle|\leq 
		C|\varphi|_{\mathscr{C}^{0}(X)}.
	\end{equation}
	By considering a countable $\mathscr{C}^{0}$-dense family of 
	$\varphi$'s in $\Omega^{(n-1,n-1)}_{0}(X)$, and applying Theorem 
	\ref{thm:2.15}, we get \eqref{eq:2.47}.
\end{proof}

\begin{Remark}
	The extra assumptions in the above corollary are necessary in our 
	approach to the conclusion \eqref{eq:2.47};  it is, however, an 
	interesting question whether these extra assumptions can actually be removed.
\end{Remark}

\begin{Remark}\label{rk:3.9a}
	For each $p\in \N_{>0}$, we can take a sequence of i.i.d.\ standard complex 
	Gaussian random variables $\eta^{p}=\{\eta^{p}_{j}\}_{j=1}^{d_{p}}$, and 
	assume that they are mutually independent for different $p.$ 
	We define the flat Gaussian random sections 
	\begin{equation}
		\psi^{S_{p}}_{\eta^{p}}=\sum_{j=1}^{d_{p}}\eta^{p}_{j}S^{p}_{j},
		\label{eq:2.46}
	\end{equation}
	where $S_{p}=\{S^{p}_{j}\}_{j=1}^{d_{p}}$ is an orthonormal 
	(Hilbert) basis of $H^{0}_{(2)}(X,L^{p})$. Then the statements in 
	Theorems \ref{thm:asym}, \ref{thm:2.15},
	Proposition \ref{prop:boundedasym} and Corollary \ref{cor:3.8a} 
	still hold true for the sequence of random sections 
	$\psi^{S_{p}}_{\eta^{p}}, p\geq 1$.
\end{Remark}

\subsection{Large deviation estimates and hole probability}\label{subsection3.2}
In this subsection, we study the large deviation estimates for random 
zeros in a given domain with respect to the high tensor powers as in 
\cite{SZZ}, \cite{DMS} and \cite{DLM:21}. In particular, we prove 
Theorems \ref{thm:6.2} and \ref{thm:6.3}.
A key intermediate result in the approach to the above theorems is the 
proposition as follows, whose proof is deferred to the next 
subsection.

\begin{Proposition}\label{prop:1.3.3}
	Let $U$ be a relatively compact open subset in $X$. For any $\delta>0$, 
	there exists $C_{U,\delta}>0$ such that for all $ p\gg 0$,
	\begin{equation}
		\P\Big(
		\int_{U}\Big|\log{\big|\psi^{S_{p}}_{\eta}(x)\big|_{h^{p}}}\Big| 
		\,\mathrm{dV}(x)\geq 
		\delta p \Big)\leq e^{-C_{U,\delta}p^{n+1}}.
		\label{eq:1.3.5}
	\end{equation}
\end{Proposition}

\begin{proof}[Proof of Theorem \ref{thm:6.2}]
The Poincar\'{e}-Lelong formula \eqref{eq:2.25} shows that 
	\begin{equation}
\frac{\sqrt{-1}}{\pi}\partial\overline{\partial}
\log{|\psi^{S_{p}}_{\eta}|_{h^{p}}}=[\Div(\psi^{S_{p}}_{\eta})]-pc_{1}(L,h)
	\label{eq:3.4.13}
	\end{equation}
as an identity of $(1,1)$-currents on $X.$	
Now fix $\varphi\in \Omega^{(n-1,n-1)}_{0}(U)$.
Then
\begin{equation}
\begin{split}
\Big(\frac{1}{p}[\Div(s_{p})],\varphi\Big)
-\int_{X}c_{1}(L,h)\wedge \varphi=\frac{\sqrt{-1}}{p\pi}\int_{X}
\log{|\psi^{S_{p}}_{\eta}|_{h^{p}}}\,\partial\overline{\partial}\varphi.
		\end{split}
		\label{eq:3.4.14}
	\end{equation}
Since $\varphi$ has a compact support in $U$, so has 
$\partial\overline{\partial}\varphi$. Set
	\begin{equation}
S_{\varphi}=\max_{x\in U}\left|\frac{\sqrt{-1}\partial
\overline{\partial}\varphi(x)}{\mathrm{dV}(x)}\right|.
		\label{eq:3.4.16}
	\end{equation}
	We can and we may assume that $S_{\varphi}>0$. Then
	\begin{equation}
	\begin{split}
	\left|\frac{\sqrt{-1}}{p\pi}\int_{X}\log{|\psi^{S_{p}}_{\eta}|_{h^{p}}}\,
		\partial\overline{\partial}\varphi\right|\leq 
		\frac{S_{\varphi}}{p\pi}\int_{U}
		\big|\log{|\psi^{S_{p}}_{\eta}(x)|_{h^{p}}}\big|\,\mathrm{dV}(x).
		\end{split}
		\label{eq:3.4.17}
	\end{equation}
Applying Proposition \ref{prop:1.3.3} to right-hand side of 
\eqref{eq:3.4.17} we get \eqref{eq:6.0.7}. 
\end{proof}

\begin{proof}[Proof of Theorem \ref{thm:6.3}]
Estimate \eqref{eq:1.6.14paris} is a direct consequence of 
\eqref{eq:6.1.10} by taking 
$\delta=n\Vol^{L}_{2n}(U)$. 
Hence, it is sufficient to prove \eqref{eq:6.1.10}. For this purpose, let 
$\chi_{U}$ denote the characteristic function of $U$ on $X$.
Let $\delta>0$ be arbitrary, and take $\psi_{1}$, $\psi_{2}\in 
\mathcal{C}^{\infty}_{0}(X,\R)$ such that
$0\leq \psi_{1}\leq\chi_{U}\leq \psi_{2}\leq 1,$ and
\begin{equation}
		\int_{X} \psi_{1}\frac{c_{1}(L,h_{L})^{n}}{n!}\geq 
		\Vol^{L}_{2n}(U)-\delta, \quad
		\int_{X} \psi_{2}\frac{c_{1}(L,h_{L})^{n}}{n!}\leq 
		\Vol^{L}_{2n}(U)+\delta\,.
	\label{eq:5.1.30paris}
\end{equation}
Note that the existence of such functions is guaranteed by the 
assumption that $\partial U$ has measure $0$ with respect to 
$\mathrm{dV}$, hence also to $\frac{1}{n!}c_{1}(L,h_{L})^{n}$.
For $j\in \{1,2\},$ set $\varphi_{j}=\frac{1}{(n-1)!}{\psi_{j}}c_{1}(L,h_{L})^{n-1}$. 
By applying Theorem \ref{thm:6.2} to $\varphi_{j}$ separately, we 
get exactly \eqref{eq:6.1.10}. 
\end{proof}

\subsection{Proof of Proposition \ref{prop:1.3.3}}\label{subsection3.3}
Let $U\subset X$ be a relatively compact open subset. For $s_{p}\in H^{0}(X, L^{p})$, we set
\begin{equation}
	\mathcal{M}^{U}_{p}(s_{p})=\sup_{x\in 
	U}|s_{p}(x)|_{h^{p}}<+\infty.
	\label{eq:5.2.1}
\end{equation}
Before proving Proposition \ref{prop:1.3.3}, we need to investigate 
the probabilities for both, $\mathcal{M}^{U}_{p}(\psi^{S_{p}}_{\eta})$ taking atypically large and small values, respectively.

\begin{Proposition}\label{prop:3.14ss}
	For any $\delta>0$, there exists a constant $C_{U,\delta}>0$ such 
	that for $p\in\N_{>1}$, 
	\begin{equation}
		\P\big (\mathcal{M}^{U}_{p}(\psi^{S_{p}}_{\eta})\geq 
		e^{\delta p}\big)\leq e^{-\delta 
		p^{n+1}+C_{U,\delta}p^{n}\log{p}}\,.
		\label{eq:5.2.2}
	\end{equation}
\end{Proposition}
\begin{proof}
The basic idea of the proof is that the local $\cLL$-norm of a holomorphic 
function is bounded by its local sup-norm as in \eqref{eq:2.2a}. 
	We fix $\delta>0$ and let $r>0$ be sufficiently small so that we can 
	choose a finite set of points $\{x_{j}\}_{j=1}^{\ell}\subset U$ such 
	that the geodesic open balls $B^{X}(x_{j},r)$, $j=1,\ldots, \ell$ 
	form an open covering of $\overline{U}$. Since $r$ is 
	sufficiently small, then we can assume that each larger ball
	$B^{X}(x_{j},2r)$ lies in a complex chart (hence viewed as an 
	open subset of $\C^{n}$), and that for each $j$, we can fix a 
	local holomorphic frame $e_{L,j}$ of $L$ on a neighborhood of $B^{X}_{x_{j}, 2r}$ 
	with $\sup_{x\in B^{X}(x_{j},2r)} |e_{L,j}(x)|_{h_{L}}=1$.
	Set
	\begin{equation}
		\nu=\min\big\{ \inf_{x\in 
		B^{X}(x_{j},2r)}|e_{L,j}(x)|_{h_{L}} \,:\,j=1 , \ldots, \ell\big\}.
		\label{eq:3.14}
	\end{equation}
	It is clear that $0<\nu\leq 1$. By fixing $r$ small enough, we 
	can and do assume that 
	 \begin{equation} \label{eq:nuBd}
	 -\log{\nu}\leq 
	\frac{\delta}{6}\,\cdot
	\end{equation} 	
		As in \eqref{eq:2.2a}, since $U$ is relatively compact, there 
	exists a constant $C>0$ such that for each $j=1,\ldots, \ell$, if 
	$f$ is a holomorphic function on a neighborhood of 
	$B^{X}(x_{j},2r)$, then
	\begin{equation}
		\sup_{x\in B^{X}(x_{j},r)}|f(x)|\leq 
		C\|f\|_{\mathcal{L}^{2}(B^{X}(x_{j},2r))},
		\label{eq:3.16}
	\end{equation}
	where the volume form $\mathrm{dV}(x)$ on $X$ is used in the norm 
	$\|\cdot\|_{\mathcal{L}^{2}(B^{X}(x_{j},2r))}$.
	Note that the choices of $x_{j}$, $r$, $\ell$, and the constants 
	$\nu$, $C$
	are independent of the tensor power $p$. Set 
		$\widetilde{U}=\cup_{j} B^{X}(x_{j},2r)\supset U$.
	For $p\in\N, s_{p}\in H^{0}(X,L^{p})$, on each $B^{X}(x_{j},2r)$, we write
	\begin{equation}
		s_{p}|_{B^{X}(x_{j},2r)}=f_{j}e_{L,j}^{\otimes p},
		\label{eq:3.15}
	\end{equation}
	where $f_{j}$ is a holomorphic function on the chart in $\C^{n}$ 
	corresponds to $B^{X}(x_{j},2r)$. Then we have
	\begin{equation}
		\begin{split}
			\mathcal{M}^{U}_{p}(s_{p})=\sup_{x\in U}|s_{p}(x)|_{h^{p}}&\leq \max_{j} \sup_{x\in 
			B^{X}(x_{j},r)}|f_{j}(x)|\\
			&\leq C\max_{j} 
			\{\|f_{j}\|_{\mathcal{L}^{2}(B^{X}(x_{j},2r))}\}\\
			&\leq 
			\frac{C}{\nu^{p}}\max_{j}\{\|s_{p}\|_{\mathcal{L}^{2}(B^{X}(x_{j},2r),L^{p})}\}\\
			&\leq 
			\frac{C}{\nu^{p}}\|s_{p}\|_{\mathcal{L}^{2}(\widetilde{U},L^{p})}.
		\end{split}
		\label{eq:3.17bis}
	\end{equation}

	The next step is to estimate the quantity 
	$\E[\|\psi^{S_{p}}_{\eta}\|^{2p^{n}}_{\mathcal{L}^{2}(\widetilde{U},L^{p})}]$ for $p\geq 2$. Applying H\"{o}lder's inequality with $\frac{1}{p^{n}}+\frac{p^{n}-1}{p^{n}}=1$, we get 
	\begin{equation}
		\E\big[\|\psi^{S_{p}}_{\eta}\|^{2p^{n}}_{\mathcal{L}^{2}(\widetilde{U},L^{p})}\big]\leq \Vol(\widetilde{U})^{p^{n}-1}\E\Big[\int_{\widetilde{U}}|\psi^{S_{p}}_{\eta}(x)|^{2p^{n}}_{h^{p}}(x)\mathrm{dV}\Big].
		\label{eq:holderpn}
	\end{equation}
	As in \eqref{eq:3.16}, on a neighborhood of $B^{X}(x_{j}, 2r)$, write
	\begin{equation}
		S^{p}_{i}=f^{p}_{i}e_{L,i}^{\otimes p}.
		\label{eq:3.18bis}
	\end{equation}
	If $x\in B^{X}(x_{j},2r)$, set
	\begin{equation}
		F_{j}(x)=\sum_{i=1}^{d_{p}} \eta_{i} f^{p}_{i}(x).
		\label{eq:3.19bis}
	\end{equation}
	Then $F_{j}(x)$ is a complex Gaussian random variable with 
	(total) variance $\sum_{i=1}^{d_{p}}|f^{p}_{i}(x)|^{2}$. 
	By our 
	assumption on the local frame $e_{L,j}$, we get
	\begin{equation}
		\sum_{i=1}^{d_{p}}|f^{p}_{i}(x)|^{2}\leq \frac{1}{\nu^{2p}}P_{p}(x,x).
		\label{eq:3.20bis}
	\end{equation}
	Then we 
	have
	\begin{equation}
		\E\big[|F_{j}(x)|^{2p^{n}}\big]=p^{n}!\Big(\sum_{i=1}^{d_{p}}|f^{p}_{i}(x)|^{2}\Big)^{p^{n}}.
		\label{eq:3.21}
	\end{equation}
	As a consequence, we get that for $x\in \widetilde{U}$,
	\begin{equation}
		\begin{split}
				\E\big[|\psi^{S_{p}}_{\eta}(x)|_{h^{p}}^{2p^{n}}\big]&\leq 
				\frac{1}{\nu^{2p^{n+1}}}\E\big[|F_{j}(x)|^{2p^{n}}\big]
				\leq 
				\frac{p^{n}!}{\nu^{4p^{n+1}}}(P_{p}(x,x))^{p^{n}}.
		\end{split}
	\end{equation}

	Since we are in the context of $\sigma$-finite measures and 
	the integrands are non-negative, Tonelli's Theorem applies, so that
	\begin{equation}
		\E\Big[\int_{\widetilde{U}}|\psi^{S_{p}}_{\eta}(x)|^{2p^{n}}_{h^{p}}\mathrm{dV}(x)\Big]\leq \frac{p^{n}!}{\nu^{4p^{n+1}}}\int_{\widetilde{U}} (P_{p}(x,x))^{p^{n}}\mathrm{dV}(x).
		\label{eq:3.24}
	\end{equation}
Moreover, by the on-diagonal estimate for the Bergman kernel on a given compact 
subset, there exists a 
	constant $C_{\widetilde{U}}>0$ (independent of $p$) such that for 
	$p\in\N$, $x\in\widetilde{U}$,
	\begin{equation}
		P_{p}(x,x)\leq C_{\widetilde{U}}p^{n}.
		\label{eq:3.25}
	\end{equation}
		Combining \eqref{eq:holderpn} with the above inequalities, we infer that
	\begin{equation}
		\E\Big[\|\psi^{S_{p}}_{\eta}\|^{2p^{n}}_{\mathcal{L}^{2}(\widetilde{U},L^{p})} \Big]\leq \big(C_{\widetilde{U}}\Vol(\widetilde{U})\big)^{p^{n}}\frac{p^{n}!}{\nu^{4p^{n+1}}}(p^{n})^{p^{n}}.
		\label{eq:3.26}
	\end{equation}
	By applying \eqref{eq:3.17bis} to $\psi^{S_{p}}_{\eta}$, we get
	\begin{equation}
		\begin{split}
			\E\big[\mathcal{M}^{U}_{p}(\psi^{S_{p}}_{\eta})^{2p^{n}}\big]&\leq 
		\Big(\frac{C}{\nu^{p}}\Big)^{2p^{n}}\E\big[\|\psi^{S_{p}}_{\eta}\|^{2p^{n}}_{\mathcal{L}^{2}(\widetilde{U},L^{p})}\big]
		\leq \frac{(\widetilde{C}p^{n})^{2p^{n}}}{\nu^{6p^{n+1}}},
		\end{split}
		\label{eq:3.18}
	\end{equation}
	where $C>0$, $\widetilde{C}>0$ are constants independent of $p$.

	Then \eqref{eq:5.2.2} follows from Chebyshev's inequality and the inequality
	$\frac{1}{\nu}\leq e^{\frac{\delta}{6}}$ from \eqref{eq:nuBd}.
\end{proof}

\begin{Remark}
The choice to consider the $p^{n}$-th moment of 
$\|\psi^{S_{p}}_{\eta}\|^{2}$ leads to the exponent $p^{n+1}$ in the exponential of the resulting
probability estimate. One can consider arbitrary $N$-th moments to obtain a more
general statement on this probability upper bound. 	

When $X$ is compact, or if $X$ is noncompact but $d_{p}$ is 
bounded polynomially in $p$, then the upper bound 
$Ce^{-cp^{n+1}}$ can be obtained in a simpler way as in \cite{SZZ} 
and in \cite{DLM:21} (and of course with a much sharper upper bound). 
\end{Remark}

Now we consider the probabilities of small values of 
$\mathcal{M}^{U}_{p}(\psi^{S_{p}}_{\eta})$, and we will adapt the 
ideas in \cite{SZZ}, 
\cite{DLM:21}. At first, we introduce a result 
on the near-diagonal estimate of Bergman kernel.

Recall that $\dot{R}^{L}$ is defined in \eqref{eq:3.3a}. Now fix a 
point $x\in X$. Let $\{\mathbf{f}_{j}\}_{j=1}^{n}$ be an orthonormal basis of 
$(T_{x}^{1,0}X, g_{x}^{TX}(\cdot,\overline{\cdot}))$ such that
\begin{equation}
	\dot{R}^{L}_{x}\,\mathbf{f}_{j}=\mu_{j}(x)\mathbf{f}_{j},
\end{equation}
where $\mu_{j}(x)$, $j=1, \ldots, n$, are the eigenvalues of 
$\dot{R}^{L}_{x}$. Then by the first inequality in \eqref{eq:intro3.0.1}, we have
\begin{equation}
	\mu_{j}(x)\geq\varepsilon.
\end{equation}
Set 
$\mathbf{e}_{2j-1}=\frac{1}{\sqrt{2}}(\mathbf{f}_{j}+\overline{\mathbf{f}}_{j})$, 
$\mathbf{e}_{2j}=\frac{\sqrt{-1}}{\sqrt{2}}(\mathbf{f}_{j}-\overline{\mathbf{f}}_{j})$, $j=1, \ldots, n$. Then they form an orthonormal basis of the (real) tangent vector space $(T_{x}X, g_{x}^{TX})$. 
If $v=\sum_{j=1}^{2n} v_{j}\mathbf{e}_{j}\in T_{x}X$, we can write
\begin{equation}
	v=\sum_{j=1}^{n}(v_{2j-1}+\sqrt{-1}v_{2j})\frac{1}{\sqrt{2}} 
	\mathbf{f}_{j} + \sum_{j=1}^{n}(v_{2j-1}-\sqrt{-1}v_{2j})\frac{1}{\sqrt{2}} 
	\overline{\mathbf{f}}_{j}.
	\label{eq:5.1.8paris}
\end{equation}
Set $z=(z_{1},\ldots,z_{n})$ with $z_{j}=v_{2j-1}+\sqrt{-1}v_{2j}$, 
$j=1,\ldots, n$. We call $z$ the complex coordinate of $v\in T_{x}X$. 
Then by \eqref{eq:5.1.8paris},
\begin{equation}
	\frac{\partial}{\partial z_{j}} = \frac{1}{\sqrt{2}} 
	\mathbf{f}_{j} , \; \frac{\partial}{\partial \overline{z}_{j}} = \frac{1}{\sqrt{2}} 
	\overline{\mathbf{f}}_{j},
	\label{eq:5.1.9paris}
\end{equation}
so that
\begin{equation}
	v=\sum_{j=1}^{m}\Big(z_{j} \frac{\partial}{\partial z_{j}} 
	+\overline{z}_{j} \frac{\partial}{\partial \overline{z}_{j}}\Big).
	\label{eq:5.1.10paris}
\end{equation}
Note that $|\frac{\partial}{\partial z_{j}}|^{2}_{g^{TX}}=|\frac{\partial}{\partial 
\overline{z}_{j}}|^{2}_{g^{TX}}=\frac{1}{2}$.
For $v, v'\in T_{x}X$, let $z, z'$ denote the corresponding 
complex coordinates. 
	
Define a weighted distance function $\Phi^{TX}_{x}(v,v')$ as follows,
\begin{equation}
\Phi^{TX}_{x}(v,v')^{2}=\sum_{j=1}^{n}\mu_{j}(x)|z_{j}-z'_{j}|^{2}.
\end{equation}

For sufficiently small $\delta_{0}>0$, we identify the small open 
ball
$B^{X}(x,2\delta_{0})$ in $X$ with the ball $B^{T_{x}X}(0, 
2\delta_{0})$ in $T_{x}X$ via the geodesic coordinate. 	
Let $\mathrm{dist}(\cdot,\cdot)$ denote the Riemannian distance 
of $(X, g^{TX})$. There exists $C_{1}>0$ such that for $v,v'\in 
B^{T_{x}X}(0, 2\delta_{0})$, we 
have
\begin{equation}
	C_{1}\mathrm{dist}(\exp_{x}(v),\exp_{x}(v'))\geq \Phi^{TX}_{x}(v,v')\geq 
	\frac{1}{C_{1}}\mathrm{dist}(\exp_{x}(v),\exp_{x}(v')).
	\label{eq:5.1.15DLM}
\end{equation}
In particular,
\begin{equation}
	\Phi^{TX}_{x}(0,v)\geq \varepsilon^{1/2} 
	\mathrm{dist}(x,\exp_{x}(v)).
	\label{eq:5.1.15paris}
\end{equation}
Moreover, if we consider a compact subset $K\subset X$, the 
constants $\delta_{0}$ and $C_{1}$ can be chosen uniformly for all 
$x\in K$.

For $p\in \N$, $x,y\in X$, the normalized Bergman kernel is defined as
\begin{equation}
	N_{p}(x,y)=\frac{|P_{p}(x,y)|_{h^{p}_{x}\otimes 
	h_{y}^{p,\ast}}}{\sqrt{P_{p}(x,x)}\sqrt{P_{p}(y,y)}}.
	\label{eq:1.3.1}
\end{equation}

The following result is proved in \cite[Theorem 5.1]{DLM:21}, where 
we use essentially the near-diagonal expansion of Bergman kernel in \cite[Theorems 4.2.1 
\& 6.1.1]{MM07}.
\begin{Theorem}\label{thm:5.1.1}
	Let $U$ be a relatively compact open subset of $X.$ Then the 
	following uniform estimates on the normalized Bergman kernel hold 
	for $x,y\in U$:
	For $k\geq 1$ and $b>\sqrt{16k/\varepsilon}$ fixed, we have for $p\gg 0\,$  
	(such that $b\sqrt{\frac{\log 
		p}{p}}\leq 2\delta_{0}$) that
	\begin{equation}
		N_{p}(x,y)=\begin{cases} 
		&\big(1+o(1)\big)\exp\Big(-\dfrac{p}{4}\,\Phi_{x}(0,v')^{2}\Big),  \\
		&\qquad\qquad\text{uniformly for } 
		\mathrm{dist}(x,y)\leq b\sqrt{\frac{\log 
		p}{p}}\,,\;\text{with\;}y=\exp_{x}(v'), v'\in T_{x}X; \\
		&\mathcal{O}(p^{-k}), \;\;   \text{uniformly for } 
		\mathrm{dist}(x,y) \geq 
		b\sqrt{\frac{\log p}{p}}.
	\end{cases}
	\label{eq:5.2.4}
\end{equation}
\end{Theorem}

\begin{Proposition}\label{prop:3.17}
	There exist constants $C_{U}>0, C_{U}'>0$ such 
	that for all $\delta>0$ and $p\in\N$, 
	\begin{equation}
		\P\big(\mathcal{M}^{U}_{p}(\psi^{S_{p}}_{\eta})\leq 
		e^{-\delta p}\big)\leq e^{-C_{U}\delta p^{n+1}+C_{U}' 
		p^{n}\log{p}}\,.
		\label{eq:5.2.3}
	\end{equation}
\end{Proposition}
\begin{proof}
For $x\in X$ we fix some $\lambda_{x}\in L_{x}$ with 
 $|\lambda_{x}|_{h}=1$, and set
 \begin{equation}
 \xi_{x}=\frac{\langle \lambda_{x}^{\otimes p}, 
 \psi^{S_{p}}_{\eta}(x)\rangle_{h^{p}}}{\sqrt{P_{p}(x,x)}}\,\cdot
\label{eq:53.2.15}
 \end{equation}
Then $\xi_{x}$ is a complex Gaussian random variable. Moreover, for 
any two points $x,y\in X$, we have 
 \begin{equation}
	\big|\mathbb{E}[\xi_{x}\overline{\xi}_{y}]\big|
=N_{p}(x,y).
\label{eq:53.2.19}
 \end{equation}
	
	Then by the asymptotic equations in \eqref{eq:5.2.4}, using the similar arguments in \cite[Subsection 
	3.2]{SZZ} or the proof of \cite[Theorem 1.13]{DLM:21}, we can 
	prove a more general version of \eqref{eq:5.2.3} as follows: for a sequence of positive numbers 
$\{\lambda_{p}\}_{p\in\N}$,
\begin{equation}
	\P\big(\mathcal{M}^{U}_{p}(\psi^{S_{p}}_{\eta})\leq 
	\lambda_{p}\big)\leq e^{Cp^{n}\log \lambda_{p} +C'p^{n}\log 
	p},\quad  p\gg 0.
	\label{eq:53.2.16}
\end{equation}
Then, for any 
$\delta>0$, choosing $\lambda_{p}=e^{-\delta p}$ in \eqref{eq:53.2.16}, we recover 
\eqref{eq:5.2.3}. This completes our 
proof.
\end{proof}
Combining Propositions \ref{prop:3.14ss} and \ref{prop:3.17}, we arrive at the following.
\begin{Corollary}\label{cor:3.18ss}
	For any relatively compact open subset $U\subset X$, and for 
	$\delta>0$, there exists a constant $C=C(U,\delta)>0$ such that 
	for $p\gg 1$,
	\begin{equation}
		\P\left(\left|\log{\mathcal{M}^{U}_{p}(\psi^{S_{p}}_{\eta})}\right|\geq \delta p\right) \leq e^{-Cp^{n+1}}.
		\label{eq:3.79ss}
	\end{equation}
\end{Corollary}

\begin{proof}[Proof of Proposition \ref{prop:1.3.3}]
	The proof of Proposition \ref{prop:1.3.3} 
follows by combining from the arguments in \cite[Subsection 4.1]{SZZ}
with Corollary \ref{cor:3.18ss}. Here, we just sketch the proof.

For $t>0$, set
\begin{equation}
	\log^{+}{t}=\max\{\log{t},0\},\; 
	\log^{-}{t}:=\log^{+}(1/t)=\max\{-\log{t},0\}.
	\label{eq:3.3.4}
\end{equation}
Then
\begin{equation}
	|\log{t}|=\log^{+}{t}+\log^{-}{t}.
	\label{eq:3.3.5}
\end{equation}

Let $U$ be a relatively compact nonempty open subset in $X.$ Then for any 
nonzero holomorphic section $s_{p}\in H^{0}(X,L^{p})$, we have that $\big|\log|s_{p}|_{h^{p}}\big|$ is integrable on 
$\overline{U}$ with respect to $\mathrm{dV}$. We now start with showing that
\begin{equation}
	\P\left(
	\int_{U}\log^{+}{|\psi^{S_{p}}_{\eta}(x)|_{h^{p}}}\,\mathrm{dV}(x)\geq 
	\frac{\delta}{2} p\right)\leq e^{-C_{U,\delta}p^{n+1}}.
	\label{eq:3.3.6}
\end{equation}
For this purpose, observe that on $U$ we have
\begin{equation}
	\log^{+}{|\psi^{S_{p}}_{\eta}|_{h^{p}}}\leq 
	\big|\log\mathcal{M}^{U}_{p}(\psi^{S_{p}}_{\eta})\big|,
	\label{eq:3.3.7}
\end{equation}
which then supplies us with
\begin{equation}
	\begin{split}
		&\P\left(
		\int_{U}\log^{+}{|\psi^{S_{p}}_{\eta}(x)|_{h^{p}}}\,\mathrm{dV}(x)\geq 
	\frac{\delta}{2} p\right)\\
		&\leq \P\left(
		\left|\log\mathcal{M}^{U}_{p}(\psi^{S_{p}}_{\eta})\right|\geq 
		\frac{\delta}{2 \Vol(U)} p\right),
	\end{split}
	\label{eq:3.3.9}
\end{equation}
where $\Vol(U)$ denotes the volume of $U$ with respect to 
$\mathrm{dV}$. In combination with Corollary \ref{cor:3.18ss}, this immediately implies \eqref{eq:3.3.6}. 

The next step is to prove that 
\begin{equation}
	\P\left(
	\int_{U}\log^{-}{|\psi^{S_{p}}_{\eta}(x)|_{h^{p}}}\,\mathrm{dV}(x)\geq 
	\frac{\delta}{2} p\right)\leq e^{-C_{U,\delta}\,p^{n+1}}.
	\label{eq:3.3.10}
\end{equation}
Suppose that $U$ contains an annulus $B(2,3):=\{z\in\C^{n}\;:\;2< 
|z|<3\}$ (possibly after rescaling of coordinates), and the line bundle $L$ 
on $B(1,4)$ (still contained in $U$) has a holomorphic local frame 
$e_{L}$. Set $\alpha(x)= 
\log|e_{L}(x)|^{2}_{h}$.  We can then write
\begin{equation}
	\psi^{S_{p}}_{\eta}=F_{p}e^{\otimes p}_{L},
	\label{eq:3.4.11DLM}
\end{equation}
where $F_{p}$ is a random holomorphic function on $B(1,4)$. Then
\begin{equation}
	\log|\psi^{S_{p}}_{\eta}|_{h^{p}}=\log|F_{p}|+\frac{p}{2}\alpha.
	\label{eq:3.4.12DLM}
\end{equation}
In the following estimates, each $K_{i},$ $i \in \N,$ denotes a sufficiently 
large positive constant. Then by \eqref{eq:3.3.5} and \eqref{eq:3.3.9}, we have
\begin{equation}
	\P\left(
	\int_{B(2,3)}\log^{+}{|F_{p}|}\,\mathrm{dV}\geq 
	K_{1} p  \right)\leq e^{-C_{U,K_{1}}\, p^{n+1}}.
	\label{eq:3.4.13DLM}
\end{equation}
Using the Poisson kernel and the sub-mean inequality for 
$\log(|F_{p}|)$, we can improve \eqref{eq:3.4.13DLM} to get
\begin{equation}
	\P\left(
	\int_{B(2,3)}\log{|F_{p}|}\,\mathrm{dV}\geq 
	K_{2} p  \right)\leq e^{-C_{U,K_{2}}\, p^{n+1}}.
	\label{eq:3.4.14DLM}
\end{equation}
From this point we proceed as in \cite[Subsection 
4.1, pp.\ 1992]{SZZ}. For $\delta\in\; 
]0,\frac{1}{2}]$, we get a finite set of (almost uniformly distributed) points 
$\{z_{j}\}_{j=1}^{q}$ in $B(2,3)$ such that for all $s_{p}\in 
H^{0}(X,L^{p})$, $p\in\N$, $s_{p}=f_{p}e_{L}^{\otimes p}$ on 
$B(1,4)$,
we have
\begin{equation}
	\begin{split}
		&-\int_{B(2,3)} \log|s_{p}|_{h^{p}}\mathrm{dV}\\
		&\leq 
		-\sum_{j=1}^{q} \mu_{j} 
		\log|s_{p}|_{h^{p}}(z_{j})+K_{3}\delta\int_{B(2,3)} 
		\big|\log{|f_{p}|}\big|\,\mathrm{dV}+p\delta K_{3} \sup_{z\in 
		B(2,3)}|d\alpha(z)|_{g^{T^{\ast}X}},
	\end{split}
\end{equation}
where the quantities $q$ and $\mu_{j}>0$ only depend on $\delta$, and 
we have $\sum_{j=1}^{q}\mu_{j}\simeq 1$. Note that the constant 
$K_{3}$ does not depend on $\delta$.
Applying the above inequality to $\psi^{S_{p}}_{\eta}$ and $F_{p}$, using 
Corollary \ref{cor:3.18ss} for each term 
$\log|\psi^{S_{p}}_{\eta}|_{h^{p}}(z_{j}),$ and taking advantage of \eqref{eq:3.4.13DLM}, we infer that
\begin{equation}
	\P\left(
	-\int_{B(2,3)}\log{|\psi^{S_{p}}_{\eta}|_{h^{p}}}\,\mathrm{dV}\geq 
	K_{4}\delta p\right)\leq e^{-C_{U,\delta}p^{n+1}},\;\forall\;p\gg 0.
	\label{eq:3.4.16DLM}
\end{equation}

Noting that $\log^{-}=-\log+\log^{+}$ and that a finite set of annuli of the form 
$B(2,3)$ covers $U$, we can infer \eqref{eq:3.3.10} from 
\eqref{eq:3.3.9} and
\eqref{eq:3.4.16DLM}. This completes our proof.
\end{proof}

\begin{Remark}\label{rm:asympnormality}
With results for the regimes of the law of large numbers as well as of large deviations at our disposal, a naturally ensuing question is that of central limit type behavior. In fact, the asymptotic normality of (functionals of) the zeros of random holomorphic functions 
or sections has been introduced and proved by Sodin-Tirelson 
\cite[Main Theorem]{STr} for certain random holomorphic functions on 
$\C$ or $\mathbb{D}$ and by Shiffman-Zelditch 
\cite[Theorem 1.2]{MR2742043} for the random holomorphic sections of 
line bundles on a compact K\"{a}hler manifold. An extension to 
general random polynomials on $\C^{n}$ was also proved by Bayraktar  \cite{B9}. One key ingredient in 
their approaches is the normalized Bergman kernel defined in 
\eqref{eq:1.3.1} viewed as 
the covariance function of a normalized Gaussian process on $\C$ or $X$, as 
constructed in the proof of Proposition \ref{prop:3.17}. Then using the 
estimates given in Theorem \ref{thm:5.1.1} and the seminal result proved by Sodin and Tirelson 
in \cite[Theorem 2.2]{STr}, one could obtain an 
extension of \cite[Main Theorem]{STr} \cite[Theorem 1.2]{MR2742043} 
to our noncompact setting. 
\end{Remark}

\subsection{Remark on the lower bound for the hole 
probabilities}\label{ss3.4ss}

To obtain a lower bound of matching order $e^{-cp^{n+1}}$ for the hole probability in 
\eqref{eq:1.6.14paris} is generally
more complicated. For the case of scaled Bargmann-Fock spaces 
(cf.\ Example \ref{ex:3.6}), we can provide a lower bound and we sketch its proof in the sequel.

Recall that for any $p\in\N$, the family
$\{S^{p}_{\alpha}\}_{\alpha\in \N^{n}}$ denotes an orthonormal basis of 
$H^{0}_{(2)}(\C^{n},L^{p})$. For $K>0$, define the index set
\begin{equation}
	I(K)=\big \{\alpha=(\alpha_{1},\ldots,\alpha_{n})\in\N^{n}\;:\; 0\leq 
	\alpha_{j}\leq K, j=1,\ldots, n\big\},
\end{equation}
 set $I^{\ast}(K)=I(K)\backslash \{(0,\ldots, 0)\},$ and put
\begin{equation}
	q_{p}:=\sharp I(Kp)=(1+\lfloor Kp\rfloor)^{n}=\mathcal{O}(p^{n}).
	\label{eq:3.4.2ss}
\end{equation}

For this canonical family of orthonormal bases, we 
can verify directly the following local concentration condition: for 
any relatively compact subset $U\subset \C^{n}$ and for any 
$c>0$, there exist constants $K=K(U,c)>0, C'=C'(U,c)>0$ such that
\begin{equation}
	\sup_{z\in\bar{U}} \sum_{\alpha \notin I(Kp)}|S^{p}_{\alpha}(z)|^{2}_{h^{p}}\leq C' e^{-cp}.
	\label{eq:3.92ss}
\end{equation}
Let $\psi^{S_{p}}_{\eta}$ be the random holomorphic section 
(actually, function) on $\C^{n}$ 
constructed in Example \ref{ex:3.6}.

\begin{Lemma}\label{lem:3.17lower}
	For any relatively compact open subset $U\subset \C^{n}$, there exists a constant $C'_{U}>0$ such that for $p\gg 
	1$,
	 	\begin{equation}
 		\P\left(\Div(\psi^{S_{p}}_{\eta})\cap 
		U=\varnothing\right)\geq 
 		e^{-C'_{U}p^{n+1}}.
 		\label{eq:1.2.4}
 	\end{equation}
\end{Lemma}
\begin{proof}
	For $U=\varnothing$ the statement is trivial, so assume
	 $U$  nonempty is as in the assumptions. 
	Fix a relatively compact open neighborhood $U'$ of $\overline{U}$ and
	 define the strictly positive quantity 
	\begin{equation}
		\widetilde{M}:=\min_{z\in \overline{U'}} e^{-\frac{|z|^{2}}{2}} \in (0, 1).
	\end{equation}
	Let the 
	constants $K$ and $C'$ be the ones in \eqref{eq:3.92ss} for the 
	constant $c=-2\log{\widetilde{M}}>0$ and for $U'$.
	For $p\in \N$, write $S^{p}_{0}\equiv p^{n/2}$ for the unit constant
	section in $H^{0}_{(2)}(\C^{n},L^{p})$ corresponding to 
	$\alpha=(0,\ldots,0)\in\N^{n}$. Then
	\begin{equation}
		\min_{z\in \overline{U'}} 
		|S^{p}_{0}(z)|_{h^{p}}=p^{n/2}\widetilde{M}^{p}.
	\end{equation}
	Defining the random holomorphic sections
	\begin{equation}
		\begin{split}
		&\psi^{S_{p}}_{\eta,\mathrm{I}}(z):=\sum_{ 
		\alpha \in I^{*}(Kp)}\eta_{\alpha} S^{p}_{\alpha}(z) \quad \text{and} \\
		&\psi^{S_{p}}_{\eta,\mathrm{II}}(z):=\sum_{\alpha \notin I(Kp)}\eta_{\alpha} 
		S^{p}_{\alpha}(z).
		\end{split}
	\end{equation}
	we can decompose
	\begin{equation}
		\psi^{S_{p}}_{\eta}=\eta_{0} S^{p}_{0} + 
		\psi^{S_{p}}_{\eta,\mathrm{I}} + 
		\psi^{S_{p}}_{\eta,\mathrm{II}}.
		\label{eq:3.97ss}
	\end{equation}
	Note that the three random sections on the right-hand side of
	\eqref{eq:3.97ss} are independent from each other.
	
	In the remaining part of the proof, we view the above sections as 
	holomorphic functions on $\C^{n}$, and let $|\cdot|$ denote the 
	standard modulo on $\C$ (instead of considering the norm $|\cdot|_{h^{p}}$ 
	on line bundle).
	
	Applying \eqref{eq:2.2a} to the function $\psi^{S_{p}}_{\eta,\mathrm{II}}$ and using the estimate  \eqref{eq:3.92ss}, we 
	arrive at the upper bound
	\begin{equation}
		\begin{split} \label{eq:secMomII}
					\mathbb{E}\Big[\sup_{z\in 
		U}\Big|\psi^{S_{p}}_{\eta,\mathrm{II}}(z)\Big|^{2}\Big]&\leq 
		C_{U'}\sigma^{2}\int_{U'}\sum_{\alpha \notin I(Kp)}|S^{p}_{\alpha}(z)|^{2}\mathrm{dV}(z)\\
		&\leq \widetilde{C}_{U'}\Vol(U')\sigma^{2} 
		\widetilde{M}^{-2p} e^{-cp}\\
		&= \widetilde{C}_{U'}\Vol(U')\sigma^{2} =:\widetilde{C}',
		\end{split}
	\end{equation}
	where the last equality follows from our choice 
	$c=-2\log{\widetilde{M}}$. 
	For any $\lambda>0,$ as a consequence of Chebyshev's inequality in combination with \eqref{eq:secMomII},  we have
	\begin{equation}
		\P\Big( \sup_{z\in 
		U}\Big|\psi^{S_{p}}_{\eta,\mathrm{II}}(z)\Big| < 
		\lambda\Big)\geq 
		1-\frac{\widetilde{C}'}{\lambda^{2}}.
	\end{equation}
	We define the good event
	\begin{equation}
		\Omega_{p}=\left\{|\eta_{0}|\geq 1\,;\, |\eta_{\alpha}|\leq 
		\frac{1}{3\sqrt{q_{p}-1}}\widetilde{M}^{p}, \alpha \in I^{\ast}(Kp)\,;\, \sup_{z\in 
		U}\left|\psi^{S_{p}}_{\eta,\mathrm{II}}(z)\right| < 
		\frac{1}{3}p^{n/2} \right\}.
	\end{equation}
	For all sufficiently large $p\in\N$, we have
	\begin{equation}
		\begin{split}
			\P(\Omega_{p})&=\P(|\eta_{0}|\geq 1) \cdot \P\left(\sup_{z\in 
		U}\left|\psi^{S_{p}}_{\eta,\mathrm{II}}(z)\right| < 
		\frac{1}{3}p^{n/2} 
		\right)\\
		&\qquad\cdot\P\Bigg(|\eta_{\alpha}|\leq 
		\frac{1}{3\sqrt{q_{p}-1}}\widetilde{M}^{p}, \alpha \in 
		I^{\ast}(Kp)\Bigg)\\
		&\geq e^{-1}\Big(1-\frac{9\widetilde{C}'}{
	p^{n}}\Big)\cdot \Pi_{\alpha \in 
	I^{\ast}(Kp)}\left(\frac{1}{18(q_{p}-1)}\widetilde{M}^{2p}\right).
		\end{split}
		\label{eq:3.102ss}
	\end{equation}
	Then by \eqref{eq:3.4.2ss}, there exists $C'_{U}>0$ such that for 
	$p\gg 1$,
	\begin{equation}
		\P(\Omega_{p})\geq e^{-C'_{U}p^{n+1}}.
	\end{equation}
	Our lemma then follows once we show the inclusion
	\begin{equation}
		\Omega_{p}\subset \left\{\Div(\psi^{S_{p}}_{\eta})\cap 
		U=\varnothing\right\}.
		\label{eq:3.105ss}
	\end{equation}	
	Indeed, if $|\eta_{\alpha}|\leq \frac{1}{3\sqrt{q_{p}-1}}\widetilde{M}^{p}$, 
	$\alpha \in I^{\ast}(Kp)$, then for $z\in U$,
	\begin{equation}
		\begin{split}
				\big|\psi^{S_{p}}_{\eta,\mathrm{I}}(z)\big|^{2}&\leq 
		\Big(\sum_{\alpha\in 
		I^{\ast}(pK)}|\eta_{\alpha}|^{2}\Big)\Big(\sum_{\alpha\in 
		I^{\ast}(pK)}|S^{p}_{\alpha}(z)|^{2}\Big)\\
		&\leq \frac{1}{\widetilde{M}^{2p}}\Big(\sum_{\alpha\in 
		I^{\ast}(pK)}|\eta_{\alpha}|^{2}\Big)P_{p}(z,z)\\
		&\leq \frac{1}{9}p^{n}.
		\end{split}
	\end{equation}
	As a consequence, on $\Omega_{p}$ and for $z\in U$, we get
	\begin{equation}
		\begin{split}
			\big|\psi^{S_{p}}_{\eta,\mathrm{I}}(z) + 
		\psi^{S_{p}}_{\eta,\mathrm{II}}(z)\big|&\leq 
		\big|\psi^{S_{p}}_{\eta,\mathrm{I}}(z)\big|+
		\big|\psi^{S_{p}}_{\eta,\mathrm{II}}(z)\big|
		\leq \frac{1}{3}p^{n/2}+\frac{1}{3}p^{n/2}\\
		&< p^{n/2}\leq \left|\eta_{0}S^{p}_{0}(z)\right|.
		\end{split}
	\end{equation}
	The above strict inequality implies that \eqref{eq:3.105ss} is fulfilled. 
	This finishes the  proof of the lemma.
\end{proof}

We now shortly explain how by applying our results to the 
special case of the Bargmann-Fock space recovers the results by 
Sodin-Tsirelson 
	(for $\C$, \cite[Theorem 1]{STr05}) 
	and Zrebiec 
	(for $\C^{n}$, \cite[Theorem 1.2]{MR2369936}) about the hole probability. They 
	proved that there exist constants $c_{1}\geq c_{2}>0$ such that for $r>0$ large,
	\begin{equation}
		\exp(-c_{1}r^{2n+2})\leq \P\big (\psi^{S}_{\eta}(z)\neq 0,\; \mathrm{for\;all\;}z\in 
		\mathbb{B}(0,r)\big)\leq \exp(-c_{2}r^{2n+2}),
		\label{eq:3.38}
	\end{equation}
	where $\mathbb{B}(0,r)=\{z\in \C^{n}\;:\; |z|<r\}$. 
	Let us now fix $r_{0}>0$. Then by \eqref{eq:1.6.14paris} and 
	\eqref{eq:1.2.4}, we get
	\begin{equation}
		\exp(-c \sqrt{p}^{2n+2})\leq \P\big( \psi^{S_{p}}_{\eta}(z)\neq 0,\; \mathrm{for\;all\;}z\in 
		\mathbb{B}(0,r_{0})\big)\leq \exp(-c'\sqrt{p}^{2n+2}).
		\label{eq:3.40}
	\end{equation}
	By using \eqref{eq:flatguassian}, the inequality 
	\eqref{eq:3.40} is equivalent to
		\begin{equation}
		\exp(-c \sqrt{p}^{2n+2})\leq \P\big( \psi^{S}_{\eta}(z)\neq 0,\; \mathrm{for\;all\;}z\in 
		\mathbb{B}(0,\sqrt{p}r_{0})\big)\leq \exp(-c'\sqrt{p}^{2n+2}).
		\label{eq:3.41}
	\end{equation}
	Therefore, we recover the estimates in \eqref{eq:3.38} by 
	approximating a sufficiently large $r>0$ by $\sqrt{p}r_{0}$.

\begin{Remark} 	
In the context of a general complete K\"{a}hler manifold $X$, an 
analogue question to \eqref{eq:3.92ss} would be as follows: for any relatively compact open subset 
$U\subset X$, find a sequence of orthonormal bases 
 $\{\widetilde{S}^{p}_{j}\}_{j=1}^{d_{p}}$ of $H^{0}_{(2)}(X,L^{p})$, 
 $p\in \N,$ such that
 \begin{equation}
 	\sup_{x\in\bar{U}} \sum_{j>K' 
 	p^{n}}|\widetilde{S}^{p}_{j}(x)|^{2}_{h^{p}}\leq C e^{-cp},
 \end{equation}
 where $C$, $K'$, $c$ are certain positive constants independent of 
 $p$, and the sum in the left-hand side is taken to be $0$ if $d_{p}=\dim 
 H^{0}_{(2)}(X,L^{p})\leq 
 K'p^{n}$. This question is trivial for the cases where $d_{p}=\mathcal{O}(p^{n})$ for $p\gg 0$.
 
 The existence of such a sequence of bases suggests that, on a relatively 
 compact subset, the Bergman projections or Bergman kernels can be 
 approximated by the orthogonal projections or their 
 kernels of a sequence of finite dimensional subspaces of 
 $H^{0}_{(2)}(X,L^{p})$. Moreover, one may expect a connection between the number (or dimension 
 of the aforementioned subspace) $K'p^{n}$ and the integration of 
 dimension density
 on $U$
 \begin{equation}
 	\int_{U}P_{p}(x,x)\mathrm{dV}(x).
 \end{equation}
\end{Remark}

\section{Random $\mathcal{L}^{2}$-holomorphic sections and Toeplitz 
operators}
In this section, we always assume the same conditions on $(X,\Theta)$ 
and $(L,h_{L})$ as in Section \ref{section2}: $(X,\Theta)$ is a 
complex Hermitian manifold (without boundary), and $(L,h_{L})$ is a 
Hermitian line bundle on $X$. We do not, however, assume any completeness for $\Theta$ or positivity 
for $(L,h_{L})$. 

The goal of this section is to introduce a method of \lq canonically randomizing\rq\  the 
$\mathcal{L}^{2}$-holomorphic sections of $L$ on $X$, in particular when 
$d=\dim H^{0}_{(2)}(X, L)=\infty$. 

As mentioned in the Introduction, this is achieved by the abstract Wiener space construction from 
probability theory with an approach via 
Toeplitz operators from the theory of geometric quantization. This  
induces a  Gaussian probability measure on the space 
of $\mathcal{L}^{2}$-holomorphic sections.

\subsection{Abstract Wiener spaces}
To define a Gaussian probability measure on an 
infinite dimensional Hilbert space, we here employ the construction 
of the abstract Wiener space introduced by Gross \cite{Gross67}. We 
also refer to the article of Sheffield \cite{Sheffield} for further motivation and 
 developments on this topic.

For a (complex) vector space $\mathcal{H}$, a Hermitian norm is a norm on 
$\mathcal{H}$ which is induced by a Hermitian inner product on it.

\begin{Definition}\label{def:4.1a}
	Let $(\mathcal{H},\|\cdot\|_{\mathcal{H}})$ be a separable Hilbert space 
	of infinite dimension. A Hermitian norm $\|\cdot\|$ is called 
	measurable if for all $\epsilon>0$, there exists a finite 
	dimensional subspace $F_{\epsilon}\subset \mathcal{H}$ such that 
	for $F\subset \mathcal{H}$ a subspace of finite dimension with $F\perp 
	F_{\epsilon}$, one has
	\begin{equation}
		\mu_{F,\; \|\cdot\|_{\mathcal{H}}} \big(\{x\in F\;:\; 
		\|x\|\geq \epsilon \}\big) < \epsilon,
		\label{eq:2.1.12}
	\end{equation}
	where $\mu_{F,\; \|\cdot\|_{\mathcal{H}}}$ denotes the standard 
	Gaussian measure on $F$ with respect to the Hermitian metric 
	associated with $\|\cdot\|_{\mathcal{H}}$.
\end{Definition}

\begin{Proposition}[cf.\ {\cite{Gross67},\cite[Chapter I: Theorem 4.3]{Kuo75}}]\label{prop:4.2m}
	Let $\mathcal{H}$ be a separable Hilbert space with norm 
	$\|\cdot\|_{\mathcal{H}}$ , and $\|\cdot\|$ be a continuous (with 
	respect to $\|\cdot\|_{\mathcal{H}}$) Hermitian norm on $\mathcal{H}$. Then the following two
conditions are equivalent:
\begin{enumerate}
	\item $\|\cdot\|$ is measurable.
	\item There exists a one-to-one Hilbert-Schmidt operator $T$ of $\mathcal{H}$ such
that $\|x\|=\|Tx\|_{\mathcal{H}}$ for $x\in \mathcal{H}$.
\end{enumerate}
\end{Proposition}

Given a measurable Hermitian norm $\|\cdot\|$ on $\mathcal{H}$, let 
$\mathcal{B}$ be the completion of $\mathcal{H}$ with respect to 
$\|\cdot\|$. Then $(\mathcal{B},\|\cdot\|)$ is a separable Hilbert 
space containing $\mathcal{H}$ as a dense subspace.

Let $\mathcal{B}^{*}$ be the topological dual space of $\mathcal{B}$. 
If $\alpha\in \mathcal{B}^{*}$, then $\alpha|_{\mathcal{H}}$ is a
continuous linear functional on 
$(\mathcal{H},\|\cdot\|_{\mathcal{H}})$. If $\alpha$ vanishes 
identically on $\mathcal{H}$, then it vanishes on $\mathcal{B}$. This 
way, we can regard $\mathcal{B}^{*}$ as a (dense) subspace of 
$\mathcal{H}^{*}$, where $\mathcal{H}^{*}$ can be identified with $\mathcal{H}$ via the 
Hilbert metric associated with $\|\cdot\|_{\mathcal{H}}$.

In a slight abuse of notation we denote by $\mathcal{S}$ the Borel $\sigma$-algebra of $\mathcal{B}$. Then each 
$\alpha\in \mathcal{B}^{*}$ is a Borel-measurable function from 
$\mathcal{B}$ to $\C$.
For $F\subset \mathcal{B}^{*}\subset \mathcal{H}$ an arbitrary finite 
dimensional subspace we introduce the notation
\begin{align}
\begin{split}
	\phi_{F}: \mathcal{B}&\rightarrow F,\; \\
	\phi_{F}(b)&=\sum_{j=1}^{\dim_{\C} F} (b,v_{j})v_{j},
\end{split}
	\label{eq:2.1.12s}
\end{align}
where $\{v_{j}\}$ is an orthonormal basis of $(F, \|\cdot\|_{\mathcal{H}})$.

Then Gross \cite{Gross67} proved the following result. 
\begin{Theorem}\label{thm:Gross}
	Fix a measurable norm $\|\cdot\|$ on $\mathcal{H}$ as above. 
	There exists a unique probability measure $\mathcal{P}$ on 
	$(\mathcal{B},\mathcal{S})$ such that for  $F\subset \mathcal{B}^{*}$ any finite dimensional 
	subspace,
	\begin{equation}
		\mathcal{P}(\phi^{-1}_{F}(U))=\mu_{F,\;\|\cdot\|_{\mathcal{H}}}(U),
		\label{eq:2.1.13}
	\end{equation}
	for all Borel subset  $U$ of $F$. The triple 
	$(\mathcal{B},\mathcal{S},\mathcal{P})$ is called an abstract 
	Wiener space.
\end{Theorem}

If $\alpha\in \mathcal{B}^{*}$, then as a function on $\mathcal{B}$, 
it is an element of 
$\cLL(\mathcal{B},\mathcal{S},\mathcal{P})$. We denote this map by 
\begin{equation}
	\Phi_{0}:\mathcal{B}^{*}\rightarrow 
	\cLL(\mathcal{B},\mathcal{S},\mathcal{P}).
	\label{eq:2.1.14}
\end{equation}
Moreover, for $\alpha\in \mathcal{B}^{*}$, $\Phi_{0}(\alpha)$ is a 
Gaussian random variable with zero mean  and variance 
$\|\alpha\|^{2}_{\mathcal{H}}$. The map $\Phi_{0}$ extends to a 
continuous linear map
\begin{equation}
	\Phi:\mathcal{H}^{*}\simeq \mathcal{H}\rightarrow 
	\cLL(\mathcal{B},\mathcal{S},\mathcal{P}),
	\label{eq:2.1.15}
\end{equation}
where for $y\in \mathcal{H}$, $\Phi(y)$ is a  Gaussian random variable 
 with zero mean and variance $\|y\|^{2}_{\mathcal{H}}$.

Remark that the above construction is trivial if $\mathcal{H}$ 
is finite dimensional; indeed, in this case the Hilbert space $\mathcal{B}$ is 
reduced to $\mathcal{H}$ itself. The probability measure constructed 
in Theorem \ref{thm:Gross} is the standard Gaussian 
probability measure on $\mathcal{H}$ with respect to the norm 
$\|\cdot\|_{\mathcal{H}}$.

\subsection{Toeplitz operators on $H^{0}_{(2)}(X,L)$} 
Recall that $P$ denotes the orthogonal projection from $\cLL(X,L)$ 
onto $H^{0}_{(2)}(X,L)$, and $P(x,y)$, $x,y\in X$, denotes the 
corresponding Bergman kernel. W.l.o.g.\ we may and do always assume that $d=\dim 
H^{0}_{(2)}(X,L)\geq 1$ in the following.
\begin{Definition}\label{def:4.4}
	For a bounded function $f\in \mathscr{C}^{\infty}(X,\C)$, set
	\begin{align}
T_{f}: H^{0}_{(2)}(X, L) \rightarrow H^{0}_{(2)}(X,L),\quad
		T_{f}:=P f P,
		\label{eq:2.1.1}
	\end{align}
	where the action of $f$ is the pointwise multiplication by $f$. 
	The operator $T_{f}$ is called the Toeplitz operator associated 
	with $f$. 
\end{Definition}
The integral kernel of $T_{f}$ is provided in the representation
\begin{equation}
	T_{f}(x,x')=\int_{X} 
P(x,x^{\prime\prime})f(x^{\prime\prime})P(x^{\prime\prime},x')
\mathrm{dV}(x^{\prime\prime}).
	\label{eq:2.1.2}
\end{equation}
Note also that the Hilbert adjoint of $T_{f}$ is 
$T_{\overline{f}}$.

We introduce a class of bounded smooth functions on $X$ whose 
associated Toeplitz operators are Hilbert-Schmidt.
\begin{Definition}\label{def:5.4}
	Let $\mathcal{Q}(X,L\,;\,\C)$ 
	be the vector space of bounded smooth complex
	functions $f$ on $X$ such that
	\begin{equation}
		\int_{X}|f(x)|P(x,x)\mathrm{dV}(x)<\infty,
		\label{eq:5.8paris}
	\end{equation}
	 where $P$ is the Bergman
	kernel of $L$.
\end{Definition}

\begin{Example}
(1) It is clear that $\mathscr{C}^{\infty}_{c}(X,\C)$ is a 
		subspace of $\mathcal{Q}(X,L\,;\,\C)$. In 
		particular, if $X$ is compact, then 
		\begin{equation}
			\mathcal{Q}(X,L\,;\,\C)=\mathscr{C}^{\infty}(X,\C).
			\label{eq:5.8bis}
		\end{equation}

\noindent
(2) Let $\mathscr{C}^{\infty}_{b}(\C^{n},\C)$ 
		denote the set of bounded smooth functions on $\C^{n}$. In the case of the Bargmann-Fock space (see
		Example \ref{ex:2.9}), we have
		\begin{equation}
			\mathcal{Q}(\C^{n},L\,;\,\C)=\mathscr{C}^{\infty}_{b}(\C^{n},\C)\cap \mathcal{L}^{1}(\C^{n}, \mathrm{dV}).
			\label{eq:5.9}
		\end{equation}

\noindent
(3) In general, with the assumptions as in Section 
		\ref{section3}, if we assume further that $(X,J,\Theta)$, $(L,h_{L})$ have bounded 
		geometry, by \cite[Theorem 6]{MM15}, there exist $c>0$, $C>0$ and 
		$p_{0}\in\N^{*}$ such that 
		for $p\geq p_{0}$,
		\begin{equation}
			cp^{n}\leq \inf_{x\in X} P_{p}(x,x)\leq \sup_{x\in X} P_{p}(x,x)\leq Cp^{n},
			\label{eq:5.10}
		\end{equation}
		that is, the Bergman kernel function $P_{p}(x,x)$ is bounded from above and away from zero
		on $X$.  
		As a consequence, we get that for $p\geq p_{0}$,
		\begin{equation}
			\mathcal{Q}(X,L^{p}\,;\,\C)=\mathscr{C}^{\infty}_{b}(X,\C)\cap \mathcal{L}^{1}(X, \mathrm{dV}).
			\label{eq:5.11}
		\end{equation}
\end{Example}

\begin{Proposition}\label{prop:4.7a}
	For $f\in 
	\mathcal{Q}(X,L\,;\,\C)$, the operator $T_{f}$ on $H^{0}_{(2)}(X,L)$ has smooth 
	Schwartz kernel and is
	Hilbert-Schmidt.
\end{Proposition}

\begin{proof}
	If $d=\dim H^{0}_{(2)}(X,L)<\infty$, then the statement is 
	trivial. Hence, we assume $d=\infty$ w.l.o.g.\ in the sequel. 
	Let $\{S_{j}\}_{j=1}^{\infty}$ be a complete Hilbert basis of 
	$H^{0}_{(2)}(X,L)$.
	
Note that for any compact set $K\subset X$, the series
\begin{equation}
\sum_{j=1}^{\infty}|S_{j}(x)|_{h_{L}}^{2}
\label{eq:5.15paris}
\end{equation}
converges uniformly for $x\in K$. As a consequence, for $K_{1}, 
K_{2} \subset X$ compact, the series
\begin{equation}
\sum_{j=1}^{\infty}S_{j}(x)\otimes (S_{j}(y))^{\ast}
\label{eq:5.16paris}
\end{equation}
converges absolutely and uniformly for $x\in K_{1}$ and $y\in 
K_{2}$ \cite[Proposition (2.4)]{At:76}. 
As follows from the properties of holomorphic functions, 
	if we replace $S_{j}(x)$, $(S_{j}(y))^{\ast}$ by their 
	respective covariant derivatives, then the series in 
	\eqref{eq:5.15paris} and \eqref{eq:5.16paris} are still 
	absolutely convergent on any given compact subsets.
	
	Note that by Definition \ref{def:5.4}, for $j\in \N^{\ast}$, the 
	function
	$X\ni x\mapsto f(x)|S_{j}(x)|^{2}_{h_{L}}$ is integrable on $X$ with 
	respect to $\mathrm{dV}$. Furthermore, for $x'\in X$, $i,j\in\N^{\ast}$, we have
	\begin{equation}
		f(x')(S_{i}(x'))^{\ast}S_{j}(x')=f(x')h_{L,x'}\left(S_{j}(x'),S_{i}(x')\right),
		\label{eq:5.17}
	\end{equation}
	and 
	\begin{equation}
		\int_{X}\left|f(x')\left(S_{i}(x')\right)^{\ast}S_{j}(x')\right|\mathrm{dV}(x')\leq \big\|\sqrt{|f|}S_{i}\big\|_{\mathcal{L}^{2}(X,L)}\cdot\big\|\sqrt{|f|}S_{j}\big\|_{\mathcal{L}^{2}(X,L)}
		\label{eq:5.1.18}
	\end{equation}
	
	Now we fix two compact subsets $K_{1}$, $K_{2}\subset X$. For 
	$x\in K_{1}$, $y\in K_{2}$ and $i,j\in \N^{\ast}$ we have
\begin{equation}
	\big|S_{i}(x)\otimes (S_{i}(x'))^{\ast}f(x')S_{j}(x')\otimes 
	(S_{j}(y))^{\ast}\big|\leq |S_{i}(x)|_{h_{L}}\cdot 
	\big|f(x')(S_{i}(x'))^{\ast}S_{j}(x')\big|\cdot 
	|S_{j}(y)|_{h_{L}},
	\label{eq:5.1.19}
\end{equation}
where the norm in the left-hand side is given by $h^{p}_{x}\otimes 
h^{p,\ast}_{y}$. By \eqref{eq:5.1.18} this entails
\begin{equation}
	\begin{split}
		&\int_{X}\big|S_{i}(x)\otimes (S_{i}(x'))^{\ast}f(x')S_{j}(x')\otimes 
	(S_{j}(y))^{\ast}\big|\mathrm{dV}(x')\\
	&\leq 
	|S_{i}(x)|_{h_{L}}\left\|\sqrt{|f|}S_{i}\left\|_{\mathcal{L}^{2}(X,L)}\cdot\right\|\sqrt{|f|}S_{j}\right\|_{\mathcal{L}^{2}(X,L)} 
	|S^{p}_{j}(y)|_{h^{p}}.
	\end{split}
	\label{eq:5.1.20}
\end{equation}
Putting things together, we arrive at
\begin{equation}
\begin{split}
	\sum_{i=1}^{\infty}|S_{i}(x)|_{h_{L}}\|\sqrt{|f|}S_{i}\|_{\mathcal{L}^{2}(X,L)} &\leq  \Big(\sum_{i=1}^{\infty}|S_{i}(x)|^{2}_{h_{L}}\Big)^{1/2}\Big(\sum_{i=1}^{\infty}\|\sqrt{|f|}S_{i}\|^{2}_{\mathcal{L}^{2}(X,L)} \Big)^{1/2}\\
	&=\Big(\sum_{i=1}^{\infty}|S_{i}(x)|^{2}_{h_{L}}\Big)^{1/2}\Big(\int_{X}|f(x')|P(x',x')\mathrm{dV}(x')\Big)^{1/2}\\
	&<\infty,
\end{split}
	\label{eq:5.1.21}
\end{equation}
and the above estimates still hold if we replace $S_{i}(x)$ by its 
covariant derivatives at $x$.

Recalling the Schwartz kernel of $T_{f}$ from 
\eqref{eq:2.1.2}, the above calculations show that 
$T_{f}(x,y)$ is a smooth section on $X\times X$. For proving that 
$T_{f}$ is Hilbert-Schmidt, it only remains to show that 
\begin{equation}
	\sum_{i,j}\left|\langle 
	T_{f}S_{i},S_{j}\rangle_{\cLL(X,L)}\right|^{2}<\infty.
	\label{eq:5.1.22}
\end{equation}
Indeed, by \eqref{eq:5.1.18}, we have
\begin{equation}
	\left|\langle 
	T_{f}S_{i},S_{j}\rangle_{\cLL(X,L)}\right|^{2}\leq 
	\left\|\sqrt{|f|}S_{i}\right\|_{\mathcal{L}^{2}(X,L)}^{2}\cdot\left\|\sqrt{|f|}S_{j}\right\|_{\mathcal{L}^{2}(X,L)}^{2}
\end{equation}
Then
\begin{equation}
	\begin{split}
			&\sum_{i,j}\left|\langle 
	T_{f}S_{i},S_{j}\rangle_{\cLL(X,L)}\right|^{2}
	\leq\sum_{i,j=1}^{\infty} \left\|\sqrt{|f|}S_{i}\right\|_{\mathcal{L}^{2}(X,L)}^{2}\cdot\left\|\sqrt{|f|}S_{j}\right\|_{\mathcal{L}^{2}(X,L)}^{2}\\
	&= \Big(\int_{X} |f(x)|P(x,x)\mathrm{dV}(x)\Big)^{2}<\infty.
	\end{split}
	\label{eq:5.1.23}
\end{equation}
This completes our proof.
\end{proof}

\begin{Corollary}
	If $f\in \mathscr{C}^{\infty}(X,\C)$ is with compact 
	support, then $T_{f}$ is a Hilbert-Schmidt operator on $H^{0}_{(2)}(X, L)$. 
	Moreover, $T_{f}$ is trace class, and
	\begin{equation}
		\mathrm{Tr}[T_{f}]=\int_{X}f(x)P(x,x)\mathrm{dV}(x).
		\label{eq:2.1.7}
	\end{equation}
\end{Corollary}

\subsection{Random $\mathcal{L}^{2}$-holomorphic 
sections}\label{ss:4.3b}

Let $\mathcal{Q}(X,L\,;\,\R_{\geq 0})$ be the subspace of 
$\mathcal{Q}(X,L\,;\,\C)$ consisting of the functions valued in 
$\R_{\geq 0}$. For $f\in \mathcal{Q}(X,L\,;\,\R_{\geq 0})$, $T_{f}$ 
is a nonnegative self-adjoint 
Hilbert-Schmidt (hence compact)
operator on $H^{0}_{(2)}(X,L)$. 
\begin{Lemma}
	For $0\neq f\in \mathcal{Q}(X,L\,;\,\R_{\geq 0})$,
	the operator $T_{f}: H^{0}_{(2)}(X,L) \rightarrow 
	H^{0}_{(2)}(X,L)$ is injective.
\end{Lemma}
\begin{proof}
	Since $f\neq 0$, there exists an open subset $U$ of $X$ on which 
	$f$ is strictly positive. If $s\in H^{0}_{(2)}(X,L)$ is such that 
	$T_{f}s=0$, then
	\begin{equation}
		\begin{split}
			0&=\langle T_{f}s, s\rangle
			=\int_{X}f(x)|s(x)|^{2}_{h_L}\mathrm{dV}(x),
		\end{split}
		\label{eq:2.1.8}
	\end{equation}
	and hence $s|_{U}=0$. Since $U$ is open and $s$ is holomorphic 
	on $X$, we get $s=0$. This proves the lemma.
\end{proof}

Fix $f\in \mathcal{Q}(X,L\,;\,\R_{\geq 0})$, $f\neq 0$.
If $d<\infty$, then the above $T_{f}$ is actually an isomorphism on the 
vector space
$H^{0}_{(2)}(X,L)$. Now we focus on the case of $d=\infty$. Since $T_{f}$ is compact and injective,  it cannot be 
surjective. Hence, it does not admit a bounded inverse. Moreover, for 
any $\lambda\in \C$, $\lambda\neq 0$, the operator $T_{f}-\lambda$ 
is Fredholm with closed range and Fredholm index $0$.

Set 
$D(T^{-1}_{f})=\mathrm{Range}(T_{f}:H^{0}_{(2)}(X,L) \rightarrow 
	H^{0}_{(2)}(X,L))\subset H^{0}_{(2)}(X,L)$, which is a dense subspace. The inverse of $T_{f}$ is 
	defined as 
	\begin{equation}
		T_{f}^{-1}:D(T^{-1}_{f})\subset 
		H^{0}_{(2)}(X,L)\rightarrow 
		H^{0}_{(2)}(X,L).
	\end{equation}

Let $\sigma(T_{f})\subset \R_{\geq 0}$ denote the spectrum of 
$T_{f}$, which is a countable set consisting of two parts: the point spectrum 
$\sigma_{\mathrm{p}}(T_{f})\subset \R_{> 0}$ (eigenvalues) and the residual spectrum 
$\sigma_{\mathrm{res}}(T_{f})=\{0\}$.	In this case, the point 
spectrum of $T_{f}$ (always with finite multiplicities) is 
a decreasing sequence of strictly positive real numbers, 
\begin{equation}
	\lambda_{1}\geq \lambda_{2}\geq \ldots\geq 
	\lambda_{m}\geq \ldots\rightarrow 0.
	\label{eq:2.1.9}
\end{equation}

Since any separable (complex) Hilbert space is isometric to the 
Hilbert 
space $\ell^{2}(\C)$ by choosing an orthonormal basis, we can choose
 an orthonormal basis $\{S_{j}\}^{\infty}_{j=1}$  of 
$H^{0}_{(2)}(X,L)$ with respect to the $\cLL$-metric 
 such that
\begin{equation}
	T_{f}S_{j}=\lambda_{j}S_{j}.
	\label{eq:2.1.10}
\end{equation}
If $S\in H^{0}_{(2)}(X,L)$, we can write uniquely
\begin{equation}
	S=\sum_{j\geq 1}a_{j}S_{j},\quad a_{j}\in\C.
	\label{eq:2.1.11}
\end{equation}
Then $(a_{j})_{j}\in \ell^{2}(\C)$, yielding the 
identification between $H^{0}_{(2)}(X,L)$ and 
$\ell^{2}(\C)$.

Since $T_{f}$ is one-to-one and Hilbert-Schmidt, by Proposition 
\ref{prop:4.2m},
$\|\cdot\|_{f}:=\|T_{f}\cdot\|$ defines a Hermitian measurable norm on 
$H^{0}_{(2)}(X,L)$. We denote by $\mathcal{B}_{f}(X,L)$ the 
completion of $H^{0}_{(2)}(X,L)$ with respect to 
$\|\cdot\|_{f}$ and set
\begin{equation}
	\ell^{2}_{f}(\C)=\Big\{(a_{j}\in\C)_{j\geq 1}\;:\;\sum_{j\geq 
	1}\lambda_{i}^{2}|a_{j}|^{2}<\infty \Big\}.
	\label{eq:2.1.16}
\end{equation}
It is clearly a separable Hilbert space, and using the basis as in \eqref{eq:2.1.10}, we have
\begin{equation}
	\mathcal{B}_{f}(X,L)\simeq \ell^{2}_{f}(\C).
	\label{eq:2.1.17}
\end{equation}

\begin{Proposition}\label{prop:2.10}
	Assume $d=\infty\,$, $0\neq f\in \mathcal{Q}(X,L\,;\,\R_{\geq 0})$. 
	Then the operator $T_{f}$ extends uniquely to an isomorphism of Hilbert spaces 
	\begin{equation}
		\widehat{T}_{f}: \big(\mathcal{B}_{f}(X,L),\|\cdot\|_{f} \big)\rightarrow 
		\big (H^{0}_{(2)}(X,L),\|\cdot\|_{\cLL(X,L)} \big).
		\label{eq:5.30}
	\end{equation}
\end{Proposition}
Given $0\neq f\in \mathcal{Q}(X,L\,;\,\R_{\geq 0})$, if $d<\infty$, 
we set
\begin{equation}
	\big(\mathcal{B}_{f}(X,L),\|\cdot\|_{f} \big)= \big(H^{0}_{(2)}(X,L),\|\cdot\|_{f} \big),\quad \text{and} \quad \widehat{T}_{f}:=T_{f}.
\end{equation}
Then we unify our notation for both cases $d<\infty$ and 
$d=\infty$.

\begin{Definition}\label{def:4.11}
	Denote by $\mathcal{P}_{f}$ the probability measure from  Theorem \ref{thm:Gross} with the choice
$\mathcal{B} = \mathcal{B}_{f}(X,L).$ Let $\P_{f}$ be the
Gaussian probability measure on $H^{0}_{(2)}(X,L)$ given 
by the pushforward of $\mathcal{P}_{f}$ through the 
isomorphism in \eqref{eq:5.30}. This way, we randomize the sections 
in $H^{0}_{(2)}(X,L)$.
\end{Definition}

\begin{Lemma}\label{lm:4.13ss}
	Assume $d\geq 1\,$, $0\neq f\in \mathcal{Q}(X,L\,;\,\R_{\geq 0})$. 
	For any nonzero $S\in H^{0}_{(2)}(X,L)$, the random variable on 
	$(H^{0}_{(2)}(X,L), \P_{f})$ defined as $H^{0}_{(2)}(X,L)\ni 
	s\mapsto \langle s, S\rangle_{\cLL(X,L)}\in\C$ is a centered 
	complex Gaussian variable with 
	variance $\|T_{f}S\|^{2}_{\cLL(X,L)}$.
\end{Lemma}
\begin{proof}
	Note that $T_{f}S$ is nonzero in $H^{0}_{(2)}(X,L)$, the linear 
	form
	\begin{equation}
		H^{0}_{(2)}(X,L)\ni s'\mapsto \langle s', 
		T_{f}S\rangle_{\cLL(X,L)}\in \C
	\end{equation}
	extends to a bounded linear form on 
	$(\mathcal{B}_{f}(X,L),\|\cdot\|_{f})$, hence defines an element 
	in $\mathcal{B}_{f}(X,L)^{\ast}$, denoted by $\Psi_{S}$. Then by 
	 property \eqref{eq:2.1.13}, the random 
	variable $\Psi_{S}(s')$ with $s'$ having the law 
	$\mathcal{P}_{f}$, is a centered complex Gaussian variable with 
	variance $\|T_{f}S\|^{2}_{\cLL(X,L)}$. 
	
	Put differently, by 
	construction, for $s'\in \mathcal{B}_{f}(X,L)$,
	\begin{equation}
		\Psi_{S}(s')=\langle \widehat{T}_{f}s',S\rangle_{\cLL(X,L)}.
	\end{equation}
	Thus as a random variable, it is exactly the same as $\langle s, 
	S\rangle_{\cLL(X,L)}$ with $s$ having distribution $\P_{f}$. 
	This completes our proof.
\end{proof}

\subsection{Zeros of random $\cLL$-holomorphic sections: proof of 
Theorem \ref{thm:4.15a}}\label{ss4.4a}
We assume $d\geq 1$, and we fix $0\neq f\in 
\mathcal{Q}(X,L\,;\,\R_{\geq 0})$. Set the operator $T^{2}_{f}:= 
	T_{f}\circ T_{f}$ on $H^{0}_{(2)}(X,L)$, which is a positive self-adjoint operator of 
	trace class. Let $T^{2}_{f}(x,y)$ denote the Schwartz kernel of 
	$T^{2}_{f}$.
\begin{Lemma}\label{lm:4.13a}
	The function $X\ni x\mapsto \log T^{2}_{f}(x,x)$ is locally 
	integrable on $X$, so that the $(1,1)$-current 
	$\partial\bar{\partial}\log T^{2}_{f}(x,x)$ is well-defined on 
	$X$.
\end{Lemma}
\begin{proof}
	Let $\{S_{j}\}_{j=1}^{d}$ be the orthonormal basis of 
	$H^{0}_{(2)}(X,L)$ as given in \eqref{eq:2.1.10}. Then for $x\in 
	X$,
	\begin{equation}
		T^{2}_{f}(x,x)=\sum_{j=1}^{d}\lambda_{j}^{2}|S_{j}(x)|^{2}_{h_{L}}.
		\label{eq:4.36ss}
	\end{equation}
	If $d=\infty$, the above sum is uniformly convergent on any compact subset of 
	$X$. Similar to the proof of Lemma \ref{lem:localintegral}, we 
	get that the function $\log T^{2}_{f}(x,x)$ is quasi-plurisubharmonic function 
	on $X$, hence locally integrable. This completes our proof.
\end{proof}

As an analog to the Fubini-Study current defined in \eqref{eq:intro2.23}, 
we set
\begin{equation}
	\gamma_{f}(L,h_{L})=c_{1}(L,h^{L})+\frac{\sqrt{-1}}{2\pi}\partial\bar{\partial}\log T^{2}_{f}(x,x).
	\label{eq:4.37ss}
\end{equation}

 We can now prove Theorem \ref{thm:4.15a} for the zeros of the random 
$\cLL$-holomorphic sections constructed in last subsection.
\begin{proof}[Proof of Theorem {\ref{thm:4.15a}}]
	Note that $T^{2}_{f}(x,x)$ vanishes exactly on $\text{Bl}(X,L)$. Let $\{S_{j}\}_{j=1}^{d}$ be the orthonormal basis of 
	$H^{0}_{(2)}(X,L)$ as given in \eqref{eq:2.1.10}.
	
	By Lemma \ref{lm:4.13ss}, the complex random variables
	\begin{equation}
		\eta_{j}:=\frac{1}{\lambda_{j}}\langle s, 
		S_{j}\rangle_{\cLL(X,L)},\;j=1,\,2,\,\ldots
	\end{equation}
	form an i.i.d.\ sequence of standard centered complex Gaussian 
	variable. As a consequence, we get that for $x\in X$,
	\begin{equation}
		s(x)=\sum_{j}\eta_{j}\lambda_{j}S_{j}(x).
		\label{eq:4.40ss}
	\end{equation}
	Then we can proceed as in the proof of Theorem 
	\ref{thm:expectation}, replacing $P(x,x)$ by $T^{2}_{f}(x,x)$ 
	given in \eqref{eq:4.36ss}, and we conclude \eqref{eq:3.1.1}.
\end{proof}

\begin{Remark}
 	In the above proof, we see that the random 
 	$\cLL$-holomorphic section $s$ with probability distribution $(H^{0}_{(2)}(X,L),\P_{f})$ is equivalent to the construction 
 	given in \eqref{eq:4.40ss}, as we explained the introduction part 
	(cf.\ \eqref{eq:intro1.2}). Let 
	$S=\{S_{j}\}_{j=1}^{d}$ denote the orthonormal basis of 
	$H^{0}_{(2)}(X,L)$ as given in \eqref{eq:2.1.10}, and let 
	$\psi^{S}_{\eta}$ be the Gaussian random holomorphic section 
	defined by \eqref{eq:2.2}, which can be regarded as a random 
	variable valued in $\mathcal{B}_{f}(X,L)$. Then the probability 
	space $(H^{0}_{(2)}(X,L),\P_{f})$ gives exactly the probability 
	distribution of the random 
	section $\widehat{T}_{f}\psi^{S}_{\eta}$.
\end{Remark}

\begin{Remark}\label{rk:4.16}
	Note that in the above constructions, we consider the nonnegative real function $f$ in order 
	to guarantee the injectivity of $T_{f}$ on $H^{0}_{(2)}(X,L)$. 
	One can also consider a different setting as follows where we do not require 
	the injectivity of $T_{f}$.
	
	Take $f$ in $\mathcal{Q}(X,L;\R)$, it can be negative somewhere 
	on $X$. Set 
	\begin{equation}
		H^{0}_{(2)}(X,L,f):=(\ker T_{f})^{\perp}=\overline{T_{f}H^{0}_{(2)}(X,L)}\subset 
		H^{0}_{(2)}(X,L),
	\end{equation}
	where $\overline{(\cdot)}$ means the closure in 
	$H^{0}_{(2)}(X,L)$.
	It is a Hilbert space, and the sections in $H^{0}_{(2)}(X,L,f)$ are the $\cLL$-holomorphic 
	sections of $L$ \textit{detected} by $f$. Note that $\ker T_{f}$ 
	is always finite dimensional. We consider the (self-adjoint)
	Hilbert-Schimdt operator
	\begin{equation}
		T_{f}^{\sharp}:=T_{f}|_{H^{0}_{(2)}(X,L,f)}: H^{0}_{(2)}(X,L,f)\rightarrow 
		H^{0}_{(2)}(X,L,f).
	\end{equation}
	Then we can proceed as in Subsection \ref{ss:4.3b} to construct a 
	respectively Gaussian probability measure $\P_{f}^{\sharp}$ on 
	$H^{0}_{(2)}(X,L,f)$. Let $s^{\sharp}$ denotes the corresponding random 
	section in $H^{0}_{(2)}(X,L,f)$, then
	\begin{equation}
		\mathbb{E}^{\P^{\sharp}_{f}}\left[[\Div(s^{\sharp})]\right]=\gamma_{f}(L,h_{L}),
		\label{eq:4.43b}
	\end{equation}
	where $\gamma_{f}(L, h_{L})$ is given by the same formula in 
	\eqref{eq:4.37ss}.
	
	One step further, since $\ker T_{f}$ is finite dimensional, we 
	can equip it with the standard Gaussian probability $\P_{f}^{0}$ measure 
	associated to the $\cLL$-metric. Take the product probability space
	\begin{equation}
		(H^{0}_{(2)}(X,L),\P_{f}):=(\ker T_{f}, \P_{f}^{0})\oplus
		(H^{0}_{(2)}(X,L,f),\P_{f}^{\sharp}).
		\label{eq:4.44b}
	\end{equation}
	Set $m(f):=\dim \ker T_{f}$, and let $\{S_{j}\}_{j=1}^{m(f)}$ be 
	an orthonormal basis of $\ker T_{f}$, then the Schwartz kernel of 
	the orthogonal projection $P_{\ker T_{f}}$ is given as
	\begin{equation}
		P_{\ker 
		T_{f}}(x,y)=\sum_{j=1}^{m(f)}S_{j}(x)\otimes(S_{j}(y))^{\ast}.
	\end{equation}
	
	Let $s$ be the random section in $H^{0}_{(2)}(X,L)$ with 
	probability
	distribution $\P_{f}$ constructed in \eqref{eq:4.44b}, then we 
	have
	\begin{equation}
		\mathbb{E}^{\P_{f}}\left[[\Div(s)]\right]=c_{1}(L,h^{L})+\frac{\sqrt{-1}}{2\pi}\partial\bar{\partial}\log \big(T^{2}_{f}(x,x)+P_{\ker T_{f}}(x,x)\big).
		\label{eq:4.46b}
	\end{equation}
	Note that since $f$ is bounded on $X$, then we always have
	\begin{equation}
		T^{2}_{f}(x,x)\leq T^{2}_{f}(x,x)+P_{\ker 
		T_{f}}(x,x)\leq\max\{\|f\|^{2}_{\infty},1\}P(x,x),
		\label{eq:4.47a}
	\end{equation}
	where $\|f\|_{\infty}$ is the $\mathcal{L}^{\infty}$-norm of $f$ 
	on $X$.
	
	We will consider the above different settings in Subsection 
	\ref{ss5.5} to study the random zeros for high tensor powers of 
	a prequautum line bundle on a complete K\"{a}hler manifold.
\end{Remark}

\section{Random {$\mathcal{L}^{2}$}-holomorphic sections for high 
tensor powers}
As an analog to Section \ref{section3}, we would like to study the asymptotic 
behaviors 
of the zeros of the random {$\mathcal{L}^{2}$}-holomorphic sections 
for high 
tensor powers of a given positive line bundle on $X$. We make the 
same assumptions for $(X,\Theta)$ and $(L,h_{L})$ as in the beginning 
of Section \ref{section3} (or in Subsection \ref{s1.2a}), in 
particular, we assume \eqref{eq:intro3.0.1}.

To construct in a canonical way the sequence of random 
$\mathcal{L}^{2}$-holomorphic sections of $L^{p}$, $p\in \N_{>0}$, we 
use the Toeplitz operators $\{T_{f,p}\}_{p\in \N_{>0}}$ associated with a suitable positive 
function $f$ on $X$. Such operators $\{T_{f,p}\}_{p\in \N_{>0}}$ are 
already well-studied in the context of Berzein-Toeplitz quantization.

\subsection{Asymptotics of Toeplitz operators}\label{ss5.1a}
Recall that $P_{p}$ denotes the orthogonal projection from $\cLL(X,L^{p})$ 
onto $H^{0}_{(2)}(X,L^{p})$. For a smooth bounded function $f$ on 
$X$ and $p\in \N_{>0}$, we set
\begin{equation}
	T_{f,p}=P_{p}f P_{p}.
\end{equation}
This defines a bounded linear operator acting on $H^{0}_{(2)}(X,L^{p})$.

To obtain the asymptotic expansion of the Schwartz kernels of 
$\{T_{f,p}\}$, we need further assumptions either on the function $f$ or on the geometry of $X$ and 
$L$. We are mainly concerned with the following two cases.
\begin{enumerate}[label=(\Roman*)]
   \item\label{caseI} We keep our assumptions on $(X,\Theta)$ and 
   $(L,h_{L})$ as in \eqref{eq:intro3.0.1}. Then the function $f$ is assumed to be a smooth 
   bounded function on $X$ which is constant outside a compact 
   subset of $X$.
	\item\label{caseII} In addition to the assumptions in Case \ref{caseI}, 
	we furthermore assume  that $(X,J,\Theta)$ 
	and $(L,h_{L})$ have bounded geometry (cf.\ Subsection \ref{ss3.1a}), and we take $f$ to be a 
	bounded smooth function on $X$ with bounded derivatives (with 
	respect to $\nabla^{TX}$ and $g^{TX}$) of any order.
\end{enumerate}
It is clear that in both cases, we can always take a smooth function 
$f$ with compact support. 

\begin{Theorem}[cf.\ {\cite[Chapter 7]{MM07},\cite{MM15}, \cite[Lemmas 
3.11, 3.14 \& 4.6]{Finski22a}}]
	Assume that $\{T_{f,p}\}_{p\in\N}$ is defined either in Case 
	\ref{caseI} or in Case \ref{caseII}, then we have the following results:
\begin{itemize}
	\item For a compact subset $K\subset X$ and for every 
	$\epsilon>0$, and every $\ell,m\in \N$, there exists 
	$C_{\ell,m,\varepsilon}>0$ such that for $p\geq 1$, $x,x'\in X$ with 
	$d(x,x')>\varepsilon$, we have
	\begin{equation}
		|T_{f,p}(x,x')|_{\mathscr{C}^{m}(K\times K)}\leq C_{\ell,m,\epsilon} 
		p^{-\ell},
		\label{eq:2.1.3}
	\end{equation}
	where the $\mathscr{C}^{m}$-norm is induced by $\nabla^{TX}$, 
	and $h^{L}$, $g^{TX}$.
	\item We have the asymptotic expansion as $p\rightarrow \infty$, 
	which is uniform on any compact subset of $X$,
	\begin{equation}
		T_{f,p}(x,x)=\sum_{\ell=0}^{\infty} b_{\ell,f}(x)p^{n-l} 
		+\mathcal{O}(p^{-\infty}),
		\label{eq:2.1.4}
	\end{equation}
	where $b_{\ell,f}\in \mathcal{C}^{\infty}(X,\C)$, in particular, 
	\begin{equation}
		b_{0,f}(x)=f(x),\, x\in X.
	\end{equation}
	\item The operator norms of $T_{f,p}$, $p\in\N$, satisfy
	\begin{equation}
		\lim_{p\rightarrow \infty} 
		||T_{f,p}||=||f||_{\infty}.
		\label{eq:2.1.6}
	\end{equation}
	\item If $g$ is also a bounded smooth function on $X$ in the 
	same case as $f$ (Case 
	\ref{caseI} or \ref{caseII} from above), then on
	any given compact subset $K\subset X$, we have the uniform 
	expansion
	\begin{equation}
		(T_{f,p}T_{g,p})(x,x)=p^{n}f(x)g(x) + \mathcal{O}(p^{n-1}),
		\label{eq:4.1.22ss}
	\end{equation}
	the expansion still holds if we take the derivatives with respect 
	to $x$ of any given order on 
	both sides. 
\end{itemize}
In particular, for the Case \ref{caseII}, we can refine 
\eqref{eq:2.1.3} to an exponential decay with respect to $\sqrt{p}$, and the 
results \eqref{eq:2.1.3}, \eqref{eq:2.1.4} and \eqref{eq:4.1.22ss} 
hold uniformly on the whole manifold $X$.
\end{Theorem}

The above theorem for Case \ref{caseI} was mainly proved by Ma and Marinescu in 
\cite[Chapter 7]{MM07}. For Case \ref{caseII}, it can be proved by a 
variation of the arguments in \cite[Chapter 7]{MM07} by using the 
exponential estimate for the Bergman kernel obtained in \cite{MM15},
these proofs are explained by Finski in \cite[Sections 3 \& 
4]{Finski22a}.

Our results in the sequel will mainly employ the expansion 
\eqref{eq:4.1.22ss} with $g=f$. Note that with further geometric conditions on 
$(X,\Theta)$ and $(L,h_{L})$, we have a refined version of 
\eqref{eq:4.1.22ss}.

Let $\mathrm{Ric}$ denote the 
	Ricci curvature tensor, and set
	$\mathrm{Ric}_{\Theta}=:\mathrm{Ric}(J\cdot,\cdot)$. Let $\mathbf{r}^{X}$ denote the scalar curvature of $(X, 
	g^{TX})$, and let $\Delta$ be the (positive) Bochner Laplacian associated with 
	$g^{TX}$ acting on the functions. We will use 
	$\langle\cdot,\cdot\rangle$ to denote the $\C$-linear extension 
	of the inner product $g^{\Lambda^{\bullet}T^{*}X}$. Consider the 
	connection $\nabla^{T^{*}X}: \mathscr{C}^{\infty}(X, T^{*}X\otimes\C)\rightarrow 
	\mathscr{C}^{\infty}(X, T^{\ast}X\otimes T^{\ast}X\otimes\C)$, let 
	$D^{0,1}$, $D^{1,0}$ denote the its respective $(1,0)$, $(0,1)$ 
	components.

	The following theorem was proved in \cite{MM12} for a compact 
	K\"{a}hler manifold equipped with a prequantum line bundle, where Ma 
	and Marinescu remarked in the introduction part that the 
	computations are essentially local and then extend to the case of complete 
	(noncompact) 
	K\"{a}hler manifolds. In particular, as a consequence of
	\cite[Sections 7.4 \& 7.5]{MM07} (for the Case \ref{caseI}) and 
	\cite{MM15} \cite[Sections 3 \& 
4]{Finski22a} (for Case \ref{caseII}), these results hold for both 
our
cases \ref{caseI}, \ref{caseII}.

	\begin{Theorem}\label{thm:5.2ss}
	Assume that $(X,\Theta)$ is complete K\"{a}hler and that 
	$(L,h_{L})$ is the prequantum line bundle on $X$ (i.e. 
	$\Theta=c_{1}(L,h_{L})$). Let 
	$f,g$ be bounded smooth functions where are constants outside a 
	compact subset (Case \ref{caseI}), or if in addition $(X,\Theta), (L,h_{L})$ have the bounded 
	geometry, let $f,g$ be two bounded smooth functions on $X$ such 
	that their derivatives of any order are also bounded on $X$ (Case 
	\ref{caseII}). Then 
	for $\ell\in\N$, there exists a smooth function on $X$, denoted by 
	$\mathbf{b}_{\ell}(f,g)$, which is a polynomial in the derivatives 
	of $f,g$ with coefficients depending only on $\Theta$ and 
	$h_{L}$, such that on any compact subset $K\subset X$, we have 
	the uniform expansion as follows ($N\geq 0$),
	\begin{equation}
		(T_{f,p}T_{g,p})(x,x)=\sum_{\ell=0}^{N} 
		p^{n-\ell}\mathbf{b}_{\ell}(f,g)(x)+\mathcal{O}(p^{n-N-1}).
		\label{eq:5.7ss}
	\end{equation}
	Furthermore, we have
	\begin{equation}
		\begin{split}
			\mathbf{b}_{0}(f,g)=&fg,\\
			\mathbf{b}_{1}(f,g)=&\frac{\mathbf{r}^{X}}{8\pi}fg-\frac{1}{4\pi}\big((\Delta f)g+f(\Delta g)\big)+\frac{1}{2\pi}\langle \overline{\partial}f,\partial g\rangle,\\
			\mathbf{b}_{2}(f,g)=&\frac{1}{32\pi^{2}}\Big(f(\Delta^{2}g)+(\Delta^{2}f)g- \mathbf{r}^{X}\big(f(\Delta g)+(\Delta f)g\big) \Big)\\
			&-\frac{\sqrt{-1}}{8\pi^{2}}\langle 
			\mathrm{Ric}_{\Theta}, 
			f\partial\bar{\partial}g+g\partial\bar{\partial}f\rangle\\
			&+\frac{1}{8\pi^{2}}\Big\{\frac{1}{2}\Delta f\cdot\Delta 
			g+\frac{\mathbf{r}^{X}}{2}\langle 
			\overline{\partial}f,\partial g\rangle +\langle 
			D^{0,1}\overline{\partial}f,D^{1,0}\partial 
			g\rangle_{g^{T^{\ast}X\otimes T^{\ast}X}}\\
			&-\langle 
			\overline{\partial}\Delta f,\partial g\rangle-\langle 
			\overline{\partial}f,\partial \Delta 
			g\rangle\Big\}.
		\end{split}
		\label{eq:5.8a}
	\end{equation}
\end{Theorem}

\subsection{Random zeros on the 
support: proofs of Theorems \ref{thm:1.6a} \& \ref{thm:intro5.8a}}\label{ss:5.2b}
Fix a $p_{0}\in \N_{>0}$, set
\begin{equation}
	\mathcal{Q}_{\geq p_{0}}(X,L;\R_{\geq 0}):=\cap_{p\geq p_{0}} 
	\mathcal{Q}(X,L^{p};\R_{\geq 0}).
\end{equation}

We fix a function $f$ as follows:

\textbf{\hypertarget{AssumpA}{Assumption A}}: $f\in \mathcal{Q}_{\geq p_{0}}(X,L;\R_{\geq 0})$,  
which is nontrivial and also satisfies the condition in Case 
	\ref{caseI} or in Case \ref{caseII}.
	
Note that such function 
	always exists, for instance, the nonnegative smooth functions on $X$ 
	with compact support, and in the case of Bargmann-Fock space, we 
	can take $f$ to be a nonnegative Schwartz function on $\C^{n}$. 
	In the rest of this section, we always consider the integer $p\geq 
p_{0}$. 

Following the construction in Definition \ref{def:4.11}, let 
$\P_{f,p}$ be the corresponding probability measure on 
$H^{0}_{(2)}(X,L^{p})$. Then we will denote by $\mathbf{S}_{f,p}$ the 
random section in $H^{0}_{(2)}(X,L^{p})$ given by the probability 
distribution $(H^{0}_{(2)}(X,L^{p}),\P_{f,p})$.

By \eqref{eq:4.1.22ss}, on any compact subset $K\subset X$ and for 
$\ell\in\N$, we 
have the following identity hold uniformly in $\mathscr{C}^{\ell}$-norm for $x\in K$
\begin{equation}
	T^{2}_{f,p}(x,x)=f^{2}(x)p^{n}+\mathcal{O}(p^{n-1}),
	\label{eq:3.1.2}
\end{equation}
If we are in Case \ref{caseII}, it holds uniformly over the whole manifold $X$.

Let $U$ be an open subset of $X$, and let 
$\Omega^{(n-1,n-1)}_{0}(\overline{U})$ denote the smooth 
$(n-1,n-1)$-forms on $\overline{U}$ with compact support in $U$. For any 
$(1,1)$-current $\alpha$ on $X$, let $\alpha|_{U}$ denote its restriction on 
$U$ by acting on sections in 
$\Omega^{(n-1,n-1)}_{0}(\overline{U})$.
\begin{Theorem}\label{thm:5.4a}
	Let $U$ be an open subset of $X$ such that 
$f>0$ on $U$, then we have the weak convergence of currents on $U$ as $p\rightarrow 
	\infty$,
	\begin{equation}
		\frac{1}{p}\E^{\P_{f,p}}[[\Div(\mathbf{S}_{f,p})]|_{U}]\rightarrow c_{1}(L,h_{L})|_{U}.
		\label{eq:5.14a}
	\end{equation}
\end{Theorem}
\begin{proof}
	By \eqref{eq:4.37ss}, we get
	\begin{equation}
		\gamma_{f}(L^{p},h_{p})=pc_{1}(L,h_{L})+\frac{\sqrt{-1}}{2\pi}\partial\overline{\partial}\log{T_{f,p}^{2}(x,x)}.
	\end{equation}
	Note that by our assumption of $f$ on $U$, on any compact subset 
	of $U$, for sufficiently large 
	$p$, $\gamma_{f}(L^{p},h_{p})$ is a smooth form.
	
	Then by Theorem \ref{thm:4.15a}, we get
	\begin{equation}
		\frac{1}{p}\E^{\P_{f,p}}[[\Div(\mathbf{S}_{f,p})]|_{U}]=c_{1}(L,h_{L})|_{U}+\frac{\sqrt{-1}}{2\pi p}\partial\overline{\partial}\log{T_{f,p}^{2}(x,x)}.
		\label{eq:5.16a}
	\end{equation}

	For any $\varphi\in \Omega^{(n-1,n-1)}_{0}(\overline{U})$. Set 
	$K=\supp \varphi$, which is a compact subset of $U$.
	Set $m_{K}:=\max_{x\in K}f(x)$, $c_{K}:=\min_{x\in K}f(x)>0$, then for sufficiently 
	large $p$, $x\in K$, we have
	\begin{equation}
	2m_{K}p^{n}\geq T_{f,p}^{2}(x,x)\geq \frac{1}{2}c_{K}p^{n}.
\label{eq:5.17a}
	\end{equation}
	We can compute directly
	\begin{equation}
		\partial\overline{\partial}\log{T_{f,p}^{2}(x,x)}=\frac{1}{T_{f,p}^{2}(x,x)^{2}}\left(\partial\overline{\partial} T_{f,p}^{2}(x,x)-\partial T_{f,p}^{2}(x,x)\wedge \overline{\partial}T_{f,p}^{2}(x,x)\right).
		\label{eq:5.18a}
	\end{equation}
	Then by the uniform expansion \eqref{eq:3.1.2}, we get, as 
	$p\rightarrow\infty$,
	\begin{equation}
 		\Big \langle \frac{\sqrt{-1}}{2\pi 
 		p}\partial\overline{\partial}\log{T_{f,p}^{2}(x,x)}, 
 		\varphi \Big \rangle\rightarrow 0.
 	\end{equation}
 Then convergence in \eqref{eq:5.14a} follows.

\end{proof}
The following corollary is clear.
\begin{Corollary}
	If $f>0$ on $X$, then we have the weak convergence of currents on 
	$X$ as $p\rightarrow \infty$,
	\begin{equation}
		\frac{1}{p}\E^{\P_{f,p}}[[\Div(\mathbf{S}_{f,p})]]\rightarrow c_{1}(L,h_{L}).
		\label{eq:5.19a}
	\end{equation}
\end{Corollary}

By considering the sequence of random sections in the product 
probability space,
\begin{equation}
	(\mathbf{S}_{f,p})_{p}\in 
\Pi_{p}\big(H^{0}_{(2)}(X,L^{p}),\P_{f,p}\big),
\end{equation}
we also have the following convergence in probability one.
\begin{Theorem}\label{thm:5.5b}
	Let $U$ be an open subset of $X$ such that 
$f>0$ on $U$, then for any $\varphi\in 
\Omega^{(n-1,n-1)}_{0}(\overline{U})$, we have
\begin{equation}
	\P\left(\lim_{p\rightarrow 
	\infty}\frac{1}{p}\big\langle[\Div(\mathbf{S}_{f,p})],\varphi \big\rangle=\langle c_{1}(L,h_{L}),\varphi\rangle \right)=1.
	\label{eq:5.20a}
\end{equation}
\end{Theorem}
\begin{proof}
	Fix a nonzero $\varphi\in 
\Omega^{(n-1,n-1)}_{0}(\overline{U})$.
	Note that from the proof of Theorem \ref{thm:5.4a}, we have the convergence 
	\begin{equation}
		\lim_{p\rightarrow\infty}\Big\langle 
		\frac{1}{p}\gamma_{f}(L^{p},h_{p}),\varphi\Big\rangle=\langle 
		c_{1}(L,h_{L}),\varphi\rangle.
	\end{equation}
	
Defining the random variable
\begin{equation}
	Y_{f,p}=\frac{1}{p}\Big\langle 
	[\Div(\mathbf{S}_{f,p})]-\gamma_{f}(L^{p},h_{p}),\varphi\Big\rangle,
\end{equation}
the statement \eqref{eq:5.20a} is equivalent to proving that almost 
surely one has
\begin{equation}
	Y_{f,p}\rightarrow 0.
	\label{eq:5.23a}
\end{equation}

Note that if we use the construction from the proof of Theorem 
\ref{thm:4.15a}, we can write
\begin{equation}
	\mathbf{S}_{f,p}=\sum_{j=1}^{d_{p}}\eta^{p}_{j}\lambda^{p}_{j}S^{p}_{j},
\end{equation}
where $\{\eta^{p}_{j}\}_{j}$ is a sequence of i.i.d.\ standard complex 
Gaussian random variables, $\{\lambda^{p}_{j}\}_{j}$ is the point 
spectrum of $T_{f,p}$, and $\{S^{p}_{j}\}_{j}$ is the orthonormal 
basis of $H^{0}_{(2)}(X,L^{p})$ given by the eigensections of 
$T_{f,p}$.

Then, as explained in Remark \ref{rk:3.9a}, we can proceed as in the 
proof of Theorem \ref{thm:2.15}, so that we get
\begin{equation}
	\E[|Y_{f,p}|^{2}]=\mathcal{O}\Big(\frac{1}{p^{2}}\Big),
\end{equation}
which entails \eqref{eq:5.23a}, and hence \eqref{eq:5.20a}.
\end{proof}

It is natural to investigate a relaxations of the assumptions from Theorem \ref{thm:5.5b} as follows. For $f$ as above, consider $U$ an open subset of $\supp f$. In general,
$f$ might vanish at some points in $U,$ and it is a natural and interesting question to understand for which kind of conditions on the vanishing points of $f$ in 
$U$ we still can have the equidistribution results for the 
random zeros on $U$ as above. Since $f$ is nonnegative, if 
$f(x_0)=0$, the least 
possible 
vanishing order of $f$ at $x_{0}$ is $2$. In the sequel we will 
explain, if $f$ has only vanishing points of order $2$ at which 
$\Delta f$ does not vanish, then the above results still hold (under 
prequantum setting).

For this purpose, we will employ the results in Theorem \ref{thm:5.2ss}, so that we need to 
make the following assumption, which is stronger than \hyperlink{AssumpA}{\bf 
Assumption A}.

\textbf{\hypertarget{AssumpB}{Assumption B}}: assume that $(X,\Theta)$ is complete K\"{a}hler and that $(L,h_{L})$ is the prequantum line bundle on $X$ (i.e. 
	$\Theta=c_{1}(L,h_{L})$). Let 
	$f$ be a bounded smooth function where are constants outside a 
	compact subset, or if in addition $(X,\Theta), (L,h_{L})$ have the bounded 
	geometry, let $f$ be a bounded smooth function on $X$ such 
	that their derivatives of any order are also bounded on $X$.
 
\begin{Proposition}\label{prop:5.7a}
Assume that \hyperlink{AssumpB}{\bf Assumption B} holds for some nontrivial $f\in \mathcal{Q}_{\geq p_{0}}(X,L;\R_{\geq 0})$,
and let $U$ be an open subset of $X.$ If $f$ only vanishes up to order $2$ 
in $U$ and $\Delta f$ is nonzero at all vanishing points of $f$, then for any compact subset $K$ of $U$, there exists a 
constant $c_{K}>0$ and $p_{K}\geq p_{0}$ such that for $x\in K$, 
$p\geq p_{K}$,
\begin{equation}
	T^{2}_{f,p}(x,x)\geq c_{K}p^{n-2}.
	\label{eq:5.27a}
\end{equation}
Moreover, $\log f^{2}$ is locally integrable on $U$, and we have 
weak convergence of currents on $U$ as $p\rightarrow \infty$,
\begin{equation}
	\partial\overline{\partial} \log T^{2}_{f,p}(x,x)\rightarrow 
	\partial\overline{\partial} \log f^{2}.
	\label{eq:5.28a}
\end{equation}
Around a point $x$ where $f(x)>0$, the convergence in 
\eqref{eq:5.28a} holds in any local $\mathscr{C}^{\ell}$-norms.
\end{Proposition}
\begin{proof}
	Under the \hyperlink{AssumpB}{\bf Assumption B}, we can apply Theorem \ref{thm:5.2ss} 
	to $T^{2}_{p,f}$. Let $x_{0}\in U$ be a vanishing point of $f$, 
	by our assumption on $f$, we have
	\begin{equation}
		-\Delta f(x_{0})\neq 0.
	\end{equation}
	
	By taking a suitable geodesic normal 
	coordinate system 
	$(Y=(y_{j})_{j=1}^{2n}\in\R^{2n})$ center at $x_{0}$, we can expance
	the function $f$ near $x_{0}$ as
	\begin{equation}
		f(Y)=\sum_{j}c_{j}(x_{0})y_{j}^{2}+\mathcal{O}(|Y|^{3}),
		\label{eq:5.29a}
	\end{equation}
	where the constants $c_{j}(x_{0})\ge 0$ since $f\geq 0$. Then  
	\begin{equation}
		-\Delta f(x_{0})=\sum_{j}c_{j}(x_{0})>0.
	\end{equation}
	
	Now we compute the terms $\mathbf{b}_{\ell}(f,f)$, $\ell=1,2$, from \eqref{eq:5.8a} near $x_{0}$, 
	\begin{equation}
		\begin{split}
			&\mathbf{b}_{1}(f,f)=\frac{1}{8\pi}(\mathbf{r}^{X}f-4\Delta f)f+\frac{1}{2\pi}|\partial f|^{2},\\
			&\mathbf{b}_{2}(f,f)=\frac{1}{4\pi^{2}}\Big(\sum_{j}c_{j}(x_{0}) \Big)^{2}+\frac{1}{8\pi^{2}}|D^{0,1}\overline{\partial}f(x_{0})|^{2}_{g^{T^{\ast}X\otimes T^{\ast}X}}+\mathcal{O}(|Y|).
		\end{split}
		\label{eq:5.31a}
	\end{equation}
	Setting
	\begin{equation}
		\begin{split}
			\mu(f,x_{0})=\frac{1}{4\pi^{2}} \Big(\sum_{j}c_{j}(x_{0})\Big)^{2}+\frac{1}{8\pi^{2}}|D^{0,1}\overline{\partial}f(x_{0})|^{2}_{g^{T^{\ast}X\otimes T^{\ast}X}}>0,
		\end{split}
		\label{eq:5.32a}
	\end{equation}
	 we can choose a small open neighborhood $V_{x_{0}}$ of 
	$x_{0}$ such that for $x\in V_{x_{0}}$,
	\begin{equation}
		\mathbf{r}^{X}_{x}f(x)-4\Delta f(x)\geq 0, \quad \text{and} \quad \mathbf{b}_{2}(f,f)(x)\geq 
		\frac{1}{2} \mu(f,x_{0}),
	\end{equation}
	and so
	\begin{equation}
		\mathbf{b}_{1}(f,f)(x)\geq 0.
	\end{equation}
	Since $\mathbf{b}_{0}(f,f)=f^{2}$, then from the above 
	computations and \eqref{eq:5.17a}, 
	we get \eqref{eq:5.27a}.

	By \eqref{eq:5.29a}, on a sufficiently small open neighborhood 
	of $x_{0}$, we have
	\begin{equation}
		f(Y)\geq \frac{1}{2}\sum_{j} c_{j}(x_{0})y_{j}^{2}.
	\end{equation}
	Then it is clear that $\log f^{2}$ is integrable 
	near $x_{0}$. Then the current $\partial\overline{\partial}\log 
	f^{2}$ is well defined on $U$. Near a point where $f$ does not 
	vanish, we get the strong convergence of \eqref{eq:5.28a} by 
	means of \eqref{eq:4.1.22ss} and \eqref{eq:5.18a}.
	
	Now we focus on the point $x_{0}$ with $f(x_{0})=0$. Note that
	\begin{equation}
		p^{-n}T^{2}_{f,p}(x,x)=f^{2}+b_{1}(f,f)p^{-1}+b_{2}(f,f)p^{-2}+\mathcal{O}(p^{-3}).
	\end{equation}
	By \eqref{eq:5.31a}, we can take a small open neighborhood 
	$V'_{x_{0}}$ of 
	$x_{0}$ such that for $x\in V'_{x_{0}}$, $p\gg 0$,
	\begin{equation}
	b_{1}(f,f)(x)p^{-1}+b_{2}(f,f)(x)p^{-2}+\mathcal{O}(p^{-3})\geq 
	0,\quad \text{and} \quad f^{2}(x)\leq p^{-n}T^{2}_{f,p}(x,x)\leq 1.
	\end{equation}
	Then on $V'_{x_{0}}$, we have
	\begin{equation}
		\big|\log(p^{-n}T^{2}_{f,p}(x,x))\big|\leq |\log f^{2}(x)|.
	\end{equation}
	At the same time we have the pointwise convergence of functions 
	as $p\rightarrow\infty$,
	\begin{equation}
		\log(p^{-n}T^{2}_{f,p}(x,x))\rightarrow \log f^{2}(x).
	\end{equation}
	Since $\log f^{2}$ is integrable near $x_{0}$, by the dominated 
	convergence theorem, we get the convergence of $(1,1)$-currents 
	in \eqref{eq:5.28a} on $V'_{x_{0}}$, hence on $U$. This completes the proof.
\end{proof}

\begin{Remark}
	In the proof of Proposition \ref{prop:5.7a}, we see that if $f$ 
	has at least one vanishing point in $K\subset U$, then the power $(n-2)$ in 
	\eqref{eq:5.27a} can not be improved; otherwise, a lower bound 
	of $T_{f,p}^{2}(x,x)$ on $K$ is given as in \eqref{eq:5.17a}. When $X$ is compact, this observation 
	indicates that if $f\geq 0$ has only proper 
	vanishing points of order $2$ and at least one of such vanishing 
	point, then the lowest eigenvalue of $T_{f,p}$ 
	should behave like $\mathcal{O}(\frac{1}{p})$ as $p$ grows. For 
	this kind of results, we refer to the papers \cite{Del2019, 
	Del2020} of Deleporte. In particular, when $X$ is compact, the 
	lower bound in \eqref{eq:5.27a} can be deduced from 
	\cite{Del2019}.
\end{Remark}

As a direct consequence of Proposition \ref{prop:5.7a}, we obtain:
\begin{Theorem}\label{thm:5.8a}
	We suppose that \hyperlink{AssumpB}{\bf Assumption B} holds with a nontrivial $f\in \mathcal{Q}_{\geq p_{0}}(X,L;\R_{\geq 0})$. Let $U$ be an open subset of 
	$\supp f$ be such that $f$ only vanishes up to order $2$ 
in $U$ with nonzero $\Delta f$ at the vanishing points. Then as $p\rightarrow \infty$,
\begin{itemize}
	\item we have the convergence of $(1,1)$-currents on $U$
	\begin{equation}
		\frac{1}{p}\E^{\P_{f,p}}[[\Div(\mathbf{S}_{f,p})]|_{U}]\rightarrow c_{1}(L,h_{L})|_{U}.
	\end{equation}
	\item for any $\varphi\in 
\Omega^{(n-1,n-1)}_{0}(\overline{U})$, we have
\begin{equation}
	\P\left(\lim_{p\rightarrow 
	\infty}\frac{1}{p}\langle[\Div(\mathbf{S}_{f,p})],\varphi\rangle=\langle c_{1}(L,h_{L}),\varphi\rangle\right)=1.
\end{equation}
\end{itemize}
\end{Theorem}

\subsection{Higher fluctuation of random zeros near points of vanishing 
order two}\label{ss5.3ab}
In this subsection, we always assume \hyperlink{AssumpB}{\bf 
Assumption B} to hold, and we investigate the random zeros of $\mathbf{S}_{f,p}$ near 
a proper vanishing point $f$ with vanishing order $2$, up to a 
scale $\sim \frac{1}{\sqrt{p}}$, so called Planck scale. Note that in \cite{AF:2022}, for a compact K\"{a}hler manifold $X$ 
and under a different assumption on $f$, Ancona and Le Floch observed and proved the 
phenomenon that the random zeros fluctuate a bit more near the zeros of 
$f$. We will observe the similar situation for our setting, for this purpose, we need to refine the computations in 
\eqref{eq:5.31a} in a complex coordinate system centered at $x_{0}$ 
where $f$ vanishes with order $2$.

Suppose $f\geq 0$ and that $x_{0}$ is a vanishing point of $f$ with $\Delta 
f(x_{0})<0$. Then we can choose a holomorphic coordinate system centered 
at $x_{0}$, denoted by $z=(z_{1},\ldots,z_{n})\in\C^{n}$, such that
\begin{equation}
	g^{TX}_{z}=g^{\C^{n}}_{\mathrm{st}}+\mathcal{O}(|z|^{2}),
	\label{eq:5.41a}
\end{equation}
where $g^{\C^{n}}_{\mathrm{st}}$ denotes the standard Euclidean 
metric on $\C^{n}\simeq \R^{2n}$.

Note that we view $z$ as a column vector, and let $(\cdot)^{T}$ denote the transpose of a matrix.
In this coordinate system, we can write
\begin{equation}
	f(z)=z^{T}A\bar{z}+z^{T}Bz+\bar{z}^{T}B\bar{z}+\mathcal{O}(|z|^{3}),
	\label{eq:5.42b}
\end{equation}
where the matrix $A$ is Hermitian and semipositive definite, $B$ is 
symmetric complex matrix, they are determined uniquely by the Hessian 
of $f$ at $x_{0}$.
Set
\begin{equation}
	\hat{f}_{x_{0}}(z)=z^{T}A\bar{z}+z^{T}Bz+\bar{z}^{T}\bar{B}\bar{z}.
\end{equation}
Since $f\geq 0$, then for any $z\in \C^{n}$ with $\|z\|=1$,
\begin{equation}
	z^{T}A\bar{z}\geq 2|\Re(z^{T}Bz)|,
\end{equation}
where $\Re(\cdot)$ denotes the real part.
In particular, $\hat{f}_{x_{0}}(z)\geq 0$.

Using this complex coordinate system, we compute
\begin{equation}
	\begin{split}
	\Delta f(z)&=-4\mathrm{Tr}[A]+\mathcal{O}(|z|),\\
	|\partial 
	f(z)|^{2}&=2|A\bar{z}+2Bz|^{2}+\mathcal{O}(|z|^{3}),\\
	|D^{0,1}\overline{\partial}f(z)|^{2}_{g^{T^{\ast}X\otimes 
	T^{\ast}X}}&=16\mathrm{Tr}[B\bar{B}^{T}]+\mathcal{O}(|z|).
	\end{split}
\end{equation}
Note that $\mu(f,x_{0})$ is defined in \eqref{eq:5.32a}, then we have
\begin{equation}
	\mu(f,x_{0})=\frac{1}{\pi^{2}}(\mathrm{Tr}[A])^{2}+\frac{2}{\pi^{2}}\mathrm{Tr}[B\bar{B}^{T}]>0.
\end{equation}

Then we rewrite the computations in 
\eqref{eq:5.31a} as follows,	
\begin{equation}
		\begin{split}
			&\mathbf{b}_{0}(f,f)(z)=\hat{f}_{x_{0}}^{2}(z)+\mathcal{O}(|z|^{5}),\\
			&\mathbf{b}_{1}(f,f)(z)=\frac{2}{\pi}\mathrm{Tr}[A]\hat{f}_{x_{0}}(z)+\frac{1}{\pi}|A\bar{z}+2Bz|^{2}+\mathcal{O}(|z|^{3}),\\
			&\mathbf{b}_{2}(f,f)(z)=\mu(f,x_{0})+\mathcal{O}(|z|).
		\end{split}
		\label{eq:5.45a}
	\end{equation}
	
\begin{Definition}
	Associated with the K\"{a}hler form $\Theta$ and $f$ near $x_{0}$, we 
	define a (strictly) positive function on $\C^{n}$ as follows,
	\begin{equation}
		F_{f,x_{0}}(z)=\hat{f}^{2}_{x_{0}}(z)-\frac{1}{2\pi}(\Delta 
		f)(x_{0}) \hat{f}_{x_{0}}(z)+ \frac{1}{\pi}|A\bar{z}+2Bz|^{2} 
		+\mu(f,x_{0}).
		\label{eq:5.48a}
	\end{equation}
	Note that this function does not 
	depend on the choice of the holomorphic coordinate systems 
	centered at $x_{0}$ satisfying \eqref{eq:5.41a}.
	Equivalently, we have for $z\in\C^{n}\simeq (T_{x_{0}}X, 
	J_{x_{0}})$,
	\begin{equation}
		F_{f,x_{0}}(z)=\lim_{p\rightarrow\infty} 
		\big\{p^{2}\mathbf{b}_{0}(f,f)(z/\sqrt{p})+p\mathbf{b}_{1}(f,f)(z/\sqrt{p})+\mathbf{b}_{2}(f,f)(z/\sqrt{p})\big\}.
	\end{equation}
	We also define the following positive quadratic function in 
	$z\in\C^{n}$,
	\begin{equation}
		\widehat{\mathbf{b}}_{1}(z)=\lim_{p\rightarrow\infty} 
		p\mathbf{b}_{1}(f,f)(z/\sqrt{p})=-\frac{1}{2\pi}(\Delta 
		f)(x_{0}) \hat{f}_{x_{0}}(z)+ \frac{1}{\pi}|A\bar{z}+2Bz|^{2}.
		\label{eq:5.50a}
	\end{equation}
\end{Definition}

\begin{Proposition}
	With above notation, set
	\begin{equation}
		\beta_{f,x_{0}}:=\partial\overline{\partial} 
		\widehat{\mathbf{b}}_{1}=\partial\overline{\partial} 
		F_{f,x_{0}}(0)\in 
		\Lambda^{(1,1)}T^{\ast}_{x_{0}}X,
	\end{equation}
	then it is a positive $(1,1)$-form on $\C^{n}$, more precisely,
	\begin{equation}
		\beta_{f,x_{0}}=(dz)^{T}\wedge K_{f,x_{0}}d\bar{z},
	\end{equation}
	where $K_{f,x_{0}}$ is the semipositive definite Hermitian matrix 
	given by
	\begin{equation}
		K_{f,x_{0}}=\frac{2}{\pi}\mathrm{Tr}[A]A+\frac{1}{\pi}(A^{2}+4B\bar{B}).
	\end{equation}
	We have the convergence of $(1,1)$-forms at $x_{0}$ as 
	$p\rightarrow \infty$,
	\begin{equation}
		\frac{1}{p}\partial\overline{\partial}\log 
		T^{2}_{f,p}(x,x)|_{x=x_{0}}\rightarrow 
		\frac{1}{\mu(f,x_{0})}\beta_{f,x_{0}}=\partial\overline{\partial}\log{ F_{f,x_{0}}}(0).
		\label{eq:5.54a}
	\end{equation}
\end{Proposition}
\begin{proof}
	The first part of our proposition follows directly from the 
	formulae \eqref{eq:5.48a} and \eqref{eq:5.50a}. We now prove 
	\eqref{eq:5.54a}.
	
	In the complex coordinate $z$ centered at $x_{0}$, for $|z|<1$, we have
	\begin{equation}
		p^{-n}T^{2}_{f,p}(z,z)=p^{-2}F_{f,x_{0}}(\sqrt{p}z)+\mathcal{O}(p^{-3})+\mathcal{O}(|z|^{5})+p^{-1}\mathcal{O}(|z|^{3})+p^{-2}\mathcal{O}(|z|).
	\end{equation}
	Then as a smooth differential form around $x_{0}$, we have
	\begin{equation}
		\begin{split}
			&\frac{1}{p}\partial\overline{\partial}\log(p^{-n}T^{2}_{f,p}(z,z))\\
			&=\frac{\big(\partial\overline{\partial}F_{f,x_{0}}\big)(\sqrt{p}z)+\mathcal{O}(p^{-1})+p\mathcal{O}(|z|^{3})+\mathcal{O}(|z|)}{F_{f,x_{0}}(\sqrt{p}z)+\mathcal{O}(p^{-1})+p^{2}\mathcal{O}(|z|^{5})+p\mathcal{O}(|z|^{3})+\mathcal{O}(|z|)}\\
			&\quad-\frac{\big(\partial F_{f,x_{0}}\wedge\overline{\partial} 
			F_{f,x_{0}}\big)(\sqrt{p}z)+p^{3}\mathcal{O}(|z|^{7})+p^{2}\mathcal{O}(|z|^{5})+p\mathcal{O}(|z|^{3})+\mathcal{O}(|z|)+\mathcal{O}(p^{-1})}{\big\{F_{f,x_{0}}(\sqrt{p}z)+\mathcal{O}(p^{-1})+p^{2}\mathcal{O}(|z|^{5})+p\mathcal{O}(|z|^{3})+\mathcal{O}(|z|)\big\}^{2}}.
		\end{split}
		\label{eq:5.56a}
\end{equation}
Take $z=0$ in \eqref{eq:5.56a} and then take its limit as 
$p\rightarrow\infty$, we get exactly \eqref{eq:5.54a}.
\end{proof}

\begin{Definition}
	Associated with the vanishing point $x_{0}$ of $f$ as above, for 
	$R>0$, we 
	define the linear function
	\begin{equation}
		\Phi^{R}_{f,x_{0}}:\Lambda^{(n-1,n-1)}_{x_{0}}T^{\ast}X\rightarrow \C
	\end{equation}
	as follows, for $\alpha\in \Lambda^{(n-1,n-1)}_{x_{0}}$, viewed as 
	a constant $(n-1,n-1)$-form on $\C^{n}\simeq (T_{x_{0}}X, 
	J_{x_{0}})$, then
	\begin{equation}
		\Phi^{R}_{f,x_{0}}(\alpha):=\frac{\sqrt{-1}}{2\pi}\int_{B^{\C^{n}}(0,R)}\partial\overline{\partial}\log{ F_{f,x_{0}}(z)}\wedge\alpha.
	\end{equation} 
\end{Definition}
\begin{Remark}
	It is possible to work out more concretely the quantity
	$\Phi^{R}_{f,x_{0}}(\alpha)$ using the formula \eqref{eq:5.48a}, 
	especially if $f$ has a nice shape near $x_{0}$ (for instance, 
	$B=0$). We will give a demonstration in Example \ref{Ex:5.15a}, 
	but we expect that the computations in general would be much more
	complicated, so that we will not try to do it in this paper.
\end{Remark}

\begin{Example}\label{Ex:5.15a}
	Now we assume $f$ near $x_{0}$ is given by \eqref{eq:5.42b} where 
	$B=0$ and
	\begin{equation}
		A=\mathrm{Id}_{n}
	\end{equation}
	Then
	\begin{equation}
		F_{f,x_{0}}(z)=|z|^{4}+\frac{(2n+1)}{\pi}|z|^{2}+\frac{n^{2}}{\pi^{2}}.
	\end{equation}
	Set $\omega_{0}=\sqrt{-1}\sum_{j}dz_{j}\wedge d\bar{z}_{j}$. Then 
	we have
	\begin{equation}
		\begin{split}
			&\sqrt{-1}\partial\overline{\partial}\log{ 
			F_{f,x_{0}}(z)}\wedge\frac{\omega_{0}^{n-1}}{(n-1)!}\\
			&=\pi\Big[\frac{(2n-2)\pi^{3}|z|^{6}+(6n^{2}-n-2)\pi^{2}|z|^{4}+(6n^{3}+2n^{2}-3n-1)\pi|z|^{2}+2n^{4}+n^{3}}{\pi^{4}|z|^{8}+(4n+2)\pi^{3}|z|^{6}+(6n^{2}+4n+1)\pi^{2}|z|^{4}+(4n^{3}+2n^{2})\pi|z|^{2}+n^{4}} \Big]\frac{\omega_{0}^{n}}{n!}.
		\end{split}
	\end{equation}
	In the case of $n=1$,
		\begin{equation}
		\begin{split}
			&\sqrt{-1}\partial\overline{\partial}\log{ 
			F_{f,x_{0}}(z)}\\
			&=\pi\Big[\frac{3\pi^{2}|z|^{4}+4\pi|z|^{2}+3}{\pi^{4}|z|^{8}+6\pi^{3}|z|^{6}+11\pi^{2}|z|^{4}+6\pi|z|^{2}+1}\Big] \omega_{0}.
		\end{split}
	\end{equation}
\end{Example}

\begin{Theorem}\label{thm:5.14a}
	We suppose that \hyperlink{AssumpB}{\bf Assumption B} holds with a nontrivial 
	$f\in \mathcal{Q}_{\geq p_{0}}(X,L;\R_{\geq 0})$. Let $x_{0}$ be 
	a vanishing point of $f$ with $\Delta f(x_{0})<0$. Then for any 
	fixed $R>0$, $\varphi\in 
	\Omega^{(n-1,n-1)}_{0}(X)$, and for all $p\gg 0$,
	\begin{equation}
		\frac{\sqrt{-1}}{2\pi}\int_{B(x_{0},R/\sqrt{p})}\partial\overline{\partial}\log(T^{2}_{f,p}(x,x))\wedge\varphi=p^{-n+1}\Phi^{R}_{f,x_{0}}\big(\varphi({x_{0}})\big)+\mathcal{O}(p^{-n+1/2}).
		\label{eq:5.59a}
	\end{equation}
\end{Theorem}
\begin{proof}
	Note that for $p\gg 0$, then we identify
	\begin{equation}
		B(x_{0},R/\sqrt{p})\simeq B^{\C^{n}}(0,R/\sqrt{p}).
	\end{equation}
	Then for $z\in B^{\C^{n}}(0,R/\sqrt{p})$, $l\in\N$,
	\begin{equation}
		p^{l}\mathcal{O}(|z|^{2l+1})=\mathcal{O}(p^{-1/2}).
	\end{equation}
	Also note for $z\in B^{\C^{n}}(0,R)$,
	\begin{equation}
		\varphi(z/\sqrt{p})=\varphi(x_{0})+\mathcal{O}(p^{-1/2}).
	\end{equation}
	Then \eqref{eq:5.59a} follows from \eqref{eq:5.56a}. This 
	complete our proof.
\end{proof}

As explain in Subsection \ref{s1.4a}, the formula \eqref{eq:5.59a} 
gives the different powers of $p$ in \eqref{eq:5.64a}, which shows the different fluctuations of our random zeros near a 
vanishing point or a nonvanishing point of $f$.

\subsection{Case of real functions with negative values}\label{ss5.5}
In this subsection, we would like to continue the discussion in 
Remark \ref{rk:4.16} and study the equidistribution of random zeros 
for the $\cLL$-holomorphic sections \textit{detected} by a given real function 
$f$ which is not necessary to be nonnegative.

Now we consider the case of complete K\"{a}hler manifold $(X,\Theta)$ 
equipped with a prequantum holomorphic line bundle $(L,h_{L})$. 
Recall that $\mathcal{Q}(X,L^{p};\R)$ is the subspace of 
$\mathcal{Q}(X,L^{p};\C)$ consisting of real valued functions, and 
that
\begin{equation}
	\mathcal{Q}_{\geq p_{0}}(X,L;\R):=\cap_{p\geq 
	p_{0}}\mathcal{Q}(X,L^{p};\R).
\end{equation}

\begin{Definition}
Let $f$ be a real smooth function on $X$, for $x\in X$, we say $f$ is vanishing properly at $x$ up to order $2$ if one of the following 
cases holds:
\begin{itemize}
	\item $f(x)\neq 0$, or
	\item $f(x)=0$, $df(x)\neq 0$, or
	\item $f(x)=0, df(x)=0, \Delta f(x)\neq 0$ with 
$f\Delta f \leq 0$ on an open neighborhood of $x$.
\end{itemize}
For any subset $U\subset X$, we say $f$ is vanishing properly on $U$ 
up to order $2$ if it is so for every point in $U$. Given such a 
function, we also set
\begin{equation}
	\kappa(K):=\max_{x\in K}\mathrm{ord}_{x}(f)\in\{0,1,2\}.
\end{equation}
\end{Definition}

The following proposition is an extension of Proposition 
\ref{prop:5.7a}.
\begin{Proposition}\label{prop:5.17a}
Assume that \hyperlink{AssumpB}{\bf Assumption B} holds with a nontrivial $f\in 
\mathcal{Q}_{\geq p_{0}}(X,L;\R)$. Let $U$ be an open subset of 
	$\supp f$ be such that $f$ vanishes properly on $U$ up to order 
	$2$. Then for any compact subset $K$ of $U$, there exists a 
constant $c_{K}>0$ and $p_{K}\geq p_{0}$ such that for $x\in K$, 
$p\geq p_{K}$,
\begin{equation}
	T^{2}_{f,p}(x,x)\geq c_{K}p^{n-\kappa(K)}.
	\label{eq:5.74a}
\end{equation}

Moreover, $\log f^{2}$ is locally integrable on $U$, and we have 
weak convergence of currents on $U$ as $p\rightarrow \infty$,
\begin{equation}
	\partial\overline{\partial} \log T^{2}_{f,p}(x,x)\rightarrow 
	\partial\overline{\partial} \log f^{2}.
	\label{eq:5.75a}
\end{equation}
Around a point $x$ where $f(x)\neq 0$, the convergence in 
\eqref{eq:5.75a} holds in any local $\mathscr{C}^{l}$-norms.
\end{Proposition}
\begin{proof}
	We start with proving \eqref{eq:5.74a}. For $x_{0}\in U$, if 
	$f(x_{0})\neq 0$, then $f^{2}(x_{0})>0$, \eqref{eq:5.74a} holds 
	near $x_{0}$. 
	If $f(x_{0})=0$, $df(x_{0})\neq 0$, then in a sufficiently small neighborhood of 
	$x_{0}$, there is a constant $c_{x_{0}}>0$ such that 
	have 
	\begin{equation}
		\mathbf{b}_{1}(f,f)=\frac{1}{8\pi}(\mathbf{r}^{X}f-4\Delta 
		f)f+\frac{1}{2\pi}|\partial f|^{2}\geq 
		c_{x_{0}}|df(x_{0})|^{2}_{g^{T^{\ast}X}_{x_{0}}}>0,
	\end{equation}
	so that near $x_{0}$,
	\begin{equation}
	T^{2}_{f,p}(x,x)\geq \frac{1}{2}c_{x_{0}}p^{n-1}.
	\label{eq:5.77a}
\end{equation}
If $\mathrm{ord}_{x_{0}}(f)=2$, we can adapt the proof of Proposition 
\ref{prop:5.7a}. The condition that $\Delta f(x_{0})$ is nonzero with 
$f\Delta f \leq 0$ near $x_{0}$ implies that 
on a small neighborhood of $x_{0}$,
\begin{equation}
	(\mathbf{r}^{X}f-4\Delta 
		f)f\geq 0,\; \mu(f,x_{0})>0.
\end{equation}
Then \eqref{eq:5.74a} still holds near $x_{0}$. The second part of 
our proposition also follows from the analogue arguments in the proof 
of Proposition \ref{prop:5.7a}.
\end{proof}

For $f\in \mathcal{Q}_{\geq p_{0}}(X,L;\R)$, the operator $T_{f,p}$ 
might not be injective, so that, in Remark \ref{rk:4.16}, we 
introduce a closed subspace 
$H^{0}_{(2)}(X,L^{p},f)=(\ker T_{f,p})^{\perp}$ of 
$H^{0}_{(2)}(X,L^{p})$ and the Gaussian probability measure 
$\P^{\sharp}_{f,p}$ on it. Consider the following random sections
\begin{equation}
	\begin{split}
			&(\mathbf{S}^{\sharp}_{f,p})_{p\geq p_{0}}\in \Pi_{p\geq 
	p_{0}}\big(H^{0}_{(2)}(X,L^{p},f), \; \P^{\sharp}_{f,p}\big),\\
	&(\mathbf{S}_{f,p})_{p\geq p_{0}}\in \Pi_{p\geq 
	p_{0}}\big(H^{0}_{(2)}(X,L^{p}), \;
	\P^{0}_{f,p}\otimes\P^{\sharp}_{f,p}\big).
	\end{split}
\end{equation}

From \eqref{eq:4.47a} and by Proposition \eqref{prop:5.17a}, we get
\begin{Theorem}
	We suppose that \hyperlink{AssumpB}{\bf Assumption B} holds with a nontrivial $f\in \mathcal{Q}_{\geq p_{0}}(X,L;\R)$. Let $U$ be an open subset of 
	$\supp f$ be such that $f$ vanishes properly on $U$ up to order $2$. Then as $p\rightarrow \infty$,
\begin{itemize}
	\item we have the convergence of $(1,1)$-currents on $U$
	\begin{equation}
		\begin{split}
					&\frac{1}{p}\E^{\P^{\sharp}_{f,p}}[[\Div(\mathbf{S}^{\sharp}_{f,p})]|_{U}]\rightarrow c_{1}(L,h_{L})|_{U}\\
					&\frac{1}{p}\E^{\P_{f,p}}[[\Div(\mathbf{S}_{f,p})]|_{U}]\rightarrow c_{1}(L,h_{L})|_{U}.
		\end{split}
		\label{eq:5.79a}
	\end{equation}
	\item for any $\varphi\in 
\Omega^{(n-1,n-1)}_{0}(\overline{U})$, we have
\begin{equation}
	\begin{split}
		&\P\left(\lim_{p\rightarrow 
	\infty}\frac{1}{p}\langle[\Div(\mathbf{S}^{\sharp}_{f,p})],\varphi\rangle=\langle c_{1}(L,h_{L}),\varphi\rangle\right)\\
	&=\P\left(\lim_{p\rightarrow 
	\infty}\frac{1}{p}\langle[\Div(\mathbf{S}_{f,p})],\varphi\rangle=\langle c_{1}(L,h_{L}),\varphi\rangle\right)=1.
	\end{split}
\end{equation}
\end{itemize}
\end{Theorem}

\begin{Remark}
	If $X$ is compact, then $H^{0}_{(2)}(X,L^{p})=H^{0}(X,L^{p}), p\in\N$, are finite 
dimensional, and we can take $f$ to be any real smooth 
	function vanishing properly up to order $2$ in the above theorem. 
	If $\kappa(X)\leq 1$, then the first convergence in 
	\eqref{eq:5.79a} is already proved by Ancona-Le Floch 
	\cite{AF:2022}. As mentioned in Subsection \ref{ss5.3ab}, they 
	also studied the fluctuations of the random zeros near a vanishing 
	point of $f$ with order $1$, since the computations is local, 
	then it is also applicable in our noncompact setting.
\end{Remark}

\end{document}